\documentclass[11pt]{amsart}

\usepackage[utf8]{inputenc}
\usepackage[T1]{fontenc}
\usepackage[english]{babel}
\usepackage{amsmath}
\usepackage{amsthm}
\usepackage{amssymb}
\usepackage{mathrsfs}
\usepackage{dsfont}
\usepackage{array}
\usepackage{graphicx} 
\usepackage[left=3cm,right=3cm,top=3cm,bottom=3cm]{geometry}
\usepackage{hyperref}
\usepackage{multirow}
\usepackage{tikz-cd}
\usetikzlibrary{3d}
\usetikzlibrary{shapes}
\usetikzlibrary[patterns]
\usetikzlibrary{decorations.pathreplacing}
\usepackage{caption}
\usepackage{subcaption}
\usepackage{adjustbox}
 \usepackage{comment}
\usepackage{enumitem}
\setitemize{label=$\circ$}

\usepackage{eucal} 

\usepackage{mlmodern}

\newcommand{\bcd}{\begin{center}\begin{tikzcd}}
\newcommand{\ecd}{\end{tikzcd}\end{center}}

\newcommand{\vir}[1]{[#1]^\mathrm{vir}}
\newcommand{\red}[1]{[#1]^\mathrm{red}}
\newcommand{\fvir}[2]{ [\![ #1 ]\!]^{#2} }
\newcommand{\ev}{\mathrm{ev}}
\newcommand{\pt}{\mathrm{pt}}

\newcommand{\m}[1]{\mfk\left[ #1 \right]}
\renewcommand{\d}[1]{\dfk\left[ #1 \right]}

\newcommand{\Pic}{\operatorname{Pic}}
\newcommand{\Alb}{\operatorname{Alb}}
\newcommand{\Hom}{\operatorname{Hom}}

\newcommand{\DR}{\operatorname{DR}}
\newcommand{\LogPic}{\operatorname{LogPic}}
\newcommand{\TroPic}{\operatorname{TroPic}}

\newcommand{\Gmlog}{\mathbb{G}_{\mathrm{log}}}
\newcommand{\Gmtrop}{\mathbb{G}_{\mathrm{trop}}}

\newcommand{\lcm}{\operatorname{lcm}}

\newcommand{\gp}{\mathrm{gp}}
\newcommand{\rub}{\mathrm{rub}}
\newcommand{\Ev}{\operatorname{Ev}}

\newcommand{\prim}{\mathrm{prim}}
\newcommand{\Prim}{\operatorname{Prim}}

\newcommand{\PP}{\mathbb{P}}
\newcommand{\NN}{\mathbb{N}}
\newcommand{\ZZ}{\mathbb{Z}}
\newcommand{\RR}{\mathbb{R}}
\newcommand{\CC}{\mathbb{C}}

\newcommand{\KK}{\mathbb{K}}
\newcommand{\QQ}{\mathbb{Q}}

\newcommand{\BBB}{\mathscr{B}}

\newcommand{\DDD}{\mathscr{D}}

\newcommand{\FFF}{\mathscr{F}}

\newcommand{\PPP}{\mathscr{P}}
\newcommand{\RRR}{\mathscr{R}}

\newcommand{\TTT}{\mathscr{T}}

\newcommand{\bfa}{\mathbf{a}}
\newcommand{\bfw}{ {\boldsymbol w} }
\newcommand{\bfm}{ {\boldsymbol m} }

\newcommand{\tK}{\widetilde{K}}

\newcommand{\ttD}{\mathtt{D}}

\newcommand{\bsigma}{ \boldsymbol{\sigma} }

\newcommand{\bfF}{ {\boldsymbol F} }
\newcommand{\bfG}{ {\boldsymbol G} }

\newcommand{\bfN}{ {\boldsymbol N} }
\newcommand{\bfM}{ {\boldsymbol M} }

\newcommand{\GW}{\mathtt{GW}}
\newcommand{\Dec}{\operatorname{Dec}}
\newcommand{\ft}{\mathrm{ft}}

\newcommand{\Dfk}{\mathfrak{D}}
\newcommand{\tDfk}{ {\widetilde{\mathfrak{D}}} }

\newcommand{\Ffk}{\mathfrak{F}}
\newcommand{\tFfk}{ {\widetilde{\mathfrak{F}}} }

\newcommand{\mfk}{\mathfrak{m}}

\newcommand{\dfk}{\mathfrak{d}}

\newcommand{\sfV}{\mathsf{V}}

\newcommand{\sfH}{\mathsf{H}}

\newcommand{\sfE}{\mathsf{E}}

\newcommand{\A}{\mathcal{A}}

\renewcommand{\P}{\mathcal{P}}
\renewcommand{\L}{\mathcal{L}}

\newcommand{\M}{\mathcal{M}}

\newcommand{\G}{\mathcal{G}}

\newcommand{\X}{\mathcal{X}}
\newcommand{\T}{\mathcal{T}}

\newcommand{\C}{\mathcal{C}}

\renewcommand{\O}{\mathcal{O}}
\newcommand{\F}{\mathcal{F}}
\newcommand{\Y}{\mathcal{Y}}
\newcommand{\U}{\mathcal{U}}

\renewcommand{\geq}{\geqslant}

\renewcommand{\bar}{\overline}

\newcommand{\gen}[1]{\langle #1 \rangle}

\newtheorem{theo}{Theorem}[section]
\newtheorem*{theom}{Theorem}
\newtheorem{prop}[theo]{Proposition}
\newtheorem{coro}[theo]{Corollary}
\newtheorem{lem}[theo]{Lemma}
\newtheorem{Notation}[theo]{Notation}

\newtheorem{conj}{Conjecture}[section]

\theoremstyle{definition}
\newtheorem{defi}[theo]{Definition}

\theoremstyle{remark}
\newtheorem{remark}[theo]{Remark}
\newenvironment{rem}[1]{
    \begin{remark}#1}{
    \xqed{\blacklozenge}\end{remark}
}

\theoremstyle{remark}
\newtheorem{mainexample}[theo]{Main Example}

\theoremstyle{remark}
\newtheorem{example}[theo]{Example}
\newenvironment{expl}[1]{
    \begin{example}#1}{
    \xqed{\lozenge}\end{example}
}

\newcommand{\xqed}[1]{
    \leavevmode\unskip\penalty9999 \hbox{}\nobreak\hfill
    \quad\hbox{\ensuremath{#1}}}

\def\floor (#1) at (#2,#3); {
    \node[draw,ellipse, minimum width=1cm, minimum height = 0.6 cm] (#1) at (#2,#3) {$\bullet$} ;
}

\def\ufloor (#1) at (#2,#3) (#4); {
    \node[draw,ellipse, minimum width=1cm, minimum height = 0.6 cm] (#1) at (#2,#3) {\scriptsize #4} ;
}

\def\marked (#1) to (#2) pos=#3 in=#4 out=#5; {
   \draw (#1) to[out=#5,in=#4] node[pos=#3] {$\bullet$} (#2) ;
}

\def\leftedge (#1) to (#2) pos=#3 in=#4 out=#5 w=#6; {
   \draw (#1) to[out=#5,in=#4] node[midway,left] {$#6$} (#2) ;
}
\def\leftmarked (#1) to (#2) pos=#3 in=#4 out=#5 w=#6; {
   \draw (#1) to[out=#5,in=#4] node[pos=#3] {$\bullet$} node[midway,left] {$#6$} (#2) ;
}

\def\wlmarked (#1) to (#2) pos=#3 in=#4 out=#5 w=#6; {
   \draw (#1) to[out=#5,in=#4] node[pos=#3] {$\bullet$} node[midway,left] {$#6$} (#2) ;
}

\def\rightedge (#1) to (#2) pos=#3 in=#4 out=#5 w=#6; {
   \draw (#1) to[out=#5,in=#4] node[midway,right] {$#6$} (#2) ;
}
\def\rightmarked (#1) to (#2) pos=#3 in=#4 out=#5 w=#6; {
   \draw (#1) to[out=#5,in=#4] node[pos=#3] {$\bullet$} node[midway,right] {$#6$} (#2) ;
}

\def\doublemarked (#1) to (#2) pos=#3 in=#4 out=#5; {
   \draw (#1) to[out=#5,in=#4] node[pos=#3] {$\bullet$} (#2) ;
}

\usepackage{todonotes}

\makeatletter
\def\l@subsection{\@tocline{2}{0pt}{2.5pc}{5pc}{}}
\makeatother

\makeatletter
\renewcommand{\l@section}{\@tocline{1}{0pt}{10pt}{1pc}{\bfseries}}
\makeatother
\title[MCF via correlated invariants]{Multiple cover formulas for abelian surfaces via correlated invariants} 
\author{Thomas Blomme, Francesca Carocci}

\address{Institut de Math\'ematiques de Toulouse, Universit\'e de Toulouse, route de Narbonne, Toulouse, France}
\email{thomas.blomme@math.univ-toulouse.fr}

\address{Università di Roma Tor Vergata, Via della Ricerca Scientifica 1, Roma 00133, Italy}
\email{carocci@mat.uniroma2.it}


\begin{document}

\maketitle

\begin{abstract}
We prove the multiple cover formula conjecture for Gromov--Witten invariants of abelian surfaces. This was conjectured by Oberdieck \cite{oberdieck2022gromov}, heavily inspired by the multiple cover formula conjecture for invariants of K3 surfaces \cite{oberdieck2016curve}.
The proof uses: a new reduced degeneration formula expressing the abelian surface invariants in terms of the correlated Gromov--Witten invariants introduced by the authors; a correlated refinement for Pixton's  double ramification cycle formula proved in previous work; Urundolil-Kumaran--Ranganathan log to DR comparison Theorem \cite{ranganathan2024logarithmic}; Maulik-Ranganathan theory of exotic insertion \cite{maulik2025gromov} to recover vanishing cohomology.
\end{abstract}

\setcounter{tocdepth}{2}
\tableofcontents




\section{Introduction}

\subsection{Reduced Gromov--Witten invariants and multiple cover formulas }

\subsubsection{Abelian surfaces}

Complex tori are quotients of $\CC^2$ by a rank $4$ lattice $\Lambda\subset\CC^2$. In particular, they inherit from $\CC^2$ a non-degenerate holomorphic volume form. The existence of the latter forces the vanishing of the usual virtual class on the moduli of stable maps and thus of the usual Gromov-Witten invariants \cite{kiem2013localizing}. Abelian surfaces are those complex tori which are projective.

Then the vanishing of the usual Gromov-Witten invariants may also be seen as a consequence of the following fact: for a generic choice of lattice $\Lambda$, the complex torus $\CC^2/\Lambda$ is not a projective variety and does not contain any complex curve. 

To get interesting numbers, one has to consider a \textit{reduced} virtual fundamental class on the moduli of stable maps, as already done in \cite{bryan1999generating} or \cite{bryan2018curve}, and also similarly for GW-invariants of K3 surfaces \cite{bryan2000enumerative},\cite{maulik2010curves}.

\subsubsection{Gromov-Witten theory and multiple cover formulas}

Let $A=\CC^2/\Lambda$ be an abelian surface with $\beta\in H_2(A,\ZZ)$ a realizable curve class, Poincar\'e dual to a multiple of an ample \textit{primitive} class $c_1(\L)$. Following \cite{bryan2018curve}, the moduli space of stable maps $\M_{g,n}(A,\beta)$ is endowed with a reduced virtual class $\red{\M_{g,n}(A,\beta)}$ of dimension $g+n$. 
Reduced Gromov-Witten cycles are obtained by capping the latter with cohomology classes pulled back along the evaluation map, and pushing forward via the forgetful map to the moduli of stable curves:
$$\ev\colon\M_{g,n}(A,\beta)\to A^n,\quad \mathrm{ft}\colon\M_{g,n}(A,\beta)\to\overline{\M}_{g,n}.$$
For $\gamma\in H^*(A^n,\QQ)$, we set
$$\gen{\gamma}_{g,\beta} = \ft_*(\red{\M_{g,n}(A,\beta)}\cap\ev^*\gamma).$$
To get numbers, called GW-invariants, it suffices to intersect with a class $\alpha\in H^*(\overline{\M}_{g,n},\QQ)$ of complementary dimension:
$$\gen{\alpha;\gamma}_{g,\beta} = \int_{\red{\M_{g,n}(A,\beta)}}\mathrm{ft}^*\alpha\cup\ev^*\gamma.$$

It is well-known that the deformation invariance property of these numbers implies that the invariants depend on $\beta$ only through its self-intersection $\beta^2$ and its divisibility $\ell(\beta)$.

\medskip

The divisibility $1$ case, also known as \textit{primitive case}, is easier to handle and several computations and conjectures are known in this case \cite{bryan1999generating,bryan2018curve}.
Concerning divisible classes, in \cite[Appendix B]{oberdieck2022gromov}, G. Oberdieck conjectures a very simple multiple cover formula expressing the reduced invariants for non-primitive classes in terms of the invariants for primitive classes.

\medskip

To state the conjecture, we need some elementary results on automorphisms of $H^*(A,\CC)$. Given $g,n$ and $\beta$, and eventually $\alpha\in H^*(\overline{\M}_{g,n},\QQ)$,  GW-cycles (resp. GW-invariants) may be seen as a linear function of $\gamma\in H^*(A,\CC)^{\otimes n}$, the $n$th-tensor power of $H^*(A,\CC)$, which is itself the exterior algebra of $H^1(A,\CC)$. 

Furthermore, any curve class $\beta\in H_2(A,\ZZ)\subset H_2(A,\CC)$ can be interpreted as a symplectic form over $H^1(A,\CC)$. By \cite[Chapter 2.6]{griffiths2014principles}, the moduli of abelian varieties with given polarization of fixed type, is the quotient of the Siegel upper half-plane by the Siegel modular group $\operatorname{Sp}_\ZZ(\beta)$ of integral automorphisms of $H^1(A,\ZZ)$ preserving the symplectic form $\beta$. 

Deformation invariance of the reduced virtual class then implies that GW-cycles (and GW-invariants) are invariant by the action of $\operatorname{Sp}_\ZZ(\beta)$ on the argument $\gamma$. The group $\operatorname{Sp}_\ZZ(\beta)$ being Zariski dense in the symplectic group $\operatorname{Sp}_\CC(\beta)$ preserving $\beta$ (with $\CC$-coefficient), we deduce the $\operatorname{Sp}_\CC(\beta)$-invariance. This group is a subgroup of the bigger group $SL(H^1(A,\CC))$ of automorphisms preserving the point class $\pt\in H^4(A,\CC)=\Lambda^4 H^1(A,\CC)$.

The group $SL(H^1(A,\CC))$ acts by base change on the space of skew-symmetric form $H_2(A,\CC)$, and we denote the action of $\varphi\in SL(H^1(A,\CC))$ by $\varphi^*$.

If $\beta,\beta'\in H_2(A,\ZZ)$ are two curve classes, which are in particular skew-symmetric forms, a necessary and sufficient condition to find an element $\varphi\in SL(H^1(A,\ZZ))$ such that $\beta'=\varphi^*\beta$ is that they share the same self-intersection and divisibility. Working instead with $\CC$-coefficient, we do not have a divisibility anymore and it is then possible to find $\varphi\in SL(H^1(A,\CC))$ such that $\varphi^*\beta=\beta'$ if and only if $\beta^2=(\beta')^2$. We can now state the conjecture. In the remainder of the paper, we abbreviate ``multiple cover formula'' by MCF.

\begin{conj}\label{conj:oberdieck}
    If $\gamma\in H^*(A^n,\QQ)$ and $\alpha\in H^*(\overline{\M}_{g,n},\QQ)$. For each $k|\beta$, we denote by $\varphi_k\in SL(H^1(A,\CC))$ an isomorphism such that $\varphi^*_k(\beta/k)$ is an integral and primitive curve class (with the same self-intersection) in $A$. We then have
    $$\gen{\alpha;\gamma}_{g,\beta} = \sum_{k|\beta} k^{3g-3+n-\deg\alpha}\gen{\alpha;\varphi_k(\gamma)}_{g,\varphi_k^*(\beta/k)},$$
\end{conj}

It is possible to check that if $\beta'=\varphi^*\beta$, then $\operatorname{Sp}_\CC(\varphi^*\beta)=\varphi\operatorname{Sp}_\CC(\beta)\varphi^{-1}$. In particular, each of the terms on the right-hand side is indeed invariant by the action of $\operatorname{Sp}_\CC(\beta)$. Furthermore, it proves that $\gen{\alpha;\varphi_k(-)}_{g,\varphi_k^*(\beta/k)}$ does not depend on the precise choice of $\varphi_k$.

\medskip

Oberdieck's MCF for abelian surfaces was conjectured in \cite{oberdieck2022gromov} following the MCF for reduced GW-invariants in the K3 case, which has a slightly more ancient history. For general insertions, the formula in the K3 setting was conjectured by Oberdieck and Pandharipande \cite{oberdieck2016curve}, but for  special cases, namely for genus $0$ curves, and for genus $g$ curves with a $\lambda_g$ insertion, the  MCF for K3 surfaces goes back to Aspinwall-Morrison and Yau-Zaslow in the 90's \cite{yau1996bps}, and  to Katz-Klemm-Vafa \cite{katz1999m} respectively. 
\medskip
In the table below we recap what is known about GW-invariants (primitive and non) for both the K3 and the abelian surface case. For divisible classes, the computations rely on a MCF, reducing the computation to the primitive case. We therefore highlight which instances of the multiple cover formula conjecture have been proven for either geometry. This in particular explains how the current work fits in the landscape of multiple cover formulas conjectures and results for both K3 and abelian surfaces. Essentially, in the K3 case, only a handful of specific cases are known, beside the divisibility 2 case \cite{bae2021curves} which is valid for general insertions. This leaves most cases of the conjecture open. 

\medskip

In the abelian setting, the current paper solves the conjecture.

\medskip

\renewcommand{\arraystretch}{2}
\begin{center}\label{table}
\scalebox{0.7}{
\begin{tabular}{|c||l|l||l|l||l|l||}
  \hline
  \multirow{3}{*}{} 
      & \multicolumn{2}{c||}{{\bf{g=0}}} 
          & \multicolumn{2}{c||}{{\bf{ genus g , only a $\lambda_g$  }}} 
          &  \multicolumn{2}{c||}{{\bf{ genus g,  other insertions}}}\\   \hline           \cline{3-7}
  &$\rm{div}(\beta)=1$ & $\rm{div}(\beta)>1$ & $\rm{div}(\beta)=1$ & $\rm{div}(\beta)>1$ & $\rm{div}(\beta)=1$& $\rm{div}(\beta)>1$ \\  \hline
{\bf{K3}} &\cite{bryan2000enumerative} &  $\rm{div}(\beta)=2,$ \cite{lee2005yau} & \cite{maulik2010curves} & \cite{pandharipande2016katz} & points,\cite{bryan2000enumerative} &  $\rm{div}(\beta)=2$, \cite{bae2021curves} \\   
& &  any, \cite{klemm2010noether,pandharipande2016katz} & & &\cite{maulik2010curves} & any div, \;\;   \begin{color}{red}\huge ?\end{color} \\ \hline
  $\mathbf{A=\CC^2/\Lambda}$ & $\empty$ & $\empty$ & \cite{bryan2018curve} & $g=2$, \cite{bryan2018curve}  & points,\cite{bryan1999generating} &  points, \cite{blomme2022abelian3,blomme2025short}\\ 
  & & & & any, Theorem~\ref{theo-MCF-abelian-lambda} & \cite{bryan2018curve}&   full conjecture:  Theorem~\ref{theo-MCF-abelian-general}\\\hline
\end{tabular}}

\end{center}


\subsubsection{Proof techniques}

It is somewhat expected that in order to uncover multiple cover phenomena one should exploit correspondence results to move the computation to the sheaf theory side, as the experience suggests that PT invariants and DT invariants for moduli of one dimensional sheaves are much closer to BPS states than stable maps (we refer for example to \cite{pandharipande2010stable,maulik2018gopakumar}). 
This is precisely the approach of \cite{pandharipande2016katz}, and 
at the time of writing this paper, we were aware that G. Oberdieck and R. Pandharipande were also working to prove the multiple cover formula by exploting a GW/PT correspondence. 

Indeed, in their recent paper \cite{oberdieck2026multiple} 
they propose a proof for the descendent multiple cover formulas for both abelian and K3 surfaces, conditional to the proof of the  conjectural families GW/PT
correspondence for semipositive relative 3-folds with primary insertions.

\medskip

The first author's proof of the conjecture for point insertions \cite{blomme2022abelian3,blomme2025short} via tropical methods (relying on Nishinou's correspondence theorem \cite{nishinou2020realization} in abelian surfaces) suggests that the MCF behaves appropriately with respect to degenerations, hinting towards a more general MCF for refined counts in log-CY surfaces. To argue in favor of this conjecture, the present paper also provides a MCF for the log-CY $(E\times\PP^1|E\times\{0,\infty\})$, where $E$ is an elliptic curve. Furthermore, \cite{BCUK2026corrDDR} provides a MCF for the logGW-invariants of a toric surface relative to its full toric boundary. Another clue that MCF phenomena for log-CY surfaces should be more general is provided by the proof of Takahashi conjecture \cite{bousseau2019takahashi}, giving an MCF for a refined count of curves maximally tangent at a one point to a cubic in the plane.
Once again, Bousseau's argument requires moving the computation on moduli of one dimensional sheaves.

\medskip

In the present paper, we introduce a new technique  to study multiple cover formula phenomena \emph{staying within the realm of Gromov-Witten theory}. We show that it is effective by proving  the  MCF for abelian surfaces, as well as MCF  for the log-CY $E\times\PP^1$ relative to its $0$ and $\infty$-sections, where $E$ is an elliptic curve.


\medskip

The present paper may be seen as a milestone in some ongoing project tackling the computation of reduced GW-invariants of surfaces. On the way to prove Conjecture \ref{conj:oberdieck}, the authors have been introducing a distilled version of relative GW invariants, called \textit{correlated GW-invariants} \cite{blomme2024correlated}, which, as advertised in \textit{loc.cit.} and proven here, naturally appear in the degeneration formula for the reduced class $\red{\M_{g,n}(A,\beta)}$.

\medskip

We refer the reader to Section~\ref{sec-future-directions} below for a more detailed account of some of the applications of the ideas introduced here we envision.

\subsection{Correlated GW invariants and reduced decomposition formula}

\subsubsection{Correlation}\label{sec:correlation}

Let $E$ be an elliptic curve. We consider the moduli space of stable maps to $Y=E\times\PP^1$ relative to its \textit{boundary divisor} $D=D^+ + D^-=E\times\{0\}+E\times\{\infty\}$. Let $\bfw=(w_1,\dots,w_n)$ be a $n$-tuple of non-zero integers with $\sum w_i=0$. For a fixed $a>0$, we consider the homology class $\beta_{a,\bfw}=a[E]+(\sum w_i)^+[\PP^1]\in H_2(Y,\ZZ)$.
Let $\M_{g,m}(Y|D,a,\bfw)$, denote the moduli space of relative stable maps with \textit{tangency profile} $\bfw=(w_1,\dots,w_n)$ of class $\beta_{a,\bfw}.$  This is endowed with its natural evaluation map
$$\ev\colon\M_{g,m}(Y|D,a,\bfw)\to E^n\times (E\times \PP^1)^m,$$
where the evaluation at the boundary marked point is defined by post-composition with the projection to $E.$

\medskip

Given $f\colon (C,p_j,q_i)_{\substack{1\leqslant i\leqslant n \\ 1\leqslant j\leqslant m}}\to E\times\PP^1$  a relative stable map, and  denoting by $x_i=\ev_{q_i}(f)\in E$ the image of  a boundary evaluation, we have the relation
\begin{equation}\label{eq:correlation}
    \sum_{i=1}^n w_ix_i = 0\in E.
\end{equation}
Let $|\bfw|=\gcd(w_i)_{1\leqslant i\leqslant n}$ be the gcd of the intersection orders.  Then, we can consider the morphism $\kappa$ from our moduli space to $E[|\bfw|],$ the subgroup of $|\bfw|$-torsion elements in $E$, mapping a relative stable map to
$$ \sum_{i=1}^n \frac{w_i}{|\bfw|}x_i \in E[|\bfw|].$$
We denote by $\theta=\sum_{i=1}^n \frac{w_i}{|\bfw|}x_i$ the image, and call it the \textit{correlator} of $f$, as it contains an information on the position of the boundary points, relative to each other. We will refer to $\kappa$ as the \emph{correlator function}. We refer to \cite{blomme2024correlated} for further details.

Since the target $E[|\bfw|]$ is discrete, the correlator function splits $\M_{g,m}(Y|D,a,\bfw)$ as a union of open and closed substacks $\M_{g,m}^\theta(Y|D,a,\bfw)$ indexed by the correlators $\theta$. We can thus define a refinement of the GW-invariants, called \textit{correlated invariants}, obtained restricting the virtual class to a given component: for $\gamma\in H^*(E^n\times(E\times\PP^1)^m,\QQ)$ and $\alpha\in H^*(\bar{\M}_{g,n+m},\QQ)$, we set
$$\gen{\gamma}_{g,a,\bfw}^{\theta} = \ft_*(\vir{\M_{g,m}^\theta(Y|D,a,\bfw))}\cap\ev^*\gamma) \; \text{and} \; \gen{\alpha;\gamma}_{g,a,\bfw}^{\theta} = \int_{\vir{\M_{g,m}^\theta(Y|D,a,\bfw))}} \mathrm{ft}^*\alpha\cup\ev^*\gamma.$$

\medskip

The computations performed in \cite{blomme2024correlated} suggest that to handle the correlated GW-invariants, one should consider classes with coefficients in the group algebra of $E$, and support given by the correlators, i.e. $|\bfw|$-torsion elements. We thus define the full correlated virtual class
$$\fvir{\M_{g,m}(Y|D,a,\bfw)}{} = \sum_{|\bfw|\theta\equiv 0} \vir{\M_{g,m}^\theta(Y|D,a,\bfw)} \cdot (\theta) \in H^*(\M,\QQ)\otimes\QQ[E].$$
Integrating over the latter, correlated GW-invariants are naturally elements of $\QQ[E]$:
$$\gen{\gen{\alpha;\gamma}}_{g,a,\bfw} = \sum_{|\bfw|\theta\equiv 0} \gen{\alpha;\gamma}^{\theta}_{g,a,\bfw}\cdot (\theta).$$

\subsubsection{Reduced decomposition formula}
Let $a,d\geqslant 1$ be integers, $E,F$ be elliptic curves, and $u\in E[d]$ a $d$-torsion element in $E$. We consider one parameter semi-stable degenerations $A_t(u)\to \Delta$ of families of elliptically fibered surfaces $E\hookrightarrow A_t(u)\xrightarrow{p}F$. 

The parameter $u$ is called the \textit{monodromy} and determines the translation when going around a specified loop in $F$. For instance, the choice $u=0$ corresponds to the trivial product $E\times F$. The  abelian surface $A_t(u)$ comes with a natural polarization $\beta(u)$ fully determined by its intersection with a section $a$, and intersection $b$ with a fiber. We refer to $d$ as the ``norm'' of the curve class and write it $|\beta(u)|$ to highlight the dependence on the curve class. We refer to Section \ref{sec-construction-abelian-families} for details of the construction.

The self-intersection of $\beta(u)$ is equal to $2ad$ and does not depend on the choice of $u$, while the divisibility $\ell(\beta(u))=\gcd(a,d,\frac{d}{\mathrm{ord}(u)})$ depends on $u\in E[|\beta|]$. For instance, the divisibility is $\gcd(a,d)$ when $u=0$, and $1$ if $u$ is chosen to be a primitive $d$-torsion element (i.e. of order $d$). These are explicitly constructed in Section~\ref{sec-construction-abelian-families}.
The central fiber $A_0(u)$ is obtained gluing of $N$ copies of $E\times\PP^1$ along their boundary divisors, where the first and last copies are glued with a $u$-shift. 

\medskip

The degeneration formula allows one to write the reduced class $\red{\M_{g,n}(A_0(u),\beta(u))}$ as a sum of classes of (virtual) components $\M_{\Gamma}$ indexed by certain decorated graphs $\Gamma$, which are the possible (generic) dual graphs of curves mapped to the central fiber, where:
vertices $v\in V(\Gamma)$ correspond to irreducible components, and edges $e=(h,h')\in E(\Gamma)$ to nodes mapped to the singular locus of $A_0(u)$. Here $h,h'$ denote the two half edges associated with $e,$ corresponding to the pre-images of the node in the normalization.

\medskip

The main step when proving a degeneration formula is to prove that the virtual class of each $\M_{\Gamma}$ is obtained a follows:
\begin{itemize}[leftmargin=0.6cm]
    \item to each vertex $v\in\sfV(\Gamma)$ is associated a moduli space $\M_v$ of log stable maps to $(Y_v|D_v)$, where $Y_v$ is one of the components of $A_0(u)$; consider the product of virtual classes $\vir{\M_v}$,
    \item take the cap product with the Poincar\'e dual class of diagonal $\Delta\colon E^{\sfE(\Gamma)}\to E^{\sfH(\Gamma)},$ which geometrically means imposing that evaluation of boundary points corresponding to a common node match. Here, we denoted by $\sfH(\Gamma)$ (resp., $\sfE(\Gamma)$) the set half-edges (resp., edges) of $\Gamma$.
\end{itemize}


\medskip

Equation (\ref{eq:correlation}) tells us that the intersection with the diagonal is not transverse. An easy particular case when we have a unique component, i.e. a curve in $E\times\PP^1$ with a symmetric tangency profile $\bfw$, meaning there is an involution $\tau$ of $[\![1;n]\!]$ with $w_{\tau(i)}=-w_i$. The diagonal constraints forces $x_{\tau(i)}=x_i$ for each $i$ (with a translation in case of non-trivial monodromy). Meanwhile, Equation (\ref{eq:correlation}) takes the following form:
$$\sum_{\{i,\tau(i)\}}w_i(x_i-x_{\tau(i)})=0\in E.$$
In particular, the constraints are not linearly independent: if all gluings except one are satisfied, the last one is satisfied \emph{up to torsion}. There are, therefore, two problems that arise:
\begin{itemize}[leftmargin=0.6cm]
    \item One has to take care of the transversality issue, correcting the class of the diagonal to a \textit{reduced} diagonal; this amounts to consider the embedding of $E^{\sfE(\Gamma)}\xrightarrow{\Delta_u^{\mathrm{red}}}{\bf{H}}_{\Gamma}^u\subseteq  E^{\sfH(\Gamma)}$ of the diagonal into a smaller target. We refer to Section~\ref{sec:evhyperplane} for the definition and the geometric properties.
    \item One needs to handle the bookkeeping of the components that actually glue together and not only up to torsion, which is where correlated invariants are involved. This takes the form of a compatibility relation between vertex correlators $\underline{\theta}=(\theta_v)_{v\in \sfV(\Gamma)}$ and the monodromy $u$.
\end{itemize}

In the end, the reduced decomposition formula takes the following shape (See Theorem~\ref{theo:degene-formula-reduced}):
\[\red{\M_0^u}=\sum_{\widetilde{\Gamma}}\frac{l_{\Gamma}}{|\sfE(\Gamma)|!}\mu_{\widetilde{\Gamma},*}\varphi_{\Gamma}^*\left (\sum_{\underline{\theta}} \Delta_u^{\mathrm{red},!}(\prod_v \vir{\M_v^{\theta_v}}) \right)\]
    where the second sum is over the $\underline{\theta}=(\theta_v)$ satisfying the above mentioned compatibility condition.

\subsection{Strategy}

In a nutshell, we provide not one but three MCF, in three different settings. In order of motivation, these are: abelian surfaces, the log-CY $E\times\PP^1$, the correlated \textit{double ramification cycle}. Through an adaptation of \cite{ranganathan2024logarithmic}, the MCF for the correlated DR-cycle implies the MCF for $E\times\PP^1$, which in turn implies the MCF for abelian surfaces via the reduced decomposition formula. The technical core of the paper could thus be seen as the MCF for the correlated DR-cycle, as the other follow through already known geometric techniques. Meanwhile, the MCF for the correlated DR-cycle stems from combinatorics of coverings of an elliptic curve, i.e. enumeration of subgroups of $(\ZZ/n\ZZ)^2$.

We now present with more care the key steps of the paper.

\subsubsection{Strategy I: degeneration formula and group algebra}

One of the key insights of our approach is that, in analogy to what one experiences with the correlated invariants, the reduced decomposition formula takes a nice form when considered as an element in the group algebra $\QQ[E]$, where the elements in the support, playing the role of correlators, are the various possible monodromies $u$ of the families of abelian surfaces mentioned above.

Let $B=(d,a)$ for which Section \ref{sec-construction-abelian-families} constructs families of elliptically fibered abelian surfaces with polarization $\beta(u)$ indexed by $u\in E[d]$. For $u=0$, $A(0)$ is merely the product of two elliptic curves $E\times F$. Moreover, the construction also provides isomorphisms $\varphi_u$ mapping $\beta(0)=a[E]+d[F]$ to $\beta(u)$. Given $\gamma\in H^*((E\times F)^n,\QQ)$, we can thus consider the following element in the group algebra $\G=\QQ[E]$ with cycle coefficients: 
$$\gen{\gen{\gamma}}_{g,B} = \sum_{u\in E[d]} \ft_*\left(\red{\M_{g,n}(A(u),\beta(u))}\cap\ev^*(\varphi_u(\gamma)) \right) \in H^*(\bar{\M}_{g,n},\QQ)\otimes\QQ[E].$$
Intersecting with $\alpha\in H^*(\overline{\M}_{g,n},\mathbb Q)$ of complementary dimension, one gets group algebra elements:
$$\gen{\gen{\alpha;\gamma}}_{g,B} = \sum_{u\in E[d]} \gen{\alpha;\varphi_u(\gamma)}_{g,\beta(u)}\cdot (u) \in\QQ[E].$$
The $(u)$-coefficient is the reduced GW-invariant for the abelian surface $A(u)$, with polarization $\beta(u)$. The invariant encompasses all possible values of $u$; again, the explicit description of Section~\ref{sec-construction-abelian-families} shows
that all  the polarization  $\beta(u)$ have the same self-intersection but different divisibilities depending on the order of $u$.

\medskip

The reduced GW-cycles $\gen{\gen{\gamma}}_{g,B}$ (and GW-invariants $\gen{\gen{\alpha;\gamma}}_{g,B}$)
can now be seen as functions $\bfF\colon\NN^2\to\QQ[E]$ such that $\bfF(d,a)$ has support in the $d$-torsion of $E$.
The $(u)$-coefficient of $\bfF(d,a)$ is the reduced GW-invariant for a class whose self-intersection only depends on $d$ and $a$, and divisibility depending instead on the order of $u$.

\medskip

Using this group algebra technology, the MCF can be reformulated as an elementary linear functional equation (which we momentarily call \emph{algebraic-MCF}) on this generating function. This abstract algebra setting is the topic of Section \ref{sec-group-algebra}.

\medskip

Furthermore, the correlated reduced decomposition formula expresses the previous generating function with values in $\QQ[E]$ as a sum over \textit{degeneration graphs}/\textit{rigid tropical curves} of \emph{correlated multiplicity functions}, 
which also take values in $\QQ[E]$.

Restricting to the \textit{multiples} of a given decomposition diagram (obtained by scaling the decorations), the correlated multiplicity gives a function $\NN\to\QQ[E]$, already satisfying the above mentioned abstract algebraic-MCF. Summing over graphs without having to concretely enumerate them, we deduce the MCF.

\medskip

Going even deeper, the correlated multiplicity functions are themselves products of correlated invariants functions coming from the moduli spaces of log maps to $(Y|D)$. Functions satisfying the algebraic-MCF turn out to be stable by product. Therefore, the MCF for the multiples of a given graph actually stems from the MCF for each vertex, which amounts to an MCF in the log-CY $E\times\PP^1$, justifying our study of $E\times\PP^1$.

\subsubsection{Strategy II: correlated invariants of $(Y|D)$ from correlated DR-cycles}

To prove the algebraic MCF for the correlated GW-invariants of $(Y|D)$ we first adapt the results of \cite{ranganathan2024logarithmic} expressing logGW-invariants of toric varieties in terms of products of (logarithmic) DR-cycles \cite{holmes2021extending,holmes2025logDR} to correlated invariants of $(Y|D)$. More precisely, in Theorem \ref{theo:rigid-to-rubber}, we prove that the correlated logGW-invariants of $(Y|D)$ can be expressed in terms of the correlated logarithmic DR-cycle for the trivial line bundle with target variety $E$.

\medskip

In \cite{blommecarocci2025DR}, we proved that the correlated DR-cycle admits an explicit Pixton's type formula applying the universal DR-formula proved in \cite{bae2023pixton}. The correlated formula lift, in a straightforward way in the case at hand, to a formula for the logarithmic cycle. Using this explicit tautological expression, we prove a cycle version of the MCF for the correlated logarithmic DR-cycle with target $E$. The MCF for $(Y|D)$ then follows from the MCF for the logDR-cycle and Theorem~ \ref{theo:rigid-to-rubber}.

\subsubsection{Strategy III: exotic insertion and vanishing cohomology }
Following the strategy illustrated so far, one gets the MCF for abelian surfaces for insertions allowed by the degeneration formula: the class $\gamma\in H^*(A)^{\otimes n}$ has to extend to the total space $A_T$ of a semi-stable degeneration of $A$ to $A_0$ on some smooth curve $T$.

In order to retrieve the vanishing cohomology, we apply the results of  \cite[Section~2]{maulik2025gromov}. In the geometric setting of interest, \cite[Theorem~2.4.2]{maulik2025gromov} states that any non-vanishing Gromov-Witten class $\mathrm{ft}_*(\red{\M_{g,m}(A,\beta)}\cap\ev^*\gamma)$ can be recovered from the correlated GW-cycles of $E\times\PP^1$ with so-called \textit{exotic insertions}. The computation in the presence of exotic insertions is also possible using the correlated DR-cycle.

\medskip

With this result at end, the general MCF-formula follows from the MCF for $(Y|D)$ with exotic insertion, which in turn follows from the MCF for the correlated logDR-cycle.

\begin{rem}
   The first step of the strategy, combined with the proof of the MCF for the correlated DR-cycle with target $E$, is where (we feel) the key insight happens. Once taken for granted the correlated DR-formula, it can be appreciated without essentially any background in logarithmic techniques in enumerative geometry.

   By contrast, step II and III of the strategy  heavily rely on recent techniques in logarithmic enumerative geometry \cite{ranganathan2024logarithmic,maulik2025gromov, MaulikRanganathan2025}. 
\end{rem}



\subsection{Results}

We now state our multiple MCF. Each possesses a \textit{group algebra version}, and a \textit{correlator version}. As the group algebra setting has not yet been well introduced, we prefer to  emphasize the correlator version in the introduction. Though the MCF has been stated for GW-invariants in \cite{oberdieck2022gromov}, it is actually true at the cycle level in $H^*(\bar{\M}_{g,n},\QQ)$. We provide both the cycle and the invariant version.

\subsubsection{MCF for correlated DR}

We start with the MCF for the correlated DR-cycle, for which we use the \textit{correlated log DR-cycle formula}, recalled in Section \ref{sec-correlated-DR}. Deleting $g$ and $(Y|D)$ to keep a light notation, we have the following decomposition of the correlated DR-cycle with target $E$ in the group algebra:
$$\mathbf{DR}(a,\bfw) = \sum_{|\bfw|\theta\equiv 0} \operatorname{DR}^\theta(a,\bfw) \cdot (\theta).$$
We state the MCF as a relation between its coefficients.

\begin{theom}[\textbf{\ref{theo:MCF-for-correlated-DR}, MCF for DR}]
   Let $\gamma\in H^*(E^n,\QQ)$ be a geometric insertion. Then, the function $(a,\bfw)\mapsto \ft_*(\mathbf{DR}(a,\bfw)\cap\ev^*\gamma)$ satisfies the $(2g-1+\deg\gamma)$-MCF:
   $$\ft_*(\operatorname{DR}^\theta(a,\bfw)\cap\ev^*\gamma) = \sum_{k|a,\bfw,\theta}k^{2g-1+\deg\gamma}\cdot\ft_*(\operatorname{DR}^{\theta_\prim}(\frac{a}{k},\frac{\bfw}{k})\cap\ev^*\gamma)\in H^*(\bar{\M}_{g,n},\QQ),$$
   where $k|\theta$ means that $\theta$ is of $(|\bfw|/k)$-torsion, and for the $k$-summand, $\theta_\prim$ means a choice of primitive correlator. If furthermore $\alpha\in H^*(\overline{\M}_{g,n},\QQ)$, then $(a,\bfw)\mapsto \mathbf{DR}(a,\bfw)\cap(\ft^*\alpha\cup\ev^*\gamma)$ satisfies the $(3g-3+n-\deg\alpha)$-MCF:
   $$\operatorname{DR}^\theta(a,\bfw)\cap(\ft^*\alpha\cup\ev^*\gamma) = \sum_{k|a,\bfw,\theta}k^{3g-3+n-\deg\alpha}\cdot\operatorname{DR}^{\theta_\prim}(\frac{a}{k},\frac{\bfw}{k})\cap(\ft^*\cup\ev^*\gamma)\in\QQ.$$
\end{theom}

The above formulas are to be understood as follows: the cycle/invariant with a given correlator may be expressed in terms of smaller classes but with primitive correlators. The formula is also true replacing the correlated DR-cycle with target by its logarithmic lift.


\subsubsection{MCF for $E\times\PP^1$}

\begin{theom}[\textbf{\ref{theo-MCF-EP1-general},\ref{coro:MCF-for-EtimesP}, MCF for $E\times\PP^1$}]
Let $\gamma\in H^*(E^n\times (E\times\PP^1)^m,\QQ)$ be a geometric insertion. The function $(a,\bfw)\mapsto \gen{\gen{\gamma}}_{g,a,\bfw}$ satisfies the $(2g-2+\deg\gamma)$-MCF:
$$\gen{\gamma}^\theta_{g,a,\bfw} = \sum_{k|a,\bfw,\theta} k^{2g-2+\deg\gamma}\gen{\gamma}^{\theta_\prim}_{g,a/k,\bfw/k} \in H^*(\bar{\M}_{g,n},\QQ).$$
If furthermore, $\alpha\in H^*(\bar{\M}_{g,n+m},\QQ)$, the function $(a,\bfw)\mapsto \gen{\gen{\alpha;\gamma}}_{g,a,\bfw}$ satisfies the $(3g-3+n+m-\deg\alpha)$-MCF:
    $$\gen{\alpha;\gamma}_{g,a,\bfw}^{\theta} = \sum_{k|a,\bfw,\theta} k^{3g-3+n+m-\deg\alpha}\gen{\alpha;\gamma}_{g,a/k,\bfw/k}^{\theta_\prim} \in\QQ.$$
\end{theom}

The theorem is also valid replacing $\gamma$ by an \textit{exotic insertion} $\gamma^\diamond$ \cite{maulik2025gromov}. We refer to Section~\ref{sec:logintermezzo} for the definition of exotic insertions.

\begin{rem}
The above MCF may be seen as a variant of the MCF satisfied by relative invariants of $(\PP^2,E)$ and the Takahashi conjecture \cite{takahashi1996curves,takahashi2001log}, proven in \cite{bousseau2019takahashi}. In \textit{loc.cit.}, the authors consider counts of rational curves in $\PP^2$ maximally tangent to a smooth cubic $E\subset\PP^2$. For fixed degree, only a finite number of points of the cubic $E$ show up as tangency points, and we can refine the count according to the given point. This is a form of correlator in the $(\PP^2,E)$ setting. The refined count also satisfies an MCF of similar flavor. 

The above theorem suggests that such MCF should hold for more general relative invariants of $(\PP^2,E)$ and other log-CY geometries as well.
\end{rem}

\subsubsection{MCF for Abelian Surfaces}

To state a cycle version of the MCF in abelian surfaces, we adopt the following notation: for $A$ an abelian surface with polarization $\beta$, and $\gamma\in H^*(A^n,\QQ)$, we set
$$\gen{\gamma}_{g,\beta} = \ft_*(\red{\M_{g,n}(A,\beta)}\cap\ev^*\gamma.$$

\begin{theom}[\textbf{\ref{theo-MCF-abelian-general}, MCF for Abelian Surfaces}]
    Let $(A,\beta)$ be a polarized abelian surface. For each $k|\beta$, let $\varphi_k\in SL(H^1(A,\CC))$ be an ismorphism such that $\varphi_k(\beta/k)$ is primitive. Then, we have the following MCF:
    $$\gen{\gamma}_{g,\beta} = \sum_{k|\beta} k^{2g-3+\deg\gamma}
    \gen{\varphi_k(\gamma)}_{g,\varphi^*(\beta/k)} \in H^*(\bar{\M}_{g,n},\QQ).$$
    For $\alpha\in H^*(\overline{\M}_{g,n},\QQ)$, the reduced GW-invariants satisfy the following MCF:
    $$\gen{\alpha;\gamma}_{g,\beta} = \sum_{k|\beta} k^{3g-3+n-\deg\alpha}
    \gen{\alpha;\varphi_k(\gamma)}_{g,\varphi^*(\beta/k)} \in\QQ.$$
\end{theom}

\begin{rem}
Since the abelian surfaces we consider are elliptic fibrations, the isomorphisms $\varphi_k$ (see Remark~\ref{rem:phikdescription} for an explicit description) can be chosen to be shear transforms, acting as the identity on the classes coming from the fiber $E$, but twisting the classes coming from the sections.
\end{rem}

\subsubsection{Illustration}

To conclude, we present as an illustration the relatively easier case where $\gamma$ consists of point insertions. This particular case is presented as a guideline to the whole paper in Section \ref{sec-MCF-point-case}, where we prove the MCF for $E\times\PP^1$ and abelian surfaces. The proof relies on the computation from \cite{blomme2024correlated} which avoids dealing with the correlated DR-cycle. Avoiding the latter, it is possible to witness the power of the group algebra techniques to handle the computations as well as its interaction with the decomposition formula. We thus recover the MCF already proven in\cite{blomme2022abelian3,blomme2025short}, but are also able to provide a refinement with the additional insertion of a $\lambda$-class.

We set $N_{g,\beta} = \gen{\pt^g}_{g,\beta}$.

\begin{theom}[\textbf{\ref{theo-MCF-abelian-points}}]
    Let $\beta$ be a curve class and for each $k|\beta$, let $\widetilde{\beta/k}$ a primitive curve class with the same self-intersection as $\beta/k$. We have the following multiple cover formulas:
    $$N_{g,\beta} = \sum_{k|\beta}k^{4g-3}N_{g,\widetilde{\beta/k}}.$$
\end{theom}


\subsection{Future directions}\label{sec-future-directions}

We expect that the framework developed here will be relevant to discover and prove multiple cover formulas in a variety of other situations. More concretely, the approach through reduced degeneration formula reducing the invariants to double ramification cycle computations is likely to work also to tackle the multiple cover formula conjectures for abelian threefolds \cite{bryan2018curve}.

In another direction, the idea of refining the invariants by keeping track of torsion information, to then exploit the group algebra structure to prove functional equations of the generating series should also apply to the study of logGW-invariants of log CY surfaces, e.g.:  toric surfaces relative to the toric boundary (this has been studied obtaining the expect results in \cite{BCUK2026corrDDR}), general $\mathbb P^1$ bundles on curves, Looijenga pairs. In particular we expect to be able to adapt the developed techniques to prove multiple cover formula in these situations.

Through degeneration formula, this would lead, in particular, to a Gromov-Witten theoretic proof of both the Takahashi conjecture \cite{bousseau2019takahashi} and the KKV formula \cite{pandharipande2016katz}, as well as to the proof of new instances of the MCF for K3 geometries \cite{oberdieck2016curve}.
These questions are object of a current  work in progress.


\subsection{Plan of the paper}

The paper is organized as follows.

\begin{itemize}[leftmargin=0.6cm]
    \item The ideas presented above heavily suggest the study of functions with values in the group algebra $\G=\QQ[E]$, having in mind an MCF structure formula. In Section \ref{sec-group-algebra}, we provide an algebraic framework fully independent of any geometric consideration for functions from a set with an $\NN$-action to $\G$. This turns out to be the right notion to define a class of functions satisfying an ``abstract multiple cover formula'' enjoying suitable structural properties. We also prove that some natural counting functions which show up later in the paper satisfy this abstract MCF.

    We suggest to read Sections~\ref{sec:algebrafunction} and \ref{sec:algebraicmultiplecov} and only go back to Subsection~\ref{sec:refined-groups-counting} during the reading of Section~\ref{sec-correlated-DR}, which is the first point where it is needed.

    \item In Section \ref{sec-reduced-decomposition}, we first recall certain degenerations of abelian surfaces, and state the decomposition formula, using the correlated invariants defined in \cite{blomme2024correlated} and also recapped in the section. We then recall the construction of the reduced class and prove the decomposition formula.

    Sections \ref{sec-reduced-class} and \ref{secred:splitting} contain the technical details  and may be skipped at first reading.

    \item Section \ref{sec-MCF-point-case} provides the proof of the MCF in the case of point insertions for abelian surfaces and $E\times\PP^1$, using the reduced (resp., correlated) decomposition formulas. This may be seen as a toy case of what happens in Sections \ref{SEC-MCFEP1} and \ref{sec:degenerationformula} and is meant to guide the reader through the key steps of the general strategy. In Section \ref{sec-MCF-lambda}, we explain how to adapt the point case to deal with the case where we have an additional $\lambda$-class.


    \item Section \ref{sec:MCFlogcorrelatedDR} deals with the MCF for the correlated (log) DR-cycle. In Section~\ref{sec-correlated-DR} we recall the correlated DR-cycle formula proven in \cite{blommecarocci2025DR} for target  $E\times\PP^1$, and provide a group valued version. We moreover show that invariants obtained from the correlated DR-cycle satisfy an MCF formula. Remarkably, in order to do so we do not need to be able to compute the invariants. The proof boils down to some clever yet elementary combinatorics for counting subgroups of $\ZZ_n^2=(\ZZ/n\ZZ)^2$ with fixed indices. This again can be handled by looking at functions with values in the group algebra and functional equations that they satisfy. These group combinatorics are handled in Section \ref{sec:refined-groups-counting}.

    \item In Section \ref{SEC-MCFEP1}, we adapt the results from \cite{ranganathan2024logarithmic} to compute (correlated) GW-invariants of $E\times\PP^1$ via the (correlated) DR-cycle. Doing so, we are able to prove the MCF fo correlated invariants of $E\times\PP^1$ exploiting the result of the previous section.
   
    \item In Section \ref{sec:degenerationformula}, we adapt the results of \cite{maulik2025gromov} so that the reduced decomposition formula enables one to recover all possible insertions and not only the non-vanishing cohomology classes. Combining this enhanced reduced decomposition formula for abelian surfaces with the MCF formula for $E\times\PP^1$, we prove the MCF for abelian surfaces and any possible insertion.
\end{itemize}

\medskip

\textit{Acknowledgements.} The authors would like to thank Ajith Urundolil Kumaran, Ilia Itenberg, Georg Oberdieck , Rahul Pandharipande and Dhruv Ranganathan for several useful conversations about the topic of this paper. In particular we thank Georg Oberdieck for several useful comments on a first draft of this paper.
Part of this work was realized during respective visits in Neuchâtel and Rome, where both authors would like to thank for the nice working conditions. T.B. acknowledges the support of SNF grant 204125. F.C. is supported by the MIUR Excellence Department Project Mat-
Mod@TOV, CUP E83C23000330006, awarded to the Department of Mathematics,
University of
Rome Tor Vergata, and also acknowledges the support of the PRIN Project
"Moduli spaces and
birational geometry" 2022L34E7W.

\section{Counts with group algebra values}
\label{sec-group-algebra}

In this section we develop an abstract algebraic framework that is well suited to work with correlated GW-invariants (with given insertions), seen as functions of the curve class $\beta$. This framework allows for a  convenient formulation of multiple cover formulas. We also show that some well-known functions admit natural refinements with values in the group algebra and prove that these refinements satisfy multiple cover formulas.

\subsection{Group algebras and functions}\label{sec:algebrafunction}

Let $G$ be the group $(\RR/\ZZ)^r$ for some $r\geqslant 1$, and let $\G$ be its group algebra, with coefficients in $\CC$ or $\QQ$. The generator associated to $\theta\in G$ is denoted by $(\theta)$. The set of torsion points $(\QQ/\ZZ)^r$ can also be seen as the direct limit $\lim\ZZ_n^r$. We are mostly interested in the case where $G$ is  an elliptic curve, isomorphic as a group to $(\RR/\ZZ)^2$ with torsion points $(\QQ/\ZZ)^2$.

\subsubsection{Some operators on $\G$}

\begin{defi}
    We define the following operators on $\G$:
    \begin{itemize}
        \item The multiplication operator defined on generators as follows: if $d\in\ZZ$,
        $$\m{d}(\theta) = (d\cdot\theta).$$
        \item The division operator defined on generators as follows: if $d\neq 0$,
        $$\d{\frac{1}{d}}(\theta) = \frac{1}{d^r}\sum_{d\tau\equiv\theta} (\tau),$$
        where $d^r$ is the cardinality of the subroup of $d$-torsion elements in $G.$ The division operator takes the average over $d$-roots in $(\RR/\ZZ)^r$.
    \end{itemize}
\end{defi}
\begin{prop}
    The multiplication and division operators respect the multiplication: if $\mathbf{x},\mathbf{y}\in\G$, we have
    $$\m{d}(\mathbf{x}\mathbf{y}) = \m{d}(\mathbf{x})\m{d}(\mathbf{y}) \text{ and }\d{\frac{1}{d}}(\mathbf{x}\mathbf{y}) = \d{\frac{1}{d}}(\mathbf{x})\d{\frac{1}{d}}(\mathbf{y}).$$
\end{prop}

We denote by $T_n = \d{\frac{1}{n}}\left((0)\right)$, the average of $n$-torsion elements. These satisfy the following properties:
    \begin{enumerate}
        \item $T_nT_m=T_{\lcm(n,m)}$,
        \item $\m{d}T_n=T_{n/\gcd(d,n)}$,
        \item $\d{\frac{1}{d}}T_n=T_{nd}$.
    \end{enumerate}
In particular, they are idempotent elements. Notice that the image of $\d{\frac{1}{d}}$ is invariant by multiplication by $T_d$. Moreover, if $\theta_0$ is any $d$-root of $\theta$,
$$\d{\frac{1}{d}}(\theta)=(\theta_0)T_d.$$
The division operator gives an isomorphism  $\G\xrightarrow{\d{\frac{1}{d}}}\G\cdot T_d$, whose inverse is given by $\m{d}$.

\begin{expl}
    If $G=\CC/\gen{1,\tau}$, then $T_2$ is the average of $2$-torsion elements:
    $$T_2=\frac{1}{2^2}\left[ (0)+(1/2)+(\tau/2)+((1+\tau)/2)\right].$$
\end{expl}

Let $\G_n$ be the finite dimensional subalgebra of $\G$ generated by the $T_k$ for $k|n$, and $\G_\infty$ the subalgebra generated by all $T_n$. The $(T_k)_{k|n}$ form a basis of $\G_n$ as a vector space (over $\mathbb{K}=\mathbb C$ or $\mathbb Q)$.

\medskip

We denote by $\Prim_n$ the function on $\G_n$ defined by
$$\begin{array}{rccl}
  \Prim_n\colon&\G_n&\longrightarrow & \mathbb K  \\
   & \sum_{k|n}c_k T_k &\longmapsto & c_n
\end{array}$$
that maps an element to its $T_n$-coefficient, which we call \textit{$n$-primitive coefficient}. 

\subsubsection{Algebra of sequences in $\G$}

We consider the space $\G^{\NN}$ of sequences indexed by the multiplicative monoid $\NN$ with values in the group algebra $\G$; i.e. $F\in\G^{\NN}$ is a function from $\mathbb N$ to the group algebra $\G.$
This is a ring with the component-wise addition and multiplication. Since $\NN$ has a monoid structure, $\G^{\NN}$ is also endowed with the \emph{convolution product} defined as follows:
$$(F\ast G)(\delta) = \sum_{kl=\delta}F(k)G(l).$$

We now put our focus on a specific type of sequences called of \textit{  diagonal type}.

\begin{defi}
We say that a sequence $F\in\G^{\NN}$ is of \emph{diagonal type} if $F(\delta)\in\G_\delta$: i.e., $F(\delta)$ belongs to the span of $T_k$ for $k|\delta$ for each $\delta\in\mathbb N$. We denote the subset of  $\G^{\NN}$  of diagonal type sequences by $\A$.
\end{defi}

It follows immediately from the fact that $\G_{\delta}$ is a sub-algebra that:

\begin{prop}
    The subset of diagonal type sequences $\A$ is a subalgebra of $\G^{\NN}$, namely it is stable by component-wise sum, component-wise product and convolution.
\end{prop}

We also define $\Prim\colon\A\to\KK^{\NN}$ as follows:
$$\Prim(F)(\delta) = \Prim_\delta(F(\delta))\in\KK.$$
This map is of course linear but it is not a morphism of algebra as it does not respect the product. Computing the primitive coefficient of a product in $\G_\delta$ actually requires some care. The reason is that we may have $T_kT_l=T_n$ even if $k$ and $l$ are strict divisors of $n$.

\subsubsection{$\NN$-modules}
\label{sec:N*modules}

We now introduce a slight generalization of the previous set-up, where instead of having functions $\NN\to\G$ we have functions from what we call a $\mathbb N$-module $X$ with values in $\G$.

\begin{defi}
    A $\NN$-module is the data of a set $X$ endowed with the following structure:
    \begin{enumerate}[label=(\roman*)]
    \item a free action $\NN\times X\to X$ such that each $x\in X$ only has a finite number of divisors,
    \item a map $x\in X\mapsto |x|\in\NN$ called \textit{norm} compatible with the $\NN$-action in the sense that $|k\cdot x|=k|x|$.
\end{enumerate}

Elements in $X$ which have no non-trivial divisors are called \textit{primitive}. As any $x\in X$ has only a finite number of divisors and the action is free, any element $x$ can uniquely be expressed as $\delta\cdot\widetilde{x}$ where $\widetilde{x}$ is a primitive element.

The integer $\delta=\ell(x)$ is called the integral length of $x$. Furthermore, we have $|x|=\ell(x)|\widetilde{x}|$. Notice that a primitive element needs not to have norm $1$. A morphism of $\NN$-module is a map $f:X\to Y$ compatible with the monoid action and the norm:
$$f(k\cdot x)=k\cdot f(x) \text{ and }|f(x)|=|x|.$$
\end{defi}

\begin{expl}\label{ex:Nmodules}
We give some examples of $\NN$-modules which will be relevant in the later sections.
\begin{itemize}[label=$\circ$,leftmargin=0.6cm]
    \item We can consider $\NN$ as a $\NN$-module with the action on itself given by multiplication. However, we also have to choose a norm, fully determined by $|1|$. We thus define the $d$-norm by $|k|_d=dk\in\NN$.

    In fact, if $x\in X$ is a general $\NN$-module, the orbit $\NN\cdot x$ of $x$ is a $\NN$-submodule of $(X,|\cdot |)$
    isomorphic to $\NN$ with the $|x|$-norm.

    \item The set of \textit{diagram degrees} $\BBB$ which, as a set, is just $\NN^2$ can be equipped with an $\NN$-module structure: $\NN$ acts by scaling each factor, and the norm is defined by the projection onto the first factor.


    \item The set of decorated diagrams $\DDD$ as defined in \cite{blomme2022abelian3} and reviewed in Section \ref{sec:pearl-diagrams}. The action of $\NN$ is by scaling edge weights $w_e$ and vertex degrees $a_v$. The norm is the degree of the associated tropical cover. In other words, we have a degree map $\DDD\to\BBB$ compatible with the $\NN$-action, and the norm on $\DDD$ is just the pull-back of the norm of $\BBB$: if $\deg\Dfk=B$, we have $|\Dfk|=|B|$.

    \item The ramification $\NN$-module $\RRR=\{(a,\bfw)\in\NN\times\ZZ^n\backslash\{0\} \text{ s.t. }\sum_1^n w_i=0\}$, with norm $|(a,\bfw)|$ given by the g.c.d. of the $w_i$.

    \item The $\NN$-module $\PPP_r=\{(d,n)\in\NN^2 \text{ such that }d|n^r\}$, of possible indices for subgroups of $\ZZ_n^r$. The $\NN$-action is given by the usual action on each coordinate and norm is given by the second coordinate. Notice that primitive elements are not the same as in $\NN^2$. For instance, $(4,2)\in\PPP_2$ but $(2,1)$ does not so that $(4,2)$ is primitive. 
\end{itemize}
\end{expl}

Until the end of the section, we consider a $\NN$-module $X$ with a norm $|\cdot|$. Let $\G^X$ be the set of functions $X\to\G$. This is a group endowed with an action of $\G^{\NN}$ by convolution: if $F\in\G^{\NN}$ and $U\in\G^X$, we set
$$(F\ast U)(x) = \sum_{k|x}F(k)\cdot U(x/k).$$

We now define a notion of diagonal type function similar to the case of sequences.

\begin{defi}
We say that $U\in\G^X$ is of diagonal type if $U(x)\in\operatorname{Span}(T_{|x|/k})_{k|x}$. We denote by $\X$ the subset of diagonal type functions in $\G^X$. The latter is endowed with an action of $\A$ by restricting the convolution operation.
\end{defi}

We also have a primitive coefficient function $\Prim\colon\X\to\CC^X$ defined by
$$\Prim(U)(x) = \Prim_{|x|}(U(x))\in\CC.$$

\subsection{Algebraic multiple cover formula}\label{sec:algebraicmultiplecov}

We now introduce the abstract concept of \textit{multiple cover formula} in the current algebraic setting. We are interested in diagonal type sequences which can be reconstructed from the sequence of their primitive coefficients, namely $F$ can be recovered from $\Prim(F)$, which is a sequence with values $\CC$ rather than $\G$.

\subsubsection{Multiple cover formulas in $\A$}

If $\alpha\in\ZZ$, we denote by $\epsilon_\alpha(k)=k^\alpha$ the $\alpha$-power function. If $F$ is a diagonal type sequence, we construct a new sequence from the primitive coefficients of $F$, also of diagonal type, using the convolution with a power function:
    $$(\epsilon_\alpha\ast(\Prim(F)T))(\delta) = \sum_{kl=\delta}k^\alpha\Prim_l(F(l))\cdot T_l,$$
where we denoted by $T\colon n\in\NN\mapsto T_n\in\G$.

\begin{defi}
    We say that a sequence $F\in\G^{\NN}$ of diagonal type satisfies the multiple cover formula of weight $\alpha$ (or shortly the $\alpha$-MCF) if we have
    $$F=\epsilon_\alpha\ast(\Prim(F)T).$$
\end{defi}

The term \textit{multiple cover formula} here refers to diagonal type $\G$-sequences, but the concept stems from Gromov-Witten theory. Multiple cover formulas are conjectured to be satisfied by certain generating series of GW invariants, usually considered as functions of the curve class. 

The present paper illustrates how $\G$-sequences can be used to encode GW invariants for certain geometries in such way that multiple cover formula phenomena in GW-invariants admit an equivalent reformulation in terms of multiple cover formula for diagonal type $\G$-sequences. From now on, we use the term MCF to mean multiple cover formula.

\begin{rem}
    For a function to satisfy the $\alpha$-MCF means that it can be reconstructed from its primitive coefficient. There exists other formulations that we do not develop here, for instance using instead the sequence of $(0)$-coefficients.
\end{rem}

We now prove some structural properties.

\begin{lem}\label{lem:monomial-rule}
    If $F$ satisfies the $\alpha$-MCF, then $\epsilon_r F$ satisfies $(\alpha+r)$-MCF.
\end{lem}

\begin{proof}
    Assume $F$ satisfies the $\alpha$-MCF, then we have
    \begin{align*}
       (\epsilon_r F)(\delta)=\delta^rF(\delta) = & \delta^r \sum_{kl=\delta}k^\alpha\Prim_l(F)(l) T_l \\
        = & \sum_{kl=\delta}k^rk^\alpha\cdot l^r\Prim_l(F(l)) T_l \text{ since }\delta^r=k^r l^r,\\
        = & \sum_{kl=\delta}k^{\alpha+r}\cdot\Prim_l(l^rF(l)) T_l \text{ since }\Prim_l\text{ is }\CC\text{-linear}.\\
    \end{align*}
\end{proof}

In particular, it means that every function $F$ satisfying the $\alpha$-MCF can uniquely be expressed as the product of the monomial $\epsilon_\alpha$ and a function satisfying the $0$-MCF. The following Lemma explains the advantage in working with functions satisfying the $0$-MCF.

\begin{lem}\label{lem:product-rule}
    Assume $F$ and $G$ satisfy the $0$-MCF. Then the term by term product $FG$ also satisfies the $0$-MCF.
\end{lem}

\begin{proof}
    Assume $F$ and $G$ satisfy the $0$-MCF. Then we have
    \begin{align*}
        (FG)(\delta) = & \left(\sum_{k|\delta}\Prim_k(F(k))T_k\right)\left(\sum_{l|\delta}\Prim_l(G(l))T_l\right) \\
        = & \sum_{k,l|\delta}\Prim_k(F(k))\Prim_l(G(l))T_{\lcm(k,l)} \\
        = & \sum_{d|\delta} T_d\sum_{\substack{k,l|\delta \\ \lcm(k,l)=d}} \Prim_k(F(k))\Prim_l(G(l)).
    \end{align*}
    In particular, taking the $T_\delta$-coefficient, we deduce that
    $$\Prim_\delta(F(\delta) G(\delta)) = \sum_{\substack{k,l|\delta \\ \lcm(k,l)=\delta}} \Prim_k(F(k))\Prim_l(F(l)).$$
    Furthermore, we can rewrite the inner sum as follows: summing over $k,l|\delta$ with $\lcm(k,l)=d$ (and $d |\delta$) is the same as summing over $k,l|d$ with $\lcm(k,l)=d.$ We now recognize that the inner sum is the $d$-primitive coefficient of $(FG)(d)$ and we get that the product also satisfies the $0$-MCF:
    $$(FG)(\delta) = \sum_{d|\delta} \Prim(FG)(d)T_d.$$
\end{proof}

\begin{expl}
    The function $F(\delta)=T_1$ satisfies the $0$-MCF. Indeed, we have $\Prim_\delta(F(\delta))=1$ if $\delta=1$ and $0$ else. Therefore,
    $$(\epsilon_0\ast\Prim(F)T)_\delta = \sum_{k|\delta}\Prim_k(F(k))T_k = T_1.$$
    More generally, functions mapping to a multiple of $T_1$ satisfying the $\alpha$-MCF are of the form $\delta\mapsto \lambda\delta^\alpha T_1$ where $\lambda\in\CC$.
    If $f\colon\NN\to\CC$ is any function, then $F(\delta)=\sum_{k|\delta}f(k)T_k$ is a diagonal type function tautologically satisfying the $0$-MCF. In fact, it is the unique function satisfying the $0$-MCF with primitive coefficient given by $f$.
\end{expl}

\subsubsection{Examples with sum of divisors}

Let $m\geqslant 0$. We have the well-known powered sum of divisors functions defined as follows:
$$\sigma_{m}(a) = \sum_{k|a}\left(\frac{a}{k}\right)^{m}.$$
The reason for writing $\frac{a}{k}$ instead of $k$ for the divisors will be clear in a moment. If we do not write the index, the latter is understood to be $1$ and we just make the sum of divisors.

\medskip

We consider the following deformation of the function above, which takes values in $\G_\delta$:
$$\bsigma^\delta_{m}(a) = \sum_{k|a}\left(\frac{a}{k}\right)^{m}T_{\delta/\gcd(\delta,k)}.$$
This is a refinement of the divisor function in the sense that if we consider the
degree  morphism $\G\to\CC$ mapping each generator $(\theta)$ to $1$, so that each $T_m$ is also mapped to $1,$ then 
the total degree of $\bsigma^\delta_{m}(a)$ is the classical function $\sigma_{m}(a)$. 
The next lemma summarizes some elementary properties of this refinement. We will refer to  property (i), (ii) of the Lemma as the \textit{unrefinement relations}.

\begin{lem}\label{lem:arithmeticproperties}
    The refined sum of divisors $\bsigma^\delta$ (as well as $\bsigma_{m}^\delta$) functions satisfy the following properties:
    \begin{enumerate}[label=(\roman*)]
        \item if $\gcd(\delta,d)=1$, we have $\m{d}\bsigma^\delta = \bsigma^{\delta}$, and if $d|\delta$ we have $\m{d}\bsigma^\delta=\bsigma^{\delta/d}$;
        \item $\m{d}\bsigma^\delta = \bsigma^{\delta/\gcd(d,\delta)}$; .
        \item if $d|\delta$ we have $\d{\frac{1}{\delta/d}}\bsigma^d = \bsigma^\delta T_{\delta/d}$,
        \item $\bsigma^\delta(a)=\bsigma^\delta(a)T_{\delta/\gcd(a,\delta)}$,
        \item $\bsigma^\delta(a) = \d{\frac{1}{\delta/\gcd(\delta,a)}}\bsigma^{\gcd(\delta,a)}(a)$.
        \item The function $\bsigma^\delta(a)$ is multiplicative in the sense that
        $$\bsigma^\delta(a) = \prod_p \bsigma^{p^{\nu_p(\delta)}}(p^{\nu_p(a)}).$$
    \end{enumerate}
\end{lem}

\begin{proof}
    \begin{enumerate}[label=(\roman*)]
    \item To prove the first point, recall that $T_n$ satisfies that $\m{d}T_n=T_{n/\gcd(d,n)}.$
    Then if $\gcd(d,\delta)=1$,  for any $k|\delta$ we have $\m{d}T_{\delta/\gcd(\delta,k)}=T_{\delta/\gcd(\delta,k)}$. If we have instead that $d|\delta$, then 
    $$\m{d}T_{\delta/\gcd(\delta,k)} = T_{\frac{\delta/d}{\gcd(\delta/d,k)}}.$$
    To see this last identity, we need to show that
    $$\frac{\delta}{\gcd(\delta,k)} = \gcd\left(d,\frac{\delta}{\gcd(\delta,k)}\right)\frac{\delta/d}{\gcd(\delta/d,k)},$$
    which is elementary.
    
    \item We can always write $d=\gcd(d,\delta)d'$ where $d'$ is coprime with $\frac{\delta}{\gcd(d,\delta)}$, and we apply successively both points of (i).

    \item From the first point we have that $\m{\delta/d}\bsigma^\delta = \bsigma^d$. We then apply on both sides $\d{\frac{1}{\delta/d}}$ and use that composition between multiplication and division  is multiplication by $T_{\delta/d}$.

    \item If $k|a$, we have that $\gcd(\delta,k)|\gcd(\delta,a)$, and thus $\frac{\delta}{\gcd(\delta,k)}$ is always divisible by $\frac{\delta}{\gcd(\delta,a)}$. In particular, any term in the sum is invariant when multiplying by $T_{\delta/\gcd(\delta,a)}$.

    \item Follows from (iii) and (iv).

    \item Follows from the multiplicativity of $T$ and of the sum of divisors function.
    \end{enumerate}
\end{proof}

We can now give the first family of examples of functions in $\G^{\NN}$ satisfying the $0$-MCF.

\begin{prop}
\label{prop-divisor-fct}
    For $a\in\NN$, the sequence $(\bsigma_{m}^\delta(\delta a))$ is of diagonal type and satisfies the $0$-MCF.
\end{prop}

\begin{proof}
    The diagonal type property follows from the definition. We now want to prove that
    $$\bsigma_{m}^\delta(\delta a) = \sum_{k|\delta} \Prim_k(\bsigma_{m}^k(ka))T_k.$$
    We start from the expression of $\bsigma_{m}^\delta$ to identify the primitive coefficient: only the terms with $\gcd(\delta,k)=1$ contribute to the $T_\delta$-coefficient. Therefore, we have
    $$\Prim_\delta(\bsigma_{m}^\delta(\delta a)) = \sum_{\substack{k|a \\ \gcd(k,\delta)=1}}\left(\frac{\delta a}{k}\right)^{m}.$$
    Now, we identify the the $T_l$-coefficient of $\bsigma_{m}^\delta(\delta a)$: the gcd between $k$ and $\delta$ has to be $\delta/l$, and we get
    $$\operatorname{Coeff}_{T_l}\left(\bsigma_{m}^\delta(\delta a)\right) = \sum_{\substack{k|\delta a \\ \gcd(k,\delta)=\delta/l}} \left(\frac{\delta a}{k}\right)^{m}
        = \sum_{\substack{k'|l a \\ \gcd(k',l)=1}}\left(\frac{\delta a}{(\delta/l)k'}\right)^{m}
        =  \sum_{\substack{k'|l a \\ \gcd(k',l)=1}} \left(\frac{l a}{k'}\right)^{m},$$
    where to get to the second step we wrote $k=\frac{\delta}{l}\cdot k'$ (and $\delta=\frac{\delta}{l}\cdot l$) with $\gcd(k',l)=1$. The condition $k|\delta a$ is then equivalent to $k'$ dividing $\delta a/(\delta/l) = la$. The coefficient is thus equal to $\Prim_l(\bsigma_{m}^l(l a))$ and we conclude.
\end{proof}

\begin{expl}
In particular, using Lemmas \ref{lem:monomial-rule} and \ref{lem:product-rule}, we deduce a wide family of functions satisfying MCF: we just have to multiply a monomial with a product of functions of the above type: $\delta\mapsto \delta^r\prod_v\bsigma^{\delta}(\delta a_v)$.
\end{expl}

\subsubsection{Extension to sets with a $\NN$-action}

We can now define what it means for a function $U\in\X$ to satisfy the $\alpha$-MCF, and then explain how to get back to the case of $\NN$ with diagonal type sequences.

\begin{defi}\label{def:MCFmodule}
    We say that a diagonal type sequence $U\in\X$ satisfies the $\alpha$-MCF if we have
    $$U(x) = \sum_{k|x} k^\alpha \Prim(U)(x/k)T_{|x|/k}.$$
\end{defi}


The set of functions satisfying the $\alpha$-MCF is stable by sum. Furthermore, restricting to orbits, we also have the following lemma.

\begin{lem}\label{lem:reductiontoN}
A function $U\in \X$ satisfies the $\alpha$-MCF if and only if for any primitive $\widetilde{x}\in X$, $U|_{\NN\cdot \widetilde{x}}$ satisfies the $\alpha$-MCF, where $\NN\cdot \widetilde{x}$ is the orbit of $\widetilde{x}$, seen as a $\NN$-module with the induced norm.    
\end{lem}

In other words, we can always restrict to orbits of primitive elements to check the MCF. We now restrict to such an orbit, which as observed before is a sub-$\NN$-module of $(X, |\cdot|)$ isomorphic to $(\NN,|\cdot|_d)$ with $d=|\widetilde{x}|$.

\medskip

Let $F\colon (\NN,|\cdot|_d)\to\G$ be a function. We say that it is invariant by $d$-torsion if for any $\delta$, $F(\delta)T_d=F(\delta)$.

\begin{lem}\label{lem:MCF-for-N-k}
    Let $F\colon (\NN,|\cdot|_d)\to\G$ be a function. Then $F$ satisfies the $\alpha$-MCF if and only if $F$ is invariant by $d$-torsion and $\m{d}\circ F$ satisfies the $\alpha$-MCF for usual functions on $\NN$.
\end{lem}

\begin{rem}\label{rem:invariantformT}
We notice that if $F$ satisfies the $\alpha$-MCF for $|\cdot|_d$, it has to be invariant by $d$-torsion since the expression $\sum_{k|\delta}k^\alpha\Prim_{d\cdot\delta/k}(F(\delta/k))T_{d\cdot(\delta/k)}$ is invariant, simply because $T_{d\cdot(\delta/k)}T_d=T_{d\cdot(\delta/k)}.$
\end{rem}

\begin{proof}
    Assume that $F$ satisfies the $\alpha$-MCF for $|\cdot|_d$, i.e. unraveling the definition:
    $$F(\delta) = \sum_{k|\delta} k^\alpha\cdot \Prim_{d\cdot(\delta/k)}(F(\delta/k))T_{d\cdot(\delta/k)}.$$
    It is in particular invariant by $d$-torsion. Applying $\m{d}$, we get:
    $$\m{d}(F(\delta)) = \sum_{k|\delta} k^\alpha\cdot \Prim_{d\cdot\delta/k}(F(\delta/k))T_{\delta/k}.$$
    The $\delta$-primitive coefficient is obtained taking $k=1$, and  thus we have
    $$\Prim_\delta(\m{d}F(\delta)) = \Prim_{d\cdot\delta}(F(\delta)),$$
    which substituted in the previous equation gives the $\alpha$-MCF for $\m{d}\circ F$. 
    
    Conversely, let $G(\delta)=\m{d}(F(\delta))$ and assume it satisfies the $\alpha$-MCF. We have
    $$F(\delta)=F(\delta) T_d=\d{\frac{1}{d}}\circ\m{d}(F(\delta)) = \d{\frac{1}{d}}(G(\delta)),$$
    and we can conclude with a similar computation:
    \begin{align*}
        G(\delta) & = \sum_{k|\delta}k^\alpha \Prim_{\delta/k}G(\delta/k)\cdot T_{\delta/k} , \\
        \d{\frac{1}{d}}(G(\delta)) & = \sum_{k|\delta}k^\alpha \Prim_{\delta/k}G(\delta/k)\cdot T_{d\cdot(\delta/k)}.
    \end{align*}
\end{proof}

Using Lemma~\ref{lem:reductiontoN} and Lemma~\ref{lem:MCF-for-N-k}, checking the $\alpha$-MCF for $U\in\X$ can be reduced to check the $\alpha$-MCF for certain diagonal type sequences $U\rvert_{\NN\cdot\widetilde{x}}\in\A$, and we can use the results of the previous section.

\begin{expl}
    The functions $\delta\mapsto \bsigma^{d\delta}(\delta a)\cdot T_d$ satisfy the $0$-MCF for $|\cdot|_d$ since the multiplication by $T_d$ makes them invariant by $d$-torsion and
    $$\m{d}\left(\bsigma^{d\delta}(\delta a)\cdot T_d\right) = \bsigma^\delta(\delta a),$$
    using the unrefinement relation. In particular, products of such functions also satisfy the $0$-MCF.
\end{expl}

\begin{prop}\label{prop:push-forward}
    Let $f\colon X\to Y$ be a morphism of $\NN$-modules with finite fibers. Let $U\in\X$ be a diagonal type function on $X$. We define the push-forward $f_*U\in\Y$ by
    $$f_*U(y) = \sum_{f(x)=y}U(x).$$
    If $U$ satisfies the $\alpha$-MCF, so does $f_*U$.
\end{prop}

\begin{proof}
    First, we can restrict to the orbit of a primitive element $\widetilde{y}\in Y$. If $f(x)=\widetilde{y}$, then $x$ is primitive as well as $x$ being divisible would imply $\widetilde{y}$ being divisible. Then, we can consider the preimage of $\NN\cdot\widetilde{y}$, which is a (disjoint) union of orbits of primitive elements: the action being free, if $f(kx)\in \NN\cdot\widetilde{y}$, then $f(x)\in\NN\cdot\widetilde{y}$. Restricting to a unique orbit, the map $f$ becomes injective, and the push-forward is just the extension by $0$. The $\alpha$-MCF is thus true because it is satisfied on the orbit. We recover the formula summing over the different orbits.
\end{proof}

A diagonal type function $F\colon X\to\G$ satisfies the $\alpha$-MCF if and only if its restriction to the orbit of every \textit{primitive} element does. However, the restriction to the orbit of a non-primitive element has no reason to do so, due to the existence of divisors outside of the orbit.

Alternatively, given $F\colon\NN\to\G$ satisfying the $\alpha$-MCF and $l\geqslant 1$, so that the orbit $\NN\cdot l$ is isomorphic to $(\NN,|\cdot|_l)$ as a $\NN$-module, the function $F|_{\NN\cdot l}\colon \delta\mapsto F(l\delta)$ defined on $(\NN,|\cdot|_l)$ does not usually satisfy the $\alpha$-MCF. The main obstruction is that it has no reason to be invariant by $l$-torsion. The following lemma shows that to remedy the problem, one just has to multiply by $T_l$.

\begin{lem}\label{lem:restriction-suborbit}
    Let $F\colon (\NN,|\cdot|_r)\to\G$ satisfying the $\alpha$-MCF and $l\geqslant 1$. The orbit of $l$ is isomorphic to $(\NN,|\cdot|_{rl})$. We consider
    $$G(\delta)=T_{rl}\cdot F(l\delta).$$
    It satisfies the $\alpha$-MCF as a function $(\NN,|\cdot|_{rl})\to\G$.
\end{lem}

\begin{proof}
    By construction, $G(\delta)$ is invariant by $rl$-torsion. Therefore, to check the MCF, one has to prove that $\delta\mapsto \m{rl}(G(\delta))=\m{rl}(F(l\delta))$ satisfies the $\alpha$-MCF.
    
    First of all, up to applying $\m{r}$, we assume that $r=1$. As $F(l\delta)\in\G_{l\delta}$, applying $\m{l}$, we have $\m{l}(G(\delta))\in\G_\delta$ so that it is of diagonal type. Then, up to multiplying by $\delta^{-\alpha}$, we can assume that $\alpha=0$ and $F$ satisfies the $0$-MCF. In particular, we get that
    \begin{align*}
        F(l\delta) = \sum_{k|l\delta} \Prim_kF(k)\cdot T_k  = & \sum_{u|l}\sum_{\substack{k|l\delta \\ \gcd(k,l)=u}} \Prim_kF(k)\cdot T_k \\
        = & \sum_{u|l}\sum_{\substack{k'|\delta \\ \gcd(k',\frac{l}{u})=1}} \Prim_{uk'}F(uk')\cdot T_{uk'},
    \end{align*}
    where we sort out terms according to the value of $\gcd(k,l)$. To get to the last equality, we write $k=uk'$ with $\gcd(k',\frac{l}{u})=1$. In particular, we need to have $k'|\frac{l}{u}\delta$, which is equivalent to $k'|\delta$ since $k'$ is coprime with $\frac{l}{u}$. For the same reason, we also have
    $$\m{l}T_{uk'} = \m{\frac{l}{u}}\m{u}T_{uk'} = \m{\frac{l}{u}}T_{k'} = T_{k'}.$$
    In the end, we get the following expression for $G$:
    $$\m{l}(G(\delta)) = \sum_{k|\delta}T_k\left(\sum_{\substack{u|l \\ \gcd(k,\frac{l}{u})=1}} \Prim_{uk}F(uk)\right).$$
    We deduce the expression of the primitive coefficient taking $k=\delta$:
    $$\Prim_\delta \left(\m{l}(G(\delta))\right) = \sum_{\substack{u|l \\ \gcd(\delta,\frac{l}{u})=1}} \Prim_{u\delta}F(u\delta).$$
    Now, we simply need to identify the $T_k$-coefficient as $\Prim_k \m{l}(G(k))$, which follows from its expression.
\end{proof}

\subsection{Various counts of subgroups in $\ZZ_n^2$}
\label{sec:refined-groups-counting}

Let us denote by $\ZZ_n^2=(\ZZ/n\ZZ)^2$ which we identify with the subgroup of $n$-torsion elements in $(\mathbb R/\ZZ)^2$. In this section, which may be skipped at first reading, we provide a second family of functions with values in a group algebra that satisfy the $0$-MCF. These come from a refined count of subgroups in $\ZZ_n^2$.

\subsubsection{Generalities and refined count of subgroups}

For $K$ a subgroup of $\ZZ_n^2$, we set
$$T_K=\frac{1}{|K|}\sum_{\theta\in K}(\theta),$$
the average of elements in $K$;  this is a group algebra element in $\G=\QQ[(\mathbb R/\ZZ)^2]$. Among them, we set $T_k=T_{\ZZ_k^2}$ the average of $k$-torsion elements.

We have an action of $SL_2(\ZZ)$ on $(\RR/\ZZ)^2$ and the action restricts to each torsion subgroup $\ZZ_n^2$. The results in the following lemma come from well-known facts on the structure of subgroups of $\ZZ^2$. The last point is essential as it allows to deduce results for general finite abelian groups from the study for abelian $p$-groups only.

\begin{prop}
\label{prop:elementary-stuff-on-subgroups}
We have the following:
\begin{enumerate}[leftmargin=0.6cm]
    \item $SL_2(\ZZ)$ acts transitively on torsion elements of a given order;
    \item an element of $\QQ[\ZZ_n^2]$ is invariant by $SL_2(\ZZ)$ if and only if it belongs to the span of the $(T_k)_{k|n}$.
    \item Every finite abelian group $G$ uniquely decomposes as a product of abelian $p$-groups $\prod G_p$. Furthermore, the isomorphism induces between group algebras
    $$\bigotimes_p\QQ[G_p]\simeq\QQ[G].$$
    Every subgroup $H$ can uniquely be written as a product of $p$-groups $\prod H_p$ where $H_p$ is a subgroup of $G_p$.
\end{enumerate}
\end{prop}

We denote by $\pi\colon\ZZ_n^2\to\ZZ_n$ the projection onto the first coordinate, which induces a morphism between the group algebras which we denote by $\pi_*$. Denoting by
$$\bar{T}_k=\frac{1}{k}\sum_{\theta\in\ZZ_k}(\theta)\in\QQ[\ZZ_n]\subset\QQ[\RR/\ZZ],$$
which is the average over $k$-torsion elements in $\RR/\ZZ$, we have $\pi_*T_k=\bar{T}_k$. In particular, by mapping a basis to another, $\pi_*$ induces an isomorphism from the subalgebra of $SL_2(\ZZ)$-invariants elements of $\QQ[\ZZ_n^2]$ to the subalgebra of $\QQ[\ZZ_n]$ spanned by $(\bar{T}_k)_{k|n}$. Therefore, for any $SL_2(\ZZ)$-invariant $\mathbf{g}\in\QQ[\ZZ_n^2]$, $\pi_*\mathbf{g}\in\QQ[\ZZ_n]$ uniquely determines $\mathbf{g}$. Therefore, to get a formula for an $SL_2(\ZZ)$-invariant element, it suffices to find a formula for its image by $\pi_*$.

\begin{expl}
\label{expl-computation-subgroup}
    To illustrate the previous strategy, we provide an expression for the group algebra count of subgroups of given index: for $d|n^2$, let
    $$\bfM(d,n) = \sum_{H\colon |\ZZ_n^2:H]=d} T_H.$$
    Since the action of $SL_2(\ZZ)$ preserves the index of subgroups, the element $\bfM(d,n)$ is $SL_2(\ZZ)$-invariant.

    For every subgroup $H$ of $\ZZ_n^2$, there exists a unique matrix $M_H=\left(\begin{smallmatrix} a & 0 \\ j & b \end{smallmatrix}\right)$ whose columns span $H$, with $a,b|n$ and $0\leqslant j<b$ is of $(n/a)$-torsion modulo $b$. They are obtained as follows: $\binom{0}{b}$ spans $\ker\pi|_H$, $a$ spans $\pi(H)\subset\ZZ_n$, and $\binom{a}{j}$ is a lift of $a$ by $\pi$. The integer $j$ is only defined modulo $b$ so that it may be chosen $0\leqslant j<b$, and since $(n/a)\cdot\binom{a}{j}$ belongs to $\ker\pi$, $j$ is $(n/a)$-torsion modulo $b$. Furthermore, we have $[\ZZ_n^2:H]=ab$ and $\pi_*T_H=\bar{T}_{n/a}$. 

    We can thus compute $\pi_*\bfM(d,n)$ grouping the subgroups $H$ according to the value of $a$:
    $$\pi_*\bfM(d,n) = \sum_{(a,b,j)} \bar{T}_{n/a},$$
    where the sum is over $a,b|n$ such that $ab=d$ and $0\leqslant j<b$ is of $(n/a)$-torsion. The coefficient of $T_{n/a}$ is the number of $j$ of $(n/a)$-torsion modulo $b=\frac{d}{a}$, which is $\gcd(n/a,d/a)=\frac{\gcd(n,d)}{a}$. One thus gets
    $$\pi_*\bfM(d,n) = \sum_{a|d,n \text{ s.t. }\frac{d}{a}|n} \frac{\gcd(n,d)}{a} \bar{T}_{n/a} \; \text{ and hence } \;\bfM(d,n) = \sum_{a|d,n \text{ s.t. }\frac{d}{a}|n} \frac{\gcd(n,d)}{a} T_{n/a}.$$
\end{expl}

\begin{rem}
Noticing that given $d|n^2$, one has
$$a|d,n \text{ and }\frac{d}{a}|n\quad\Leftrightarrow\quad k=\frac{\gcd(n,d)}{a}|d,\frac{n^2}{d},$$
and we recover the expression from the first version of the paper.
\end{rem}

\subsubsection{M\"obius function for lattice of subgroups}

The counting functions we require involve the M\"obius function for finite abelian groups, which we briefly recall for the reader's convenience. We refer to \cite{rota1964foundations} for a more thorough introduction to the Theory of M\"obius functions. The case of interest was also handled in \cite{weisner1935abstract}.

If $G$ is a finite abelian group, we denote by $\P(G)$ the lattice of its subgroups, i.e. the finite set of its subgroups together with the inclusion order. A lattice is naturally endowed with a M\"obius function $\mu$, allowing, thanks to the so-called  to  M\"obius  inversion Theorem, to invert the summation operator. In this case, it means the function is defined by
$$\mu(H,H)=1 \text{ and }\forall H_1\subsetneq H_2,\ \sum_{H_1\subset K\subset H_2}\mu(K,H_2)=\sum_{H_1\subset K\subset H_2}\mu(H_1,K)=0.$$
Since any finite abelian group splits as a product of $p$-groups, we have
$$\P(G)\simeq\prod\P(G_p),$$
where $G_p$ is the $p$-part of $G$. The M\"obius function of a product of lattices is the product of M\"obius function. Therefore, the M\"obius function is fully determined by the M\"obius function for a $p$-group, which is given by Weisner's theorem \cite{weisner1935abstract}: if $H\subset K$, we have $\mu(H,K)=\mu(K/H)$ where $\mu$ is defined on abelian $p$-groups by
$$\mu(G)=\left\{ \begin{array}{l}
    (-1)^kp^{\binom{k}{2}} \text{ if }G\simeq\ZZ_p^k,  \\
    0 \text{ else.} 
\end{array}\right.$$

In our case, we only care about the following values: $\mu(\{0\})=1$, $\mu(\ZZ_p)=-1$ and $\mu(\ZZ_p^2)=p$. 

\subsubsection{Refined count of pairs}

We consider the $\NN$-module $\PPP_2=\{(d,n) \text{ s.t. }d|n^2\}\subset\NN^2$, with norm given by the projection onto the second coordinate. For any $(d,n)\in\PPP_2$, we set:
$$\bfF(d,n)=\sum_{\substack{H\subset K \subset\ZZ_n^2 \\ [\ZZ_n^2:H]=d}} \mu(K/H)T_K.$$
This is a \textit{refined twisted} way of counting pairs $H\subset K$ with $H$ of given index. The twist is due to the presence of the M\"obius function, the refinement due to the presence of $T_K$. The group $SL_2(\ZZ)$ acts on the set of such pairs, from which we deduce that $\bfF(d,n)$ is $SL_2(\ZZ)$-invariant and thus in the span of $(T_k)_{k|n}$.

\begin{prop}
    We have
    $$\bfF(d,n) = \left\{\begin{array}{l}
        \sum_{k|d}\mu(\frac{d}{k})T_{n/k} \text{ if }d|n, \\
        0 \text{ else.}
    \end{array}\right. ,$$
    where $\mu$ is the standard integer M\"obius function, multiplicative with value $-1$ over primes, and $0$ if divisible by the square of a prime.
\end{prop}

\begin{proof}
    The decomposition in a product of $p$-groups implies that $\bfF$ is multiplicative:
$$\bfF(d,n) = \prod_p\bfF(p^{\nu_p(d)},p^{\nu_p(n)}).$$
Just as $\bfF$, the formula provided on the right-hand-side is multiplicative, so that it suffices to check the identity for powers of primes. To do so, we use the same strategy as in Example \ref{expl-computation-subgroup} and compute $\pi_*\bfF(p^\delta,p^n)$, where $\delta\leqslant 2n$.
    
    Let $H$ be the subgroup of $\ZZ_{p^n}^2$ spanned by the columns of $\left(\begin{smallmatrix}
        p^\alpha & 0 \\ p^{\beta-\nu}u & p^\beta
    \end{smallmatrix}\right)$, where $\alpha,\beta\leqslant n$, $\alpha+\beta=\delta$, $0\leqslant\nu\leqslant\beta$ and $0<u\leqslant p^\nu$ is invertible modulo $p$. The condition ``$j$ being $(n/a)$-torsion modulo $b$'' means that we also require $n-\alpha+\beta-\nu\geqslant\beta$, or equivalently $\nu\leqslant n-\alpha$.

    If $K\supset H$ is chosen such that $\mu(K/H)\neq 0$, then $K/H\simeq\ZZ_p^r$ for some $r$. In particular, $K\subset H[1/p]=\{ x\in \ZZ_{p^n}^2 \text{ s.t. }px\in H\}$, and $K$ is thus the preimage of $\widetilde{K}=K/H\subset H[1/p]/H$. Furthermore, assuming $\alpha>0$, we have
    $$\ZZ_{p^{n-\alpha}}=\pi(H) \; \subset\;  \pi(K)\; \subset\;  \pi(H[1/p])\; \subset\;  \{y\in\ZZ_{p^n}\text{ s.t. }py\in\pi(H)\}=\ZZ_{p^{n-\alpha+1}},$$
    where the inclusion from first to last has index $p$. Therefore, exactly one inclusion is strict, and in case $\pi(H)\subsetneq\pi(H[1/p])$, $\pi(K)$ is equal to one of them.
    
    For each $H$, we determine $H[1/p]$ and compute $\sum_{K\supset H}\pi_*T_K$.

    \begin{itemize}[leftmargin=0.6cm]
        \item If $\alpha=0$, we have $\pi(H)=\ZZ_{p^n}$ and $\pi_* T_K=T_{p^n}$ for any $K\supset H$. Therefore, one only needs to compute $\sum_{\widetilde{K}\subset H[1/p]/H}\mu(\widetilde{K})$, which is $0$, except when $H[1/p]/H$ is the trivial group (i.e. $\beta=0$ and thus $\delta=0$), where the value is $1$. From now on, we assume $\alpha>0$.
        \item If $0\leqslant \nu<\beta$, then $H[1/p]$ is spanned by $\left(\begin{smallmatrix}
        p^{\alpha-1} & 0 \\ p^{\beta-\nu-1}u & p^{\beta-1}
    \end{smallmatrix}\right)$, $H[1/p]/H\simeq\ZZ_p^2$ and $\pi(H)\subsetneq\pi(H[1/p])$. The possible $K$ are the preimages of the subgroups of $\ZZ_p^2$:
        \begin{itemize}[label=$\bullet$]
            \item $K=H$ with $\mu(H/H)=1$,
            \item $K=H[1/p]$ with $\mu(H[1/p]/H)=\mu(\ZZ_p^2)=p$,
            \item the $p+1$ cyclic subgroups for which $\mu(K/H)=\mu(\ZZ_p)=-1$.
        \end{itemize}
        Among the preimages of the $p+1$ cyclic subgroups, one projects to $\pi(H)$ and the $p$ others to $\pi(H[1/p])$. Therefore, one gets the refined count
    $$ \sum_{K\supset H}\pi_*T_K = \bar{T}_{p^{n-\alpha}} -(\bar{T}_{p^{n-\alpha}}+p\bar{T}_{p^{n-\alpha+1}}) + p\bar{T}_{p^{n-\alpha+1}} = 0.$$

    \item If $0<\nu=\beta$, then $H[1/p]$ is spanned by $\left(\begin{smallmatrix}
        p^{\alpha} & 0 \\ u & p^{\beta-1}
    \end{smallmatrix}\right)$ and $H[1/p]/H\simeq\ZZ_p$. We only have two subgroups $H$ and $H[1/p]$ which share the same projection unde $\pi$. The refined contribution is thus
    $$\sum_{K\supset H}\pi_*T_K = \bar{T}_{p^{n-\alpha}}-\bar{T}_{p^{n-\alpha}} = 0.$$
    
    \item If $\beta=\nu=0$, so that $H$ is spanned by $\left(\begin{smallmatrix}
        p^{\alpha} & 0 \\ 0 & 1
    \end{smallmatrix}\right)$, then $H[1/p]$ is spanned by $\left(\begin{smallmatrix}
        p^{\alpha-1} & 0 \\ 0 & 1
    \end{smallmatrix}\right)$ and $H[1/p]/H\simeq\ZZ_p$. The only subgroups are $H$ and $H[1/p]$, and the refined contribution is
    $$\sum_{K\supset H}\pi_*T_K = \bar{T}_{p^{n-\alpha}}-\bar{T}_{p^{n-\alpha+1}}.$$
    \end{itemize}
If $\delta=0$, yielding $\alpha=\beta=\nu=0$, the value is $\bar{T}_{p^n}$. Otherwise, the only non-zero contribution is obtained when $\beta=0$ which only happens when $\delta=\alpha\leqslant n$, and the total sum is therefore
$$\sum_{\substack{H\subset K\subset\ZZ_{p^n}^2 \\ [\ZZ_{p^n}^2:H]=p^\delta}} \pi_*T_K=\bar{T}_{p^{n-\delta}}-\bar{T}_{p^{n-\delta+1}} = \sum_{i=0}^\delta \mu(p^{\delta-i})\bar{T}_{p^{n-i}} .$$
\end{proof}

\begin{coro}
    The function $\bfF$ defined on $\PPP_2$ satisfies the $0$-MCF.
\end{coro}

\begin{proof}
    The primitive coefficient of $\bfF$ is obtained when $k=1$ in the sum, the value is $\mu(d)$ when $d|n$ and $0$ else. Assuming $d|n$, we have
    $$\sum_{k|d,n} \Prim_{n/k}\bfF\left(\frac{d}{k},\frac{n}{k}\right)T_{n/k} = \sum_{k|d,n}\mu(d/k)T_{n/k}=\bfF(d,n).$$
    If we do not have $d|n$, there is some prime number $p$ with $\nu_p(d)>\nu_p(n)$. In particular, for any $k|d,n$ we still have $\nu_p(d/k)>\nu_p(n/k)$ and thus both sides are $0$.
\end{proof}

\subsubsection{Refined count of marked pairs}
\label{sec:group-counting-Gomega}

We now introduce a second twisted counting function of pairs of subgroups $(H\subset K)$ that also satisfies the $0$-MCF. The pairs are ``marked'' by the choice of a morphism $\varphi\colon\ZZ_n\to\ZZ_n^2/K$.

\medskip

Let $\omega|n$ and $K\subset\ZZ_n^2$ be a subgroup. We consider the following $\omega$-twisted version of $T_K$:
$$T_K[\omega] = \frac{1}{[\ZZ_n^2:K]}\sum_{\varphi\colon\ZZ_n\to\ZZ_n^2/K}\frac{1}{|K|}\sum_{x\in q^{-1}(\varphi(n/\omega))}(x),$$
where $q\colon\ZZ_n^2\to\ZZ_n^2/K$ denotes the projection, meaning the second sum is over the elements in the $K$-class $\varphi(n/\omega)$. We extend the definition by $0$ in case $\omega$ does not divide $n$.

\begin{lem}
\label{lem-expr-ref-twisted-count-2}
    For any $\omega|n$ and $K\subset\ZZ_n^2$, we have
    $$T_K[\omega]=T_\omega T_K.$$
\end{lem}

\begin{proof}
    By choosing a lift of the generator $1$, each morphism $\varphi\colon\ZZ_n\to\ZZ_n^2/K$ has exactly $|K|$ lifts to a morphism $\psi\colon\ZZ_n\to\ZZ_n^2$. Therefore, using that $|K|[\ZZ_n^2:K]=n^2$, we can rewrite
    \begin{align*}
    T_K(\omega) = & \frac{1}{n^2}\frac{1}{|K|}\sum_{\substack{\psi\colon\ZZ_n\to\ZZ_n^2 \\ x\equiv \psi(n/\omega)\mod K}}(x) \\
    = & \frac{1}{n^2|K|}\sum_{\psi\colon\ZZ_n\to\ZZ_n^2}\sum_{z\in K}(\psi(n/\omega)+z) = \left(\frac{1}{n^2}\sum_{\psi\colon\ZZ_n\to\ZZ_n^2}(\psi(n/\omega))\right)\left(\frac{1}{|K|}\sum_{z\in K}(z)\right). \\
    \end{align*}
   The second factor is $T_K$. To conclude, one merely needs to notice that
   $$\frac{1}{n^2}\sum_{\psi\colon\ZZ_n\to\ZZ_n^2}(\psi(n/\omega)) = T_\omega,$$
   since $\psi(n/\omega)$ is always $\omega$-torsion, and for any $\tau$ of $\omega$-torsion, there are precisely $(n/\omega)^2$ morphisms $\psi\colon\ZZ_n\to\ZZ_n^2$ such that $\psi(n/\omega)=\tau$, obtained by mapping $1$ to a $(n/\omega)$-root of $\tau$ in $\ZZ_n^2$.
\end{proof}

We set the following function on $\PPP_2$:
$$\bfG_\omega(d,n) = \mathds{1}_{\omega|n}\cdot \sum_{\substack{H\subset K\subset\ZZ_n^2 \\ [\ZZ_n^2:H]=d}} \mu(K/H)T_K[\omega] = \mathds{1}_{\omega|n}\cdot T_\omega\cdot\bfF(d,n),$$
where the second equality comes from Lemma \ref{lem-expr-ref-twisted-count-2}.

\begin{coro}
    The function $\bfG_\omega$ on $\PPP_2$ satisfies the $0$-MCF.
\end{coro}

\begin{proof}
    The MCF follows from the MCF for $\bfF$ and abstract properties of group algebra functions. To prove the MCF, we may restrict to the orbit $\Omega=\NN\cdot\widetilde{x}$ of a primitive element $\widetilde{x}\in\PPP_2$. As a $\NN$-module, it is isomorphic to the $\NN$-module $(\NN,|\cdot|_r)$ where $r=|\widetilde{x}|$. By restriction, since $\widetilde{x}$ is primitive, $\bfF|_\Omega$ satisfies the $0$-MCF.

    Let $Z=\{x\in\Omega \text{ s.t. }\omega \mid |x| \}$ be the sub-$\NN$-module of $\Omega$ consisting of elements whose norm is divisible by $\omega$, and $\iota\colon Z\hookrightarrow\Omega$ be the inclusion, so that we have $\bfG_\omega|_\Omega(x) = \mathds{1}_{x\in Z}\cdot T_\omega\cdot\bfF|_\Omega(x)$.
    
    Under the isomorphism $\Omega\simeq (\NN,|\cdot|_r)$, by Gauss lemma, $Z$ is isomorphic to $\NN\cdot\frac{\omega}{\gcd(r,\omega)}\subset (\NN,|\cdot|_r)$: indeed, for $x=k\widetilde{x}$, $|x|=rk$ and
    $$\omega\mid rk \;\Leftrightarrow\; \frac{\omega}{\gcd(r,\omega)}\mid\frac{r}{\gcd(r,\omega)}k \;\Leftrightarrow \;\frac{\omega}{\gcd(r,\omega)}\mid k.$$

    Therefore, by Lemma \ref{lem:restriction-suborbit}, $T_{r\omega/\gcd(r,\omega)}\bfF|_Z$ satisfies the $0$-MCF on $Z$. Furthermore,
    $$T_{r\omega/\gcd(r,\omega)}\bfF|_\Omega = T_{\mathrm{lcm}(r,\omega)}\bfF|_\Omega = T_\omega T_r\bfF|_\Omega = T_\omega\bfF|_\Omega,$$
    where the last equality stems from $\bfF|_\Omega$ being invariant by $r$-torsion on $\Omega$. Last, by Lemma \ref{prop:push-forward}, the push-forward $\iota_*(T_\omega\bfF|_Z)$ satisfies the $0$-MCF. Since $\iota$ is an inclusion, the push-forward is actually the extension by $0$ outside of $Z$, and is therefore equal to $\bfG_\omega|_\Omega$, finishing the proof.
\end{proof}
\section{Reduced decomposition formula}\label{sec-reduced-decomposition}

In this section we state and prove the reduced decomposition formula for semi-stable degenerations of abelian surfaces, generalizing \cite[Section~4]{maulik2010curves}. 
Section \ref{secred:prelim} contains a description of  the one parameter degenerations to which the formula applies, and introduces all the necessary ingredients to  state the reduced decomposition formula.
Section \ref{secred:formulas} contains three expressions of the reduced degeneration formula: a version similar to \cite{kim2018degeneration}, a version with values in the group algebra, and a numerical version of the formula.
The remainder of the section  (\ref{sec-reduced-class} and \ref{secred:splitting}) is devoted to proofs  and may be skipped at first reading. 

We assume the reader to be familiar with the theory of logarithmic smooth curves \cite{KatoF} and logarithmic stable maps \cite{gross2013logarithmic, abramovich2014stable}. 


\subsection{Preliminaries to the reduced degeneration formula}
\label{secred:prelim}


\subsubsection{Families of abelian surfaces}\label{sec:familiesabelian}
Consider
$A\xrightarrow{p} T$  a normal crossing  degeneration of an abelian surface:
\begin{itemize}[leftmargin=0.5cm]
    \item $T$ is a smooth curve, which we endow with divisorial log structure at a finite number of points $t_1,\dots, t_k$;
    \item $A_t$ is a smooth abelian surface for $t\in T\setminus\{t_1,\dots, t_k\}$;
    \item for $t_i=(\operatorname{Spec}\CC,\NN)$ a point in the base with log structure, $A_{t_i}$ has normal crossing singularities and no triple points.
\end{itemize}
The total space $A$ is endowed with the natural divisorial log structure given by the singular fibers, making $p$ into a log smooth morphism.

\medskip

We describe the geometry of the family on a disc centered at one of the points $t_i,$ which we denote by $0.$
It is known \cite{persson1977degenerations} that  $A_0$ is either smooth or a cycle of elliptic ruled surfaces, namely
\[A_0=Y_1\cup Y_2\dots \cup Y_N,\]
with 
\[Y_i=\mathbb P(\mathcal O_E\oplus\mathcal L_i),\]
 and $Y_i, Y_{i+1}$ meet transversely along a smooth elliptic curve $E$ which is a section of both. 
We call $E_i^+,E_i^-$ the sections along which $Y_i$ meets $Y_{i+1}$ and $Y_{i-1}$ respectively. We thus have isomorphisms $E_i^+\simeq E_{i+1}^-$.
All the $\mathcal L_i$ are isomorphic and have degree zero; 
up to a deformation we can assume they are all trivial, which is the case of the families we construct, and therefore assume $Y_i=E\times\PP^1$. Under such assumption, we get isomorphisms $E_i^+\simeq E_i^-$. The composition of all such isomorphisms as we go around the cycle yields an automorphism of the elliptic fiber $E$, which is the translation by a certain torsion element $u\in E$ referred to as the \emph{monodromy}. We require $u$ to be torsion to ensure the existence of a polarization on the family.

We thus write $A(u)$ or $A_0(u)$ to emphasize the monodromy dependence of the family and of the central fiber from the monodromy.




\subsubsection{Families of elliptically fibered abelian surfaces with prescribed monodromy}
\label{sec-construction-abelian-families}
Let $B\in\BBB\cong \NN^2$ be a diagram degree as in Example~\ref{ex:Nmodules}, with $|B|$ defined to be the first projection of $B.$
Fix $\tau,\tau'\in \mathbb H$ in the upper half plane and set $\lambda=e^{2i\pi\tau}$. Consider the elliptic curves: $F=\CC/\gen{1,\tau}\simeq\CC^*/\gen{\lambda}, $ $E=\CC/\gen{1,\tau'}$ and $F_n=\CC/\gen{n,\tau}$, which is a degree $n$ covering of $F$. 

\medskip

Let $u=\frac{k+l\tau'}{|B|}\in\CC$ be a $|B|$-torsion element in $E$, and $\Lambda_u$ be the lattice in $\CC^2$ spanned by the columns of the following matrix:
$$\Omega_u=\begin{pmatrix}
    1 & 0 & \tau & 0 \\
    u & 1 & 0 & \tau'
\end{pmatrix}.$$
We denote by $A(u)$ the quotient of $\CC^2$ by this lattice. We have the following:
\begin{itemize}[leftmargin=0.5cm]
    \item $A(u)$ is a complex torus  (in fact an abelian surface, see below) with a map to $F=\CC/\gen{1,\tau}$ coming from the projection onto the first coordinate $\CC^2\to\CC$, and with fiber $E=\CC/\gen{1,\tau'}$, i.e. $A(u)$ is elliptically fibered;
    \item 
    $A(u)$ only depends on the class of $u$ in $E$, since $\binom{0}{1}$ and $\binom{0}{\tau'}\in\Lambda_u;$
    \item  $A(u)$ can also be thought as the quotient of $F_{|B|}\times E$ by the action of $\ZZ_{|B|}$ acting by translation via $(z_1,z_2)\mapsto (z_1+1,z_2+u)$.
    To see it, multiply the first column of $\Omega_u$ by $|B|$ and then subtract $k$ times the second column and $l$ times the last; this is not an integral change of base since the determinant is $|B|,$ showing  that  $F_{|B|}\times E$ is a cover of $A(u).$
\end{itemize}

\begin{lem}\label{lem:polarization-u}
    The complex torus  $A(u)$ is an abelian surface endowed with the following polarization:
    $$Q=\begin{pmatrix}
        0 & 0 & -|B| & 0 \\
        0 & 0 & k & -a \\
        |B|  & -k & 0 & -l \\
        0  & a & l & 0
    \end{pmatrix}.$$
    The self-intersection is $2a|B|$ and its divisibility is $\gcd(k,l,|B|,a)$.
\end{lem}

\begin{proof}
    Following \cite{griffiths2014principles}, a polarization on a complex torus is a skew-symmetric integer matrix $Q$ satisfying \textit{Riemann bilinear relations}: $\Omega Q^{-1}\Omega^\intercal=0$ and $-i\Omega Q^{-1}\overline{\Omega}^\intercal$ is hermitian positive definite. We easily check that the proposed matrix satisfies these conditions for $\Omega_u$.
    As explained in \cite[Lemma~2.5]{blomme2025short} such a polarization $Q$ is of type $(d_1,d_2)$ with $d_1= \gcd(k,l,|B|,a)=\gcd(a,|B|,\frac{|B|}{\text{ord}(u)})$ and $d_1d_2=a|B|.$
\end{proof}

\begin{Notation}
    We denote by $\beta(u)$ the polarization on $A(u)$ determined by $Q.$ Geometrically, $a$ is the intersection with a section and $\pi_*\beta(u)=|B|[F]$, where $\pi\colon A(u)\to F$ is the projection onto the base of the elliptic fibration, or alternatively $|B|=\beta(u)\cdot[E]$. The self-intersection of $\beta(u)$ and its divisibility are related to $a,|B|$ and $u$ as stated in the Lemma.
We also write $|\beta(u)|$ to mean $|B|$, or even just $|\beta|$ when the monodromy does not  play a role.
\end{Notation}

\begin{rem}\label{rem:phikdescription}
 In above notation, let $A(u)$ and $A(u')$  two elliptically fibered surfaces. The isomorphism $\varphi\colon H_*(A(u),\ZZ)\to H_*(A(u'),\ZZ)$ (inducing by Poincar\'e duality an isomorphism in cohomology and which preserves the intersection pairing) is determined by the shear transform from $H_1(A(u),\ZZ)\cong\Lambda_u$ to  $H_1(A(u'),\ZZ)\cong\Lambda_{u'}$ given by
 $$\begin{pmatrix}
        1 & 0 & 0 & 0 \\
        \frac{k'-k}{|B|} & 1 & 0 & 0 \\
        0  & 0 & 1 & 0 \\
       \frac{l'-l}{|B|}   & 0& 0 & 1
    \end{pmatrix}.$$
In particular, under this isomorphism the polarization $\beta(u)$ has image $\beta(u'),$ with same self-intersection and divisibility $\gcd(a, |B|,\frac{|B|}{\text{ord}(u')}).$ Taking $u'$ a primitive $|B|$-root we can always achieve divisibility $1.$
 
\end{rem}


\paragraph{{\bf{A one parameter semi-stable degeneration over a disc}}}
To construct a degeneration of $A(u)$ over a disc, 
we can proceed as in \cite[Section 2.4.1]{blomme2024bielliptic}. Let $\widetilde{\F}\to \mathbb A^1_t$  the partial compactification of $(\CC^*)_{z,t}^2$ such that $\widetilde{\F}_{t\neq 0}\simeq\CC^*$ and $\widetilde{\F}_0$ is an infinite chain of $\PP^1$. There is  an action of $\ZZ$  defined by $(z,t)\mapsto (\lambda t^Nz,t)$ (here $\lambda=e^{2\pi i \tau}$ with $\tau$ in the upper half plane) which is proper for $|\lambda t^N|<1$ and extends to the central fiber $t=0$. The quotient $\F$ is a family of elliptic curves over the disc $\mathbb D=\{|t|<|\lambda|^{-1/N}\}$ with $\F_t=\CC^*/\gen{\lambda t^N}$ degenerating to the singular elliptic curve which is the union of $N$ copies of $\PP^1$.

\begin{rem}
    Using the language of logarithmic geometry, the construction can be rephrased in terms of a choosing a subdivision giving a schematic model for a degeneration of elliptic curves to a logarithmic elliptic curve $\mathcal F_{\log}:=\mathbb G_{m,\log}/\ZZ$ for $\mathbb G_{m,\log}$ the logarithmic multiplicative group over $\mathbb A^1_t$ with its divisorial log structure.
\end{rem}

We can then define a degeneration of $A(u)$ by taking the quotient of $\widetilde{\F}\times E$ by the action
$$(z_1,t,z_2)\in\widetilde{\F}\times E\mapsto (\lambda t^Nz_1,t,z_2+u).$$
As before, the action is defined for $t\neq 0$ and extends over the central fiber.  
The quotient $A_t(u)\to\mathbb D$ is a semistable degeneration with $A_t(u)$ a smooth abelian surface, and a central fiber $A_0(u)$ which is the gluing of $N$ copies of $\PP^1\times E$, with \textit{monodromy} data of $u$.

\medskip

The isomorphism in cohomology of Remark~\ref{rem:phikdescription} extends to the one parameter degeneration and keep preserving the classes coming from the fiber $E.$

\medskip

\paragraph{{\bf{A one parameter semi-stable degeneration over a smooth proper curve}}}

With slightly more work, it is possible to construct families of abelian surfaces with prescribed monodromy over a proper smooth curve $T$ rather than a disc. The existence of such families will be relevant to apply the results of \cite[Section~2]{maulik2025gromov}, recalled and adapted to the reduced class case in Section~\ref{sec:vanishingcohomology}.

We sketch below the construction.
\begin{itemize}[leftmargin=0.6cm]
    \item Start with a semistable family of genus $1$ curves $F\to T$ with at least one nodal fiber $F_0$; the fibers are either smooth elliptic curves or a loop of $\PP^1$, the total space is smooth.
    \item Up to taking a ramified covering of the base, we can assume that there exists a global N\'eron model of the family \cite{bosch2012neron} and we can choose a primitive $|B|$-torsion section $\sigma$. 
    \item The section $\sigma$ acts by translation by the $|B|$-torsion element $\sigma_t$ of $F_t$ on the smooth fibers. On a singular fiber $F_0$, it acts by rotating the necklace of $\PP^1$, or each $\PP^1$. We can assume the necklace rotation is non-trivial on $F_0$.
    \item We can then proceed as before considering the product $F\times_T E$, for $E\to T$ a trivial family of smooth elliptic curves and quotient by the action translating by $\sigma$ on the first factor and by any chosen  $u$ of $|B|$-torsion in $E$ on the second coordinate.
\end{itemize}

\subsubsection{Tropical maps and rigidity}
Let $A_0(u)\to 0=(\operatorname{Spec}(\mathbb C),\mathbb N)$ be a singular fiber of a normal crossing degeneration as constructed above. We denote the  \emph{tropicalization} of $A_T(u)\to T$ (in the sense for example of \cite[Appendix~B]{gross2013logarithmic}) by \[\Sigma\xrightarrow{\mathfrak{p}}\Sigma_{T},\] 
where $\Sigma_T$ is the union of $k$-rays $\mathbb R_{\geq 0}$, one for each point of $t$ with log strucure.

Restricted to one of these rays, $\Sigma$ is the generalized cone complex obtained gluing $N$ copies of $\mathbb R^2_{\geq 0}$, corresponding to the $N$ joining divisors, along their rays, corresponding to the irreducible components $Y_i.$ 
Let us fix one such ray and let
 $\Sigma_0$ be the pre-image  $\mathfrak{p}^{-1}(1)$ of the lattice point $1\in \mathbb R_{\geq 0};$  $\Sigma_0$ is a loop with $N$ vertices $y_1,\dots y_N.$ Endowing $\Sigma_0$ with the natural orientation making the vertex labeling increasing; $\Sigma$ is simply the cone over this subdivided loop.

\begin{defi}\label{defi-degeneration graph}
A \emph{tropical map} to $\Sigma_0$ is a metric graph $\Gamma$ together with a map of complexes $\phi\colon\Gamma\to\Sigma_0$  with the following decorations:
\begin{enumerate}[leftmargin=0.5cm]
    \item  \emph{vertices} are decorated with: a genus $g_v\in\mathbb Z_{\geqslant 0}$ and a curve class $\beta_v\in H_2(Y_{\phi(v)},\ZZ)$,
     \item every \emph{edge} $e$ comes decorated with a weight  $w_e\in\ZZ_{>0}$, called \textit{slope} of $\phi$ along $e$ and making the map to $\Sigma_0$ \textit{harmonic} (i.e. balanced at each vertex).
     \item We have the compatibility $\beta_v\cdot [E]$ is equal to the degree of $\phi$ at $v$.
\end{enumerate}
The \emph{type} $\tau$  of a tropical map is the data above together with, for each $v\in \Gamma$ the polyhedron $P_{\phi(v)}\subseteq\Sigma_0$ where $v$ maps (edge or vertex).

\end{defi}

For a fixed type $\tau$, there exists a polyhedron $\sigma_\tau$ parametrizing tropical maps of fixed type; in fact, having fixed the type, we can variate the edge lengths of the bounded edges of $\Gamma$ (see for example \cite[Section~2.5]{abramovich2020decomposition}).
\begin{defi} \cite[Definition~3.6]{abramovich2020decomposition} \label{defi:rigidtropmap}
A tropical map is called \emph{rigid}  if $\sigma_\tau$ is a point.

This amounts (possibly after a subdivision of the target) to $\phi$ mapping vertices to vertices and having no edge of slope $0$.
\end{defi}

In  \cite[Section~5]{kim2018degeneration} rigid curves are called \emph{degeneration graphs} (see Figure~\ref{fig:degenerationgraph}). They are those maps for which $\phi$ maps vertices to vertices without contracted edges for the case of rigid target $\Sigma_0.$

\medskip

Taking the cone over $\sigma_{\tau}$ we get a ray which simply parametrizes global dilation along the dilation of $\Sigma_0$ (recall $\Sigma_0=\mathfrak p^{-1}(1)$).

\begin{rem}
The definition of rigid tropical map recalled above will be relevant in extending the statement and proof of the reduced degeneration formula for \emph{tropical models} (see Section~\ref{sec-tropical-models}) of the moduli space of maps to $A_0(u),$ which in turn are relevant to deal with vanishing cohomology in the degeneration formula (see \cite[Section~2]{maulik2025gromov} and Section~\ref{sec:vanishingcohomology}).
\end{rem}

We denote by $\sfE(\Gamma)$ (resp. $\sfH(\Gamma)$) the set of edges (resp. half-edges) of $\Gamma$. If $h=(e\vdash v)$ is a half-edge, then we set $w_h=\pm w_e,$ where the sign is positive if $e$ goes out of $v$ in the positive direction (with respect to the orientation on $\Sigma_0$ fixed above), and negative if it goes in the opposite direction. 
In particular, the balancing condition (i.e. the fact that $\phi$ is harmonic) means that at each vertex $\sum_{h\vdash v}w_h=0$. 

\medskip

For a fixed vertex $v$ we call the \textit{flow} through $v$ the degree of $\phi$ at $v$, i.e. the sum of positive outgoing weights. Since $\beta_v\cdot[E]$ is equal to the flow through $v$, $\beta_v$ is fully determined by its image $a_v$ after projecting in $H_2(E,\ZZ)\simeq\ZZ$.

\begin{defi}
Let $\phi\colon\Gamma\to\Sigma_0$ be a tropical map 
\begin{itemize}[leftmargin=0.5cm]
    \item The genus of $\Gamma$ is $b_1(\Gamma)+\sum_v g_v$.
    \item The bidegree of $\Gamma$ is $B=(d,a)$ where $a=\sum_v a_v$ is the sum of the degrees of vertices, and $d$ is the degree of the harmonic map to $\Sigma_0$, i.e. the sum of edge weights mapping to a given edge of $\Sigma_0$.
\end{itemize}
\end{defi}

The degree of the harmonic morphism, i.e. the first coordinate of the bidegree $B$, is denoted by $|B|$.  This is consistent with the notation fixed above, since $B\in\BBB$ and $|B|$ is its norm.

\begin{figure}\label{fig:degenerationgraph}
	\centering
		\begin{tikzpicture}
		[x={(-0.2cm,-0.4cm)}, y={(1cm,0cm)}, z={(0cm,1cm)}, 
     scale=3,
     fill opacity=0.80,
     color={gray},bottom color=white,top color=black]

 \tikzset{zxplane/.style={canvas is zx plane at y=#1,very thin}}
 \tikzset{yxplane/.style={canvas is yx plane at z=#1,very thin}}

\begin{scope}[canvas is yx plane at z=0]
     \draw (0,0) circle (0.5cm);
     \filldraw[black] (0.5,0) circle (1pt) node[right]{$y_1$};
    \filldraw[black] (0,0.5) circle (1pt) node[right]{$y_2$};
     \filldraw[black] (-0.3,-0.35) circle (1pt) node[right]{$y_3$};;
   \end{scope}

    \draw[->,thick]  (0,0,0.7)--(0,0,0.4);
    \node[right] at (0,0,0.55) {$\phi$};
\begin{scope}[canvas is yx plane at z=1]
     \draw (0,0) circle (0.5cm);
     \filldraw[black] (0.5,-0) circle (1pt) node[right]{$v_1$};
    \filldraw[black] (0,0.5) circle (1pt) node[right]{$v_2$};
      \filldraw[black] (-0.4,-0.75) circle (1pt) node[above]{$v_3'$};
     \filldraw[black] (-0.3,-0.35) circle (1pt) node[below]{$v_3$};
    
   \draw[yshift=0cm]  (0.5,0)  to [bend left=65,looseness=1.5] (-0.4,-0.75);
  \draw[yshift=0cm]  (0.5,0)  to [bend left=95,looseness=1.7] (-0.4,-0.75);
     \draw[yshift=0cm]  (-0.4,-0.75)  to [bend left=85,looseness=1.5] (0,0.55);

    \node[below] at (0.3,0.35) {\small{$8$}};
     \node[above] at (0,-0.15) {\small{$4$}};
 \node[above,right] at (0.25,-0.55) {\small{$2$}};
 \node[above,right] at (0.35,-0.95) {\small{$2$}};
  \node[left] at (-0.35,0) {\small{$4$}};
  \node[left] at (-0.77,0) {\small{$4$}};

   \end{scope}

\end{tikzpicture}
	\caption{Rigid tropical curve with $b_1(\Gamma)=3,\; \deg{\phi}=8.$ The edges are decorated with their weight.}
\label{fig:degenerationgraph}
\end{figure}



\subsubsection{Evaluation hyperplane and reduced diagonal}\label{sec:evhyperplane}

Given $\phi\colon\Gamma\to\Sigma_0$  a rigid tropical curve, for each vertex $v$ of $\Gamma$ we consider the moduli space of stable log-maps to the component $Y_{\phi(v)}=E\times\PP^1$
with discrete data prescribed by $\Gamma$ and its decoration:
$$\M_v = \M_{g,m_v}(Y_{\phi(v)}|D,a_v,\bfw_v),$$
where $m_v$ is the number of markings adjacent to $v$, $D=E^+_{\phi(v)}+E^-_{\phi(v)}$ is the boundary divisor defining the log structure on $Y_{\phi(v)}$, and $\bfw_v=(w_h)_{h\vdash v}$ is the signed collections of weights adjacent to $v$. 

This comes equipped with its natural evaluation map:
$$\ev_v\colon\M_v\longrightarrow \prod_{h\vdash v}E,$$
where the product is indexed by half-edges adjacent to $v$.

\medskip

We denote by $\Delta_u$ the diagonal inclusion, with a subscript $u$ to remember that $A_0(u)$ has a prescribed monodromy datum $u$:
$$\Delta_u\colon E^{\sfE(\Gamma)}\to E^{\sfH(\Gamma)}.$$

Concretely, this is defined as follows: labeling the edges from $1$ to $N$, the last one relating $y_N$ to $y_1$, we have that:

\[\Delta_u(x_1,x_2,\dots,x_N)=(x_N+u,x_1,x_1,x_2,x_2,\dots,x_{N-1},x_{N-1},x_N),\]
where the half-edges are ordered by vertex of adjacency. Every entry is repeated, only the first one comes with a $u$-shift due to the monodromy coming back from $N$ to $1$.

\begin{defi}
\label{defi-evaluation-hyperplane}
Let $\phi\colon\Gamma\to\Sigma_0$ be a rigid tropical curve. We define the \emph{evaluation hyperplane} $\mathbf{H}_{\Gamma}\subseteq E^{\sfH(\Gamma)}$ to be

\[\mathbf{H}_{\Gamma}=\left\{(x_{h})\in E^{\sfH(\Gamma)}\; |\; \sum_{h} w_h x_{h} \equiv 0\in E\right\}.\]
\end{defi}

We observe some geometric properties of the evaluation hyperplane.

\begin{lem}\label{lem:evaluationhyperplane}
Let $\delta=\gcd(w_e)$ be the gcd of edges weights. We have the following:
\begin{enumerate}[leftmargin=0.5cm]
    \item 
    The evaluation hyperplane is disconnected with precisely $\delta^2$ connected components. These are the pre-images of the surjective morphism
    \[\mathbf{H}_{\Gamma}\xrightarrow{\sum_h\frac{w_h}{\delta}x_h}E[\delta]\subseteq E,\]
    where $E[\delta]$ are the $\delta$-torsion points of $E$.
    
    \item The evaluation map $\ev\colon\prod_v \M_v\to E^{\sfH(\Gamma)}$ factors trough $\mathbf{H}_{\Gamma}$.
    
    \item Let $u$ be a $|B|$-torsion element. The image of the diagonal $\Delta_u(E^{\sfE(\Gamma)})$ is contained in the connected component of  $\mathbf{H}_{\Gamma}$ given by the pre-image of the $\delta$-torsion element $\frac{|B|}{\delta}u$.
\end{enumerate}
\end{lem}

\begin{proof}
    \begin{enumerate}[leftmargin=0.5cm]
        \item The hyperplane is defined by an equation in $\ZZ^{\sfH(\Gamma)}\otimes E$. So after a suitable change coordinates in $\ZZ^{\sfH(\Gamma)}$, if $\delta=\gcd(w_e)$, the equation is equivalent to $\delta x_1\equiv 0$. Therefore, there are as many component as $\delta$-torsion elements in $E$, i.e. $\delta^2$. This also ensures the same map divided by $\delta$ maps surjectively to $\delta$-torsion elements.
        \item It suffice to rewrite the sum grouping half-edges per adjacent vertex:
        \[\sum_{h} w_h x_h =\sum_v\sum_{h\vdash v} w_h x_{h}. \]
        Then, it is well-known that for a morphism to $Y_{\phi(v)}=E\times\PP^1$, $\sum_{h\vdash v} w_h x_{h}=0$. This can for instance be obtained by pushing-forward to $E$ the principal divisor $\sum w_h(p_h)$ on $C$.
        \item It is sufficient to just rewrite the sum grouping half-edges per adjacent edge and using that 
         for each edge $e=(h,h')$, we have $w_h+w_{h'}=0$. Therefore, on the diagonal, we have that $w_hx_h+w_{h'}x_{h'}$ is $0$ if $\phi(e)$ is not $e_{N\to 1}$, and $w_e u$ if it is. In the end, we thus get $|B|\cdot u$ where $|B|$ is the degree of the harmonic map to $\Sigma_0$.
         
         Taking $u$ to be $|B|$-torsion, the image lies in $\mathbf{H}_\Gamma$. The connected component in which the diagonal embeds is obtained by the same computation with coefficients divided by $\delta$, yielding $\sum \frac{w_h}{\delta}x_h=\frac{|B|}{\delta}u$.
    \end{enumerate}
\end{proof}

\begin{defi}
\label{defi-relevant-hyperplane}
  We denote by $\mathbf{H}_{\Gamma}^u$ the connected component of the evaluation hyperplane containing the image of the diagonal and call it the \emph{u-relevant hyperplane}, and we denote by $\Delta_u^\mathrm{red}$ the diagonal map with values in the latter
$$\Delta_u^\mathrm{red}\colon E^{\sfE(\Gamma)}\to\mathbf{H}_\Gamma^u\subset E^{\sfE(\Gamma)}.$$

\end{defi}

\subsubsection{Relation to correlated GW invariants for $E\times\PP^1$}


We now use the \emph{correlated Gromov-Witten setting} \cite{blomme2024correlated,blommecarocci2025DR} for each moduli space $\M_v$. Correlated invariants are the vertex contributions in the degeneration formula for the reduced class  $\red{\M(A_0(u))}$. Briefly, for each vertex $v$, the image of the evaluation map $\ev\colon\M_v\to\prod_{h\vdash v}E_h$ is contained in the evaluation hyperplane
\[\mathbf{H}_v=\left\{\sum_{h\vdash v}w_h x_h=0\right\}.\]
Taking $\delta_v=\gcd(w_h)_{h\vdash v}$, the map
$\sum_{h\vdash v}\frac{w_h}{\delta_v}x_h$ surjects onto the $\delta_v$-torsion elements $E[\delta_v]$, indexing connected components of $\mathbf{H}_v$. Composing both, we get a \textit{correlating} map
$$\kappa_v\colon\M_v\to E[\delta_v],$$
mapping a stable map $f_v\colon (C,(p_h))\to E\times\PP^1$ to the correlator $\theta_v=\sum_{h\vdash v}\frac{w_h}{\delta_v}\cdot \mathrm{pr}\circ f_v(p_h)$, which is a $\delta_v$-torsion element of $E$. 

\medskip

We get a splitting of $\mathbf{H}_v$, $\M_v$ and its virtual class, indexed by the \emph{correlators} $\theta_v\in E[\delta_v]$ (see\cite[Section~3]{blomme2024correlated}):
\[\mathbf{H}_v=\bigsqcup_{\delta_v\theta_v\equiv 0 } \mathbf{H}_{\theta_v}, \quad\M_v=\bigsqcup_{\delta_v\theta_v\equiv 0 } \M_v^{\theta_v}, \quad \text{ and }\quad \fvir{\M_v}{\delta_v} = \sum_{\delta_v\theta_v\equiv 0}  \vir{\M_v^{\theta_v}}\cdot (\theta_v),\]
which we consider as an element in the group algebra $\QQ[E]$ with cycle coefficients. 



\begin{lem}
    The evaluation map restricted to $\prod\M_v^{\theta_v}$ takes value in the component of $\mathbf{H}_\Gamma$ indexed by $\sum\frac{\delta_v}{\delta}\theta_v$. In particular, this is the same component of the image of the diagonal $\Delta_u$ if we have the compatibility relation
    \begin{equation}\label{eq:contributingtheta}
    \sum_v\frac{\delta_v}{\delta}\theta_v=\frac{|\beta|}{\delta}u.
    \end{equation}
\end{lem}

\begin{proof}
The evaluation map $\prod_v\M_v\rightarrow \mathbf{H}_{\Gamma}$ takes values in 
\[ \mathbf{H}_{\Gamma}^{\rm{small}}=\prod_v \mathbf{H}_v=\bigsqcup_{\underline{\theta}} \mathbf{H}_{\Gamma}^{\underline{\theta}} \quad \text{ where }\quad \mathbf{H}_{\Gamma}^{\underline{\theta}}=\prod_v \mathbf{H}_{\theta_v},\]
and is in particular contained in $\mathbf{H}_\Gamma$. To determine in which connected component of $\mathbf{H}_\Gamma$ a given $\mathbf{H}_{\Gamma}^{\underline{\theta}}$ lies, we just apply the morphism $\mathbf{H}_{\Gamma}\to E[\delta]$ given by $\sum_h\frac{w_h}{\delta}x_h$. Grouping  by vertex, we have
$$\sum_h\frac{w_h}{\delta}x_h = \sum_v\frac{\delta_v}{\delta}\sum_{h\vdash v}\frac{w_h}{\delta_v}x_h = \sum_v\frac{\delta_v}{\delta}\theta_v.$$
Using Lemma \ref{lem:evaluationhyperplane}(3), we thus see that $\mathbf{H}_{\Gamma}^{\underline{\theta}}$ lies in the $u$-relevant hyperplane $\mathbf{H}_{\Gamma}^u$ (and thus potentially intersects the image of the diagonal $\Delta_u^{\mathrm{red}}(E^{\sfE(\Gamma)})$) if and only if we have the advertised \textit{compatibility condition}.
\end{proof}



\subsection{Reduced decomposition formula}
\label{secred:formulas}

We are now able to state the reduced degeneration formula. We provide three versions:
\begin{itemize}[leftmargin=0.5cm]
    \item a version for a unique family $A(u)$ with monodromy $u$ (Theorem \ref{theo:degene-formula-reduced}),
    \item a version with coefficients in the group algebra $\G$ encompassing the reduced degeneration formula for all possible monodromies $u$ of $|B|$-torsion (Corollary \ref{coro:degen-gp-alg}),
    \item a numerical version (Corollary \ref{coro:decompo-num-vers}) after providing an (and in fact several) expressions for the Poincar\'e dual class to the reduced diagonal in $\mathbf{H}_\Gamma^u$.
\end{itemize}

Denote by $\M^u=\M_{g,n}(A_t(u),\beta(u))\to T$ the moduli space of log stable maps to $A_t(u)\to T$, which comes endowed with
$\red{\M^u}\in \mathrm{A}_{g+n+1}(\M^u,\QQ)$ (construction recalled below) and by 
$\M_0^u$ its fiber over $0=\operatorname{Spec}(\NN\to \CC)$, the moduli space of log stable map to $A_0(u)/0.$

\subsubsection{First version of decomposition formula for reduced class}

For each  rigid tropical curve $\phi\colon\Gamma\to\Sigma_0$
we have a (logarithmic) stratum $\M_{0,\Gamma}^u$ of $\M_0^u$ which comes with a natural map $\mu_\Gamma\colon\M_{0,\Gamma}\to\M_0^u$  which is (at least virtually) the inclusion of an irreducible component of $\M_0^u.$

Let $\widetilde{\Gamma}$ denote a rigid tropical curve together with a labeling of its edges, we denote by $\mu_{\widetilde{\Gamma}}\colon\M_{0,\widetilde{\Gamma}}^u\to\M_0^u$ the morphism obtained composing $\mu_\Gamma$ with the finite \'etale map of degree $|E(\Gamma)|!$ forgetting the edge labeling. 

Finally, we set $l_\Gamma=\mathrm{lcm}(w_e)$ denote the lowest common multiple of edge weights.

We consider the following fiber product:

\begin{equation}\label{eq:gluingdiagram}
\begin{tikzcd}
\bigodot_v \M_v\ar[r, "g"]\ar[d,"\prod_e \ev_e"] &\prod_v \M_v\ar[d,"\prod_{h} \ev_{h}"]\\
E^{\sfE(\Gamma)} \ar[r,"\Delta_u"] & E^{\sfH(\Gamma)}.
\end{tikzcd}
\end{equation}

As argued in \cite[Section~9]{kim2018degeneration}, we have an \'etale map $\varphi_{\Gamma}\colon\M_{\widetilde{\Gamma}}\to \bigodot_v \M_v$ of degree $\left(\prod_{e\in \sfE(\Gamma)} w_e\right)/l_{\Gamma}$. 
Preimages by $\varphi_\Gamma$ correspond to the possible way to garnish the map to $A_0(u)$ (obtained gluing to log stable maps at the vertices ) with a compatible log structure to make it into an honest log stable map to $A_0(u)$.

\begin{theo}
\label{theo:degene-formula-reduced}
    Let $A_0(u)\to\operatorname{Spec}(\mathbb N\to \CC)$ be the central fiber of a semi-stable degeneration of abelian surfaces with monodromy $u$ and let $\M_0^u$ the moduli space of log stable map to $A_0(u)$ considered before. We have:
    \[\red{\M_0^u}=\sum_{\widetilde{\Gamma}}\frac{l_{\Gamma}}{|\sfE(\Gamma)|!}\mu_{\widetilde{\Gamma},*}\varphi_{\Gamma}^*\left (\sum_{\underline{\theta}} \Delta_u^{\mathrm{red},!}(\prod_v \vir{\M_v^{\theta_v}}) \right)\]
    where the sum is over the $\underline{\theta}=(\theta_v)$ satisfying the compatibility condition \eqref{eq:contributingtheta}: $\sum_v\frac{\delta_v}{\delta_\Gamma}\theta_v=\frac{|B|}{\delta_\Gamma}u$, with $\delta_\Gamma$ the g.c.d. of edge weights of $\Gamma$.
\end{theo}

\begin{proof}
    The proof follows from: Proposition \ref{prop:splitting-diagrams}, which shows that the reduced class splits as a sum over rigid tropical curves of reduced classes of the components  $\M_{0,\Gamma}^u;$ Proposition \ref{prop:splitting-vertices} which shows how to relate the reduced virtual class of the component for a fixed $\Gamma$ to the class of the fiber product $\bigodot_v\M_v$. 
    Once this is done, Lemma \ref{lem:evaluationhyperplane} gives a way to describe the pre-image $\ev^{-1}(\mathbf{H}_{\Gamma}^u)\subseteq \prod_v\M_v$ in terms of the components $\prod\M_v^{\theta_v}$ for which 
    the compatibility condition \eqref{eq:contributingtheta} holds.
\end{proof}

\subsubsection{Decomposition formula with values in group algebra}
\label{sec:correlatorsindegeneration}
    
To prove the multiple cover formula, it is convenient to consider the invariants for $A(u)$ for all possible $u\in E[|B|]$ at once. 
This information is best encoded in terms of the group algebra of $E$. 
We thus introduce the \emph{correlated reduced class:}
$$[\![\M_0]\!]^{\mathrm{red}}=\sum_{|\beta|u=0}\red{\M(A_0(u))}\cdot(u),$$
which we consider as an element in the group algebra $\G=\QQ[E]$ with cycle coefficients. It becomes an honest element in the group algebra $\G$ after integrating against some cohomology classes.

The following version of the decomposition formula is obtained from Theorem \ref{theo:degene-formula-reduced} via some bookkeeping on the group algebra coefficients.

\begin{coro}
\label{coro:degen-gp-alg}
    The correlated reduced class $[\![\M_0]\!]$ gathering the reduced class of the elliptic fibrations for every $|B|$-torsion element satisfies the following decomposition formula:
    \[ [\![\M_0]\!]^\mathrm{red}=\sum_{\widetilde{\Gamma}}\frac{l_{\Gamma}}{|\sfE(\Gamma)|!}\mu_{\widetilde{\Gamma},*}\varphi_{\Gamma}^*\Delta^{\mathrm{red},!}\left( (|B|/\delta_\Gamma)^2\d{\frac{1}{|B|/\delta_\Gamma}}\left(\prod_v\fvir{\M_v}{\delta_\Gamma}\right) \right),\]
    where $\delta$ is the gcd of all edge weights of $\Gamma$, $\fvir{\M_v}{\delta}$ is the correlated class of $\M_v$ at level $\delta$, and $\Delta^{\mathrm{red},!}$ acts by $\Delta^{\mathrm{red},!}_u$ on the $(u)$-coefficient.
\end{coro}

\begin{proof}
The idea is to consider the full correlated virtual class of level $\delta_v$ for the moduli spaces $\M_v$ and then use the structure and the operators in the group algebra $\G$ to read of the compatibility relation (\ref{eq:contributingtheta}) and automatically sort out which correlated classes contribute to a given monodromy $u$. To do so, it suffices to consider the following product and expand:
  
$$\prod_v\m{\delta_v/\delta}\fvir{\M_v}{\delta_v} = \sum_{\underline{\theta}} \prod_v\vir{\M_v^{\theta_v}}\left(\sum_v\frac{\delta_v}{\delta}\theta_v\right),$$
where $\m{\delta_v/\delta}\fvir{\M_v}{\delta}$ is non other than the level $\delta$ correlated virtual class $\fvir{\M_v}{\delta}$.
The compatibility relation then states that to get the contribution for the central fiber $A_0(u)$, we have to select the $(\frac{|B|}{\delta}u)$-coefficient. 
The operator $(|B|/\delta)^2\d{\frac{1}{|B|/\delta}}$ maps a generator $(t)$ to the sum of its $(\frac{|B|}{\delta})$-roots, which are  $|B|$-torsion elements, and so after applying this operator, the coefficient of $(t)$ is now multiplied by $\sum_{\frac{|B|}{\delta}(s)=(t)} (s).$ So now  we can simply select the $(u)$-coefficient and get the stated result.

\end{proof}


\subsubsection{Numerical version of decomposition formula}
\label{sec:numericaldec}

To obtain an effective numerical version of the decomposition formula we provide a K\"unneth decomposition (in fact several) for the Poincar\'e dual class to the diagonal inside the relevant hyperplane. We need the following basic result of algebraic topology:


\begin{lem}\label{lem:factorization}
    Let $A\subset B\subset C$ be smooth compact submanifolds. Let $\mu_A,\mu_B\in H^*(C,\QQ)$ be the Poincar\'e dual classes to $[A]$ and $[B]$ respectively. Let $\eta_A\in H^*(B,\QQ)$ be the Poincar\'e dual class to $[A]$. Then we have:
    \begin{enumerate}[label=(\roman*)]
        \item If $\eta_A=j^*\widetilde{\eta}$ for some $\widetilde{\eta}$ with $j\colon B\hookrightarrow C$, then $\mu_A=\widetilde{\eta}\cup\mu_B$.
        \item Conversely, writing $\mu_A=\widetilde{\eta}\cup\mu_B$ for some $\widetilde{\eta}$, we have that $j^*\widetilde{\eta}=\eta_A$ if $j_*$ is injective.
    \end{enumerate}
\end{lem}

\begin{proof}
    This follows by push-pull and the definition of Poincar\'e duality. Let $i,j$ be the inclusions $A\hookrightarrow B\hookrightarrow C$. We have
    $$j_*i_*[A] = j_*\left( [B]\cap\eta_A \right) = j_*\left( [B]\cap j^*\widetilde{\eta} \right) = j_*[B]\cap\widetilde{\eta} = [C]\cap(\mu_B\cup\widetilde{\eta}).$$
    Conversely, we just go right to left and use the injectivity of $j^*$ to conclude.
\end{proof}

We apply the previous lemma to the case where
\begin{itemize}
    \item $A=E^{\sfE(\Gamma)}$ is the diagonal,
    \item $B=\mathbf{H}_{\Gamma}^u$ is the $u$-relevant hyperplane containing the diagonal,
    \item $C=E^{\sfH(\Gamma)}$ the ambient torus.
\end{itemize}
In particular, the injectivity assumption is satisfied since we have an inclusion of tori.

\medskip

Let $(1,\alpha,\beta,\pt=\alpha\beta)$ be a basis of $H^*(E,\ZZ)$. We denote with an index $h$ the part of the cohomology corresponding to the copy of $E$ indexed by $h\in\sfH(\Gamma)$, providing generators of $H^*(E^{\sfH(\Gamma)},\QQ)$.

\begin{prop}
    For any edge $e_0$ in $\Gamma$, we have the following decomposition of Poincar\'e dual class to the reduced diagonal in the $u$-relevant hyperplane:
    $$\Delta_u^\mathrm{red} = \frac{\delta^2}{w_{e_0}^2}\prod_{\substack{e=(h,h') \\ e\neq e_0}} (\alpha_h-\alpha_{h'})(\beta_h-\beta_{h'}).$$
\end{prop}

\begin{proof}
By Lemma \ref{lem:factorization}, the Poincar\'e dual class to the diagonal embedded via $\Delta^{\mathrm{red}}_u$ inside $\mathbf{H}^u_\Gamma$ may be obtained as restriction of a class on $E^{\sfH(\Gamma)}$ by factoring out the Poincar\'e dual class to $\mathbf{H}^u_\Gamma$ from the usual Poincar\'e dual class of the diagonal, for which we have the K\"unneth decomposition formula:
$$\Delta_u = \prod_{e=(h,h')} (\alpha_h-\alpha_{h'})(\beta_h-\beta_{h'}).$$
This is easily seen since for an elliptic curve, the diagonal class is the pull-back of the (Poincar\'e dual) of a point class via the map $(x,y)\in E^2\mapsto x-y\in E$.


We saw before that the $u$-relevant hyperplane $\mathbf{H}_{\Gamma}^u$ is the pre-image of the $\delta$-torsion point $\frac{|B|}{\delta}u$ along the morphism $E^{\sfH(\Gamma)}\to E$ given by $\sum_h \frac{w_h}{\delta}$.
Therefore, the class Poincar\'e dual class to $\mathbf{H}^u_\Gamma$ is similarly
$$\mu_{\mathbf{H}^u_\Gamma} = \left(\sum_h \frac{w_h}{\delta}\alpha_h\right)\left(\sum_h \frac{w_h}{\delta}\beta_h\right) = \left(\sum_e \frac{w_e}{\delta}(\alpha_h-\alpha_{h'})\right)\left(\sum_e \frac{w_e}{\delta}(\beta_h-\beta_{h'})\right).$$

Let us fix $e_0$  an edge of $\Gamma$. We claim that  following equality holds:
$$\prod_{e} (\alpha_h-\alpha_{h'})(\beta_h-\beta_{h'}) = \frac{\delta^2}{w_{e_0}^2}\left(\sum_e \frac{w_e}{\delta}(\alpha_h-\alpha_{h'})\right)\left(\sum_e \frac{w_e}{\delta}(\beta_h-\beta_{h'})\right)\cup\prod_{e\neq e_0} (\alpha_h-\alpha_{h'})(\beta_h-\beta_{h'}).$$

To see that, expand the sum on the right hand side. As $(\alpha_h-\alpha_{h'})^2=(\beta_h-\beta_{h'})^2=0$, due to the product over $e\neq e_0$, most terms are $0$. The only terms contributing non trivially to the cup product are obtained taking $e$ equal to $e_0$ in both sums. For the latter, the coefficients $w_{e_0}/\delta$ simplify. This finishes the proof by application of Lemma \ref{lem:factorization}.
\end{proof}

Concretely, the class can be seen as follows: one should forget about the diagonal condition at exactly one of the edges, up to a coefficient depending on the weight of the edge: $(\frac{\delta}{w_{e_0}})^2$.

\begin{rem}
    To some extent, the discrepancy of $1$ in the dimension count for the (reduced) dimension of the moduli of stable maps states that, when using the decomposition formula, the last gluing is automatic (up to torsion). This is exactly what the above class provides: one should not impose the gluing at one of the edges (as it is automatic).
\end{rem}

Given constraints $\gamma\in H^*(A_0(u)^m,\mathbb Q)$   (by abuse of notation we keep denoting by $\gamma$ the constraint obtained via the isomorphism described in Remark~\ref{rem:phikdescription}  in $H^*(A_0(u'),\mathbb Q)$)  and $\alpha\in H^*(\overline{\M}_{g,m},\mathbb Q)$, we denote by
$$\gen{\gen{\alpha;\gamma}}_{g,B} = \sum_{|B|u\equiv 0} \gen{\alpha;\gamma}_{g,\beta(u)}^{A(u)}\cdot (u),$$
the group algebra element where the $(u)$-coefficient is the reduced Gromov-Witten invariant for the elliptically fibered abelian surface $A(u)$. 

The invariant encompasses all possible values of $u$.
From the explicit description of Section~\ref{sec-construction-abelian-families} we know 
that all  the polarization  $\beta(u)$ have the same self-intersection but different divisibilities depending on the order of $u$.

\begin{coro}
\label{coro:decompo-num-vers}
    Let $\gamma_\in H^*(A_0^m,\mathbb Q)$  and $\alpha\in H^*(\overline{\M}_{g,m},\mathbb Q)$. For each graph $\Gamma$ in the decomposition formula, choose an edge $e_0$ and set $\Delta_{e_0}^\mathrm{red}=\prod_{e\neq e_0} (\alpha_h-\alpha_{h'})(\beta_h-\beta_{h'})$. We have
\begin{equation}   \label{eq-numerical-decomposition-formula}
 \gen{\gen{\alpha;\gamma}}_{g,\beta} = \sum_{\Gamma}\frac{\prod_{e}w_e}{|\operatorname{Aut}(\Gamma)|}\frac{|B|^2}{w_{e_0}^2}\d{\frac{1}{|B|/\delta}}\left (\int_{\prod\fvir{\M_v}{\delta}} \prod \ev_v^*\gamma_{v}\cup \mathrm{ft}^*\alpha\cup\ev^*\Delta_{e_0}^\mathrm{red}\right)
\end{equation}

 where $\gamma_v$ is the restriction of $\prod_{i=1}^{m_v}\gamma_i$ to $Y_{\phi}(v)^{m_v}.$
 
\end{coro}

\begin{rem}
The decomposition formula as stated here only allows to recover the GW-invariants of $A$ for those cohomological insertions extending to the total space of the degeneration family. In Section~\ref{sec:vanishingcohomology}, adapting the arguments from \cite[Section~2]{maulik2025gromov} to the reduced class, we prove that the degeneration formula actually allows for a full recovery of the GW-theory of an abelian surface, at the cost of considering \emph{exotic insertions} in the GW-theory of $E\times\mathbb P^1$.
\end{rem}

\subsection{Reduced classes for log smooth degenerations of abelian surfaces}
\label{sec-reduced-class}
We review how to construct a reduced virtual class for moduli spaces of (log) stable maps to a normal crossing  degeneration $A\xrightarrow{p} T$ of abelian surfaces. For the moment being we suppress the $u$ keeping track of the monodromy from the notation.
We also refer to \cite{kiem2013localizing,maulik2010curves,kool2014reduced,bryan2018curve}.

We write $\M:=\mathcal M_{g,m}(A/T,\beta)$ for the moduli space of (basic \cite{gross2013logarithmic}) logarithmic stable maps to $A\slash T$.
It is known by \cite{gross2013logarithmic}  that $\M$ is a Deligne-Mumford stack which is proper over $T,$ and comes endowed with a fine and saturated log structure $M_{\M}$.



We refer for example  to \cite{gross2013logarithmic,ranganathanGWexpansion, MaulikRanganathan2025, kim2018degeneration} for an explicit description of the logarithmic structure on $\M;$ the description it is in terms of tropical maps to $\Sigma\to\Sigma_T$ and their moduli space.

For $C\to S$  a log smooth curve source of a log stable map to $A/T$ we consider the $m$-markings schematic, i.e. they do not carry log structure.

\medskip


\subsubsection{Logarithmic 1-forms on abelian surfaces}

We denote by $\Omega ^{\text{log}}_{A/T}$ the sheaf of logarithmic forms in the sense of \cite[3.3]{gross2011tropicalbook} (see also \cite{olsson2005logarithmic} for the general definition of logarithmic cotangent complex of an algebraic stack endowed with a logarithmic structure).
Since $A\to T$ is log smooth, $\Omega ^{\text{log}}_{A/T}$ is a locally free sheaf.
 The log cotangent bundle enjoys the same functorial property of the usual cotangent bundle (see for example \cite[IV.1.1]{ogus}); in particular 
\[\Omega ^{\text{log}}_{A/T}\rvert_{A_0}\cong \Omega ^{\text{log}}_{A_0\slash 0},\]
where $0$ denotes the standard log point $\operatorname{Spec}(\mathbb N\to\mathbb C).$

\medskip

The sheaf $ \Omega ^{\text{log}}_{A_0\slash 0}$ can be explicitly described as  the sub-bundle of $\bigoplus_{i=1}^N \Omega_{Y_i}(E_i^- + E_i^+)$ of forms with at most logarithmic poles along the gluing divisor which have opposite residues  and agree when restricted to the divisor.

It follows from the explicit description that
the relative log canonical $K^{\text{log}}_{A/T}=\bigwedge^2\Omega ^{\text{log}}_{A/T}$ is trivial. Let $\omega$ a generator, local on the base $B$, then it induces an isomorphism between $\Omega_{A/T}^{\log}$ and its dual, the \textit{log-tangent bundle}:
\[\theta\colon T ^{\text{log}}_{A/T}\to \Omega ^{\text{log}}_{A/T}.\]

\subsubsection{Standard obstruction theory}

\paragraph{\textbf{Relative Obstruction theories}}

The moduli space $\M$
 comes equipped with a perfect obstruction theory:
\[(R\pi_*F^* T ^{\text{log}}_{A/T})^{\vee}\cong R\pi_*R\mathcal Hom(F^*\Omega ^{\text{log}}_{A/T},\mathcal O)^\vee,\]
relative to the log smooth stack of maps from log smooth curves  to the relative Artin fan $\mathcal A_{A}/T$ :
 \[\mathfrak M:=\mathfrak M_{g,m}(\mathcal A_{T}/T).\] 
We refer the reader to \cite[Section~3]{abramovich2020decomposition} and references therein for any  detail about relative Artin fans, maps to them and the relative obstruction theory we just described.
The obstruction theory for $\M$  had been constructed before, see for example
\cite{abramovich2014stable,li2001stable,gross2013logarithmic,chen2013degeneration,kim2010logarithmic}.

\medskip

\paragraph{\textbf{Absolute obstruction theory}}

The stack $\mathfrak M$ is equipped with a natural morphism to the  Olsson stack \cite{olsson2003logarithmic}  $\mathcal Log_T$ parametrizing logarithmic schemes with a log map to $T.$

Indeed, such a morphism is equivalent to the data of a log structure on $\mathfrak M.$ 
 Moreover, $\mathfrak M\to\mathcal Log_T$ is smooth (see for example \cite[Proposition~3.3]{abramovich2020decomposition}).  We can then consider a perfect obstruction theory for $\M$ relative to $\mathcal Log_T.$ This is done by a standard procedure (performed for example  in \cite[Section~4]{maulik2010curves} and \cite[Section~4.3]{chen2013degeneration}).
From the morphism 
\[F^* \Omega ^{\text{log}}_{A/T}\to \Omega ^{\text{log}}_{\mathcal C\slash \mathcal M}\to \Omega ^{\text{log}}_{\mathcal C\slash \mathcal M}(x_1+\dots +x_m) ,\] 
we get 
\begin{equation}\label{eq:cone}
    R\pi_*R\mathcal Hom(\Omega ^{\text{log}}_{\mathcal C\slash \mathcal M}(x_1+\dots +x_m),\mathcal O)\to R\pi_*R\mathcal Hom(F^*\Omega ^{\text{log}}_{A/T},\mathcal O).
\end{equation}

By  Kato \cite{KatoF}, 
we have that $R\pi_*R\mathcal Hom(\Omega ^{\text{log}}_{\mathcal C\slash \mathcal M}(x_1+\dots +x_m) ,\mathcal O)\cong\mathbb T_{\mathfrak M\slash\mathcal Log_T}[-1].$ 
The smoothness of $\mathfrak M\to\mathcal Log_T$  implies that taking the cone over \eqref{eq:cone}  we get  the dual of a perfect obstruction theory for $\M$ relative to  $\mathcal Log_T$ :
\[R\pi_*R\mathcal Hom([F^*\Omega ^{\text{log}}_{A/T}\to \Omega ^{\text{log}}_{\mathcal C\slash \mathcal M}(x_1+\dots +x_m)] ,\mathcal O).
\]
Here $[F^*\Omega ^{\text{log}}_{A/T}\to \Omega ^{\text{log}}_{\mathcal C\slash \mathcal M}(x_1+\dots +x_m)]$ is supported in degree $-1,0$.

By base change along the smooth morphism $T\to\mathcal Log_T$ \cite[Section~4]{chen2013degeneration} we can also think of this as the obstruction theory relative to $T.$ 

We denote by
\[\mathcal Ob_{\mathcal M\slash \mathfrak M}=\mathcal H^1(R\pi_*F^*T_{A/T}^{\log}),\;\;\;\mathcal Ob_{\mathcal M\slash T}=\mathcal H^1(R\pi_*R\mathcal Hom([F^*\Omega ^{\text{log}}_{A/T}\to \Omega ^{\text{log}}_{\mathcal C\slash \mathcal M}(x_1+\dots +x_m)] ,\mathcal O))\] the obstruction sheaves relative to $\mathfrak M$ and to $T.$

\begin{rem}
    There is no obstruction to deform a map from a connected curve together with a deformation of the fibers of $A\to T,$  the obstruction sheaf $\mathcal Ob_{\mathcal M\slash T}$ coincides with the absolute obstruction sheaf $\mathcal Ob_{\M}.$

\end{rem}

\begin{rem}
    In \cite[Section~4]{maulik2010curves} the authors work with logarithmic stable maps to expansion $\M^{\exp}$, as in \cite{kim2010logarithmic,ranganathanGWexpansion}. This moduli space is virtually birational to the one we consider by \cite{abramovich2014comparison}, and has a perfect obstruction theory \[R\pi_*R\mathcal Hom([F^*\Omega ^{\text{log}}_{\mathfrak A\slash \mathfrak T}\to \Omega ^{\text{log}}_{\mathcal C\slash \mathcal M}(x_1+\dots +x_m)] ,\mathcal O).
\]
where $\mathfrak A\slash \mathfrak T$ is the universal expansion, relative to the smooth Artin stack of expansion $\mathfrak T,$ which is logarithmically \'etale on $\mathcal Log_T$ (see for example \cite{abramovich2014comparison} and \cite{MaulikRanganathan2025} for a general treatment of the stach of expansions.)

This is a particular example of a \emph{tropical model}/\emph{virtual birational model} (See Section~\ref{sec:virtual_birational model}) for the moduli space of log maps to $A/T.$  
\end{rem}


\subsubsection{Cosection}
Following the steps of \cite{kiem2013localizing,maulik2010curves,kool2014reduced} we construct a cosection (i.e. a map to $\O_\M$) for the obstruction sheaf  $\mathcal Ob_{\mathcal M\slash T}\cong \mathcal Ob_{\M}$ by pull-back. We then show that it is surjective.




\begin{lem}
    We have the following:
    \begin{enumerate}[leftmargin=0.6cm]
        \item There exists a natural cosection $\sigma_\mathrm{rel}\colon \O b_{\M/\mathfrak{M}}\to\O_\M$.
        \item The cosection from (1) descends to a cosection $\sigma\colon\O b_\M\to\O_\M$.
        \item Both cosections are surjective.
    \end{enumerate}
\end{lem}

\begin{proof}
\begin{enumerate}[leftmargin=0.6cm]
\item Consider the composition, where first and last map are actually isomorphisms:
\[F^*T ^{\text{log}}_{A/T}\xrightarrow{F^*\theta} F^*\Omega ^{\text{log}}_{A/T}\to \Omega ^{\text{log}}_{\mathcal C\slash \mathcal M}\to\omega_{\mathcal C\slash \mathcal M} . \]

Taking the first cohomology of $R\pi_*$, where we recall that $\pi$
is the projection $\mathcal{C}\xrightarrow{\pi}\M,$ we thus obtain the cosection

\[\sigma_{\text{rel}}\colon \mathcal Ob_{\mathcal M\slash\mathfrak M}\cong R^1\pi_* F^*T ^{\text{log}}_{A/T}\to R^1\pi_* F^*\Omega ^{\text{log}}_{A/T}\to R^1\pi_* \Omega ^{\text{log}}_{\mathcal C\slash \mathcal M}\to R^1\pi_* \omega_{\mathcal C\slash \mathcal M} \cong\mathcal O_{\M}.\]

\item The exact triangle
\[\mathbb T_{\mathfrak M\slash T}[-1]\to R\pi_* F^*T ^{\text{log}}_{A/T} \to R\pi_* R\mathcal Hom([F^*\Omega^{\text{log}}_{A/T}\to\Omega^{\log}_{ \mathcal C\slash\M}(x_1+\dots +x_m)],\mathcal O_{\mathcal C}) \]
induces a long exact sequence of cohomology sheaves, whose last three terms are
\[R^1\pi_* R\mathcal Hom(\Omega^{\log}_{ \mathcal C\slash\M}(x_1+\dots +x_m),\mathcal O_{\mathcal C})\to \ R^1\pi_* F^*T ^{\text{log}}_{A/T} \to\mathcal Ob_{\M\slash T}\to 0.\]

Then, to show that $\sigma_{\text{rel}}$ descends to $\sigma\colon \mathcal Ob_{\M\slash T}\to\mathcal O_{\M}$ (and thus automatically to $\mathcal Ob_{\M}$ by the comment made before) it is sufficient to check that the following composition vanishes:

\[R^1\pi_* R\mathcal Hom(\Omega^{\log}_{ \mathcal C\slash\M}(x_1+\dots +x_m),\mathcal O_{\mathcal C})\to \ R^1\pi_* F^*T ^{\text{log}}_{A/T} \xrightarrow{\sigma_{\text{rel}}}\mathcal O_{\M}.\]

Since $\theta$ is the isomorphism induced by the global section $\omega$ of the log canonical,  this composition comes, as also explained in \cite{kool2014reduced}, by dualizing and pushing forward the following composition:
\[F^* K^{\log}_{A/T}\to F^*\Omega^{\log}_{A/T}\otimes F^*\Omega^{\log}_{A/T}\to \Omega^{\log}_{\mathcal C\slash\M}\otimes \Omega^{\log}_{\mathcal C\slash\M}\to \Omega^{\log}_{\mathcal C\slash\M}(x_1+\dots +x_m)\otimes \omega_{\mathcal C\slash\M}.\]
As $\omega$ is alternating and $\bigwedge^2\Omega^{\log}_{C\slash\M}=0$, we can conclude that the composition vanishes and the cosection descends to a cosection of the absolute obstruction sheaf:
\[\sigma\colon \mathcal Ob_{\M/B}=\mathcal Ob_{\M}\to\O_\M.\]

\item Surjectivity of $\sigma_{\text{rel}}$ implies the surjectivity of $\sigma$. This in turn reduces to show that for a log-curve $f\colon C\to A_t$, the following composition is surjective:
\[H^1(C,f^*T^{\log}_{A_t})\xrightarrow{H^1(f^*\theta)}H^1(C,f^*\Omega^{\log}_{A_t})\to H^1(C,\Omega^{\log}_{C})\to H^1(C,\omega_C).\]
Notice that the first and last maps are isomorphisms, and that $H^1(C,\omega_C)\cong\mathbb C$. 

We only need to argue that the composition is non zero. When $A_t$ is smooth, this is the argument of \cite[Proposition~6.4]{kiem2013localizing}. In our case, this can be done as well looking at the Serre dual:
\[H^0(C,\mathcal O_C)\to H^0(C,T_C^{\log}\otimes\omega_C)\to H^0(C,f^*T^{\log}_{A_t}\otimes\omega_C)\xrightarrow{\sim} H^0(C,f^*\Omega^{\log}_{A_t}\otimes\omega_C).\]
Since $\mathcal O_C$ is globally generated, it is enough to prove that the map of sheaves
 \[\mathcal \O_C\hookrightarrow T_C^{\log}\otimes\omega_C\xrightarrow{\mathrm{d}^{\log}f} f^*T^{\log}_{A_t}\otimes\omega_C\]
 is nontrivial,
 where the first map is induced by dualizing $\Omega_C^{\log}\to \omega_C.$ But the map ${\mathrm{d}^{\log}f\otimes\omega_C}$ is non trivial: on the locus in $C$ that maps to $A_t^{\text{sm}}$ and $f$ is not ramified on the image; at the nodes mapping to the $A_t^{\text{sing}}$ if the intersection of the image curve and the singular divisor is dimensionally transverse. This locus is non empty.
\end{enumerate}
\end{proof}

\subsubsection{Reduced class}
Let us denote by
\[ \mathbf{E}:=\left (R\pi_*R\mathcal Hom([F^*\Omega ^{\text{log}}_{A/T}\to \Omega^{\text{log}}_{\mathcal C\slash \mathcal M}(x_1+\dots +x_m)] ,\mathcal O)\right)^\vee.\]
The cosection provides a map:
\[\mathbf{E}^\vee\to \mathcal H^1(\mathbf{E}^\vee)=\mathcal Ob_{\M\slash T}\to\mathcal O_{\M};\]
by the discussion above this is also a cosection of the absolute obstruction sheaf $Ob_{\M}.$

Here the maps are to be intended as maps of complexes, thus the obstruction sheaf and the cosection are in degree $1$. 

Dualizing the latter, we get a non trivial map $\mathcal O_{\M}[1]\to \mathbf{E}$. We define $\mathbf{E}_{\text{red}}$, 
to be the cone over the morphism:
\begin{equation}\label{eq:cosectriangle}
    \mathcal O_{\M}[1]\to \mathbf{E}\to \mathbf{E}_{\text{red}}.
\end{equation}


 We can similarly  define $\mathbf{E}^{\text{rel}}_{\text{red}}$ where $\mathbf{E}^{\text{rel}}=(R\pi_*F^*T^{\log}_{A/T})^\vee$ is the obstruction theory relative to $\mathfrak {M}$, since the cosection is in fact naturally defined on the relative obstruction theory.


\begin{lem}\label{lem:relativeconered}
The cone $\mathfrak C_{\M\slash\mathcal Log_T}$ admits a closed embedding in $h^1/h^0(\mathbf{E}^{\vee}_{\mathrm{red}}).$ 
    The cone $\mathfrak C_{\M\slash\mathfrak M}$ admits a closed embedding in $h^1/h^0(\mathbf{E}^{\mathrm{rel},\vee}_{\mathrm{red}}).$ 
\end{lem}

\begin{proof}
The cosections on $\mathcal Ob_{\M/\mathfrak M}$ and 
$\mathcal Ob_{\M/\mathcal Log_T}$ are compatible with the absolute cosection by construction, hence the statement simply follows from the absolute case (\cite[Section~4]{kiem2013localizing}) and the functoriality of the $h^1/h^0(-)$ construction. 
Further details are spelled out in \cite[Lemma~2.8]{chang2017torus}.
\end{proof}
 Applying the virtual Gysin pull-back \cite{manolache2012virtual} for $\mathfrak C_{\M\slash\mathcal Log_T} \hookrightarrow h^1/h^0(\mathbf{E}^{\vee}_{\mathrm{red}})$ or equivalently for
    $\mathfrak C_{\M\slash\mathfrak M} \hookrightarrow h^1/h^0(\mathbf{E}^{\mathrm{rel}}),$  
    we  get a reduced virtual class
\[[\M]^{\text{red}}\in A_{g+m+1}(\M,\mathbb Q).\] 



For any $t\in T$ we have a cartesian square:
\bcd
\M_t:=\M_{g,m}(A_t,\beta)\ar[r,"\iota"]\ar[d, "p_t"] & \M\ar[d]\\
t\ar[r,"\iota_t"] & T
\ecd
where for $t\neq t_i$,  $\M_t$ in the usual Kontsevich moduli space of stable maps to $A_t$ and for $t=t_i,$ with $t_i$ denoting the standard log point, $\M_{t_i}$ is the moduli space of log stable maps to the semi-stable  surface $A_{t_i}.$
The obstruction theory and the cosection are compatible with restriction to $t\in T,$ and consequently also the reduced class is. Meaning the following:

\begin{prop}
   Let $\iota_t\colon t\to T$ the (regular) embedding of a point of the base of the degeneration and let $\iota_t^!\colon A_*(\M)\to A_{*-1}(\M_t)$ the Gysin pull-back, then
   \[\iota_t^!\red{\M}=\red{\M_t}.\]
   In particular, the reduced Gromov-Witten invariants are independent from $t.$
\end{prop}

\begin{proof}
    It is proved in \cite[Proposition~5.10]{behrend1997intrinsic} (see also \cite[Remark~4.10]{manolache2012virtual}) that the second line of the following diagram is a compatible triple of obstruction theories in the sense of \cite[Definition~5.8]{behrend1997intrinsic}:
\bcd
\mathcal \iota^*\mathcal O_{\M}[1]\ar[r,"="]\ar[d] & \mathcal O_{\M_t}[1]\ar[d]\\
\iota^*\mathbf{E}\ar[r] & \mathbf{E}_t\ar[r] & p_t^* (I_t\slash I_t^2)[1],
\ecd
where $I_t$ denotes the ideal sheaf of $t$ inside the base $T$. Taking the cones along the vertical maps induced by the surjective cosections and looking at the associated vector bundle stacks we get a short exact sequence of vector bundle stacks
\[p_t^*N_t\to h^1/h^0(\iota^*\mathbf{E}^{\vee,\text{red}})\to h^1/h^0(\mathbf{E}^{\vee,\text{red}}_t). \]
Since by the results of \cite{kiem2013localizing} the respective cones admit embeddings in these kernel vector bundles, the compatibility proved in \cite[Proposition~5.10]{behrend1997intrinsic} ensures the desired functoriality property for the reduced class.

By Lemma~\ref{lem:relativeconered} the same argument applies verbatim if we consider the obstruction theories and reduced obstruction theories relative to $\mathfrak M$.
\end{proof}

\begin{rem}
In the proof of the previous Proposition we didn't say that taking the cone over the vertical morphisms we obtain a new compatible triple of \emph{reduced} obstruction theories, but we only made a statement about the cycles. As we learned from \cite[Appendix~A]{maulik2010curves}  the first statement is in general stronger and more delicate to prove. (It requires a careful study of the deformation theory of the moduli spaces.) In the case at hand, we know by work of Kool-Thomas \cite{kool2014reduced}, based on classical results about the semi-regularity map \cite{bloch1972semi,ran1999semiregularity,manetti2007lie}, that $\mathbf E_t^{\text{red}}$ is indeed a reduced obstruction theory for $t\neq 0$. With some more deformation theory work to handle the case relative to a divisor the same can probably be said for the general case of log smooth families, as already predicted in \cite{maulik2010curves}.
\end{rem}

\subsection{Splitting of the reduced class}
\label{secred:splitting}

From now on we only look at $\M_0$, the moduli space of log stable maps to the central fiber $A_0,$ together with its reduced virtual cycle defined from the closed embedding:
\begin{equation}\label{eq:reduced_obstruction}
   \mathfrak{C}_{\M_0/\mathfrak M_0}\hookrightarrow h^1/h^0(\mathbf{E}^{\vee,\text{red}}_{\text{rel},0}):=\mathfrak E^{\text{red}}.
\end{equation}

Our goal is now to show that $\red{\M_0}$ can be written as a sum of virtual classes that admit a decomposition as product over vertices. To do so, we adapt the steps of \cite{kim2018degeneration} following \cite[Section~4]{maulik2010curves}. 
A degeneration formula for relative/log Gromov-Witten invariance in the semi-stable degeneration case already appeared in several other papers, where the different flavors of the theory (relative, expanded , etc) are considered \cite{li2002degeneration,chen2013degeneration,ranganathanGWexpansion,MaulikRanganathan2025,abramovich2020decomposition}.

\subsubsection{From tropical  to reduced virtual decomposition}

\begin{prop}(\cite[Section~8]{kim2018degeneration}, \cite[Section~3]{abramovich2020decomposition})
\label{prop:splitting-diagrams}
    Let $\phi\colon\Gamma\to\Sigma_0$ be a rigid tropical curve and let $l_\Gamma=\mathrm{lcm}(w_e)_{e\in \sfE(\Gamma)}$. We have a splitting
    $$\red{\M_0}=\sum_{\Gamma} l_{\Gamma}\cdot \mu_{\Gamma,*} \red{\M_{\Gamma}},$$
    where $\mu_\Gamma\colon\M_{0,\Gamma}\to\M_0$.
\end{prop}

\begin{proof}
Given a rigid tropical curve $\phi\colon\Gamma\to\Sigma_0$ we denote by $\bar{\Gamma}$ the 
data obtained forgetting the curve class decoration $\beta_v \in H_2(Y_{\phi(v)},\ZZ).$ 
Let consider the log smooth stack $\mathfrak M_0=\mathfrak M\times_{\mathcal L og_T}\mathcal Log_0.$
As explained in \cite[Section~8]{kim2018degeneration} and also in \cite[Section~3]{abramovich2020decomposition} we have a cartesian diagram:

\bcd
\bigsqcup_{\Gamma} \M_{0,\Gamma}\ar[r]\ar[d,"\sqcup \mu_{\Gamma}"] & \bigsqcup_{\bar{\Gamma}} \mathfrak M_{0,\bar{\Gamma}}\ar[d,"\sqcup \mu_{\mathfrak M_{\bar\Gamma}}"]\\
\M_0\ar[r,"p"] &\mathfrak M_0
\ecd
where the $\bar{\Gamma}$ are is one to one correspondence with the rays $\rho_{\bar{\Gamma}}$ of the cone stack $\Sigma(\mathfrak M)$ (alias with the vertices of the polyhedral stack $\Sigma(\mathfrak M_0)$ ), associated to the logarithmic smooth stack $\mathfrak M$  ($\mathfrak M_0$ respectively) as explained in \cite[Section~3]{abramovich2020decomposition}.
We notice that for each $\bar{\Gamma}$ we have 
\[p^{-1}\mathfrak M_{0,\bar{\Gamma}}=\bigsqcup_{\Gamma} \M_{0,\Gamma}\]
as $\Gamma$ runs over all the rigid tropical curves which reduce to $\bar{\Gamma}$ after forgetting the curve class decoration. 

\medskip

It is proved in \cite[Lemma~8.2]{kim2018degeneration} that for $l_{\Gamma}$ the lowest common multiple of edge weights (this is the so called base order in\cite[Section~4]{abramovich2020decomposition})  we have
\[[\mathfrak M_0]=\sum_{\bar{\Gamma}} l_{\Gamma}\cdot(\mu_{\mathfrak M_{\bar\Gamma}})_*[\mathfrak M_{0,\bar{\Gamma}}].\]

Now, Manolache's virtual pull-back construction \cite[Construction~3.6]{manolache2012virtual} gives us a \emph{reduced} virtual pull-back

\[p^!_{\mathfrak E^{\text{red}}}\colon A_*(\mathfrak M_0)\to A_*(\mathfrak{C}_{\M_0/\mathfrak M_0})\to A_*(\mathfrak E^{\text{red}})\to A_{*-\text{rk}(\mathfrak E^{\text{red}})}(\M_0) \]

such that, as discussed before, $[\M_0]^{\text{red}}=p^!_{\mathfrak E^{\text{red}}}([\mathfrak M_0]).$ Combined with the splitting of $[\mathfrak M_0]$ just discussed we get:

\begin{equation}\label{eq:reducedsplitting}
    [\M_0]^{\text{red}}=\sum_{\bar{\Gamma}} l_{\Gamma}\cdot p^!_{\mathfrak E^{\text{red}}}(\mu_{\mathfrak M_{\bar\Gamma},*}[\mathfrak M_{0,\bar{\Gamma}}])= \sum_{\Gamma} l_{\Gamma}\cdot \mu_{\Gamma,*}( p^!_{\mathfrak E^{\text{red}}_{\Gamma}} [\mathfrak M_{0,\bar{\Gamma}}])=\sum_{\Gamma} l_{\Gamma}\cdot \mu_{\Gamma,*} [\M_{\Gamma}]^{\text{red}},
\end{equation}

where $\mathfrak E^{\text{red}}_{\Gamma}$ is obtained via the $h^1/h^0$ construction from $\mu_{\Gamma}^*\mathbf{E}_{\text{rel},0}^{\text{red},\vee}$. 
The fact that $\mathfrak C_{\M_{\Gamma}/\mathfrak M_{0,\bar{\Gamma}} }$ is contained in $\mathfrak E^{\text{red}}_{\Gamma}$ simply follows from the fact that, as explained for example in \cite[Proposition~2.18]{manolache2012virtual}
there is a closed embedding of $ \mathfrak C_{\M_{\Gamma}/\mathfrak M_{0,\bar{\Gamma}} }\hookrightarrow \mu_{\Gamma}^*\mathfrak{C}_{\M_0/\mathfrak M_0}.$
The second equality in \eqref{eq:reducedsplitting} follows from the properties of the virtual pull-back  \cite[Theorem~4.1.(3)]{manolache2012virtual} since $\mu_{\bar{\Gamma}}$ are finite by \cite[Lemma~7.3,7.4]{kim2018degeneration}.
\end{proof}


\subsubsection{Marking the bounded edges}
The next step is to relate $\red{\M_{0,\Gamma}}$ to a virtual class that splits a product over the vertices of $\Gamma.$ 

This relation is via the Gysin pull-back along a suitable diagonal map. 

To do so, we need to label the edges of $\Gamma$, which is the content of the present section. From now on we fix one rigid tropical curve $\phi\colon\Gamma\to\Sigma_0$. We denote by $\widetilde{\Gamma}$ a rigid tropical curve which is furthermore endowed with a labeling of its bounded edges.

We denote by $\mathfrak M_{\widetilde{\bar{\Gamma}}}$ the stack of log smooth curves over $0$ whose dual graph admits an edge contraction to $\bar{\Gamma}$ and whose splitting nodes have been marked. 

This comes with an étale map $g_{\bar{\Gamma}}\colon \mathfrak M_{\widetilde{\bar{\Gamma}}}\to \mathfrak M_{\bar{\Gamma}}$ of finite degree $|\sfE(\Gamma)|!$. We thus have a cartesian diagram:

\bcd
\bigsqcup\M_{\widetilde{\Gamma}}\ar[r,"\sqcup G_{\Gamma}"]\ar[d] & \bigsqcup \M_{\Gamma}\ar[d]\\
\mathfrak M_{\widetilde{\bar{\Gamma}}}\ar[r, "g_{\bar{\Gamma}}"] & \mathfrak M_{\bar{\Gamma}}.
\ecd

To simplify the notation we denote by 

\[\mu_{\widetilde{\Gamma}}\colon \M_{\widetilde{\Gamma}}\xrightarrow{G_{\Gamma}}\M_{\Gamma}\xrightarrow{\mu_{\Gamma}}\M_0\]

The properties of pull-back of cones along cartesian diagrams and functoriality of the virtual pull-back ensure that:

\begin{itemize}
    \item The cone $\mathfrak C_{\M_{\widetilde{\Gamma}}/\mathfrak M_{\widetilde{\bar{\Gamma}}}}$ is contained in the vector bundle stack  $\mathfrak E^{\text{red}}_{\widetilde{\Gamma}}$, obtained via the $h^1/h^0$ construction from 
    \[\mathbf{E}_{\widetilde{\Gamma}}^{\text{red},\vee}:=\mu_{\widetilde{\Gamma}}^*\mathbf{E}_{\text{rel},0}^{\text{red},\vee}.\] 
    Notice that we have an exact triangle
    \begin{equation}\label{eq:reducedongamma}
        \mathcal O_{\M_{\widetilde{\Gamma}}}[1]\to\mu_{\widetilde{\Gamma}}^*\mathbf{E}_{\text{rel},0}\to \mu_{\widetilde{\Gamma}}^*\mathbf{E}_{\text{rel},0}^{\text{red}},
    \end{equation}
    induced by pull-back from the exact triangle \eqref{eq:cosectriangle} and by the surjective cosection.
   \item We have that
   \[  [\M_0]^{\text{red}}= \sum_{\widetilde{\Gamma}} \frac{l_{\Gamma}}{|\sfE(\Gamma)|!}\cdot (\mu_{\widetilde{\Gamma}})_* [\M_{\widetilde{\Gamma}}]^{\text{red}}\]
\end{itemize}

\subsubsection{Splitting as product on the vertices }\label{sec:splittingover vertices}
We consider the fiber product $\bigodot_v\M_v$ defined as follows:
\begin{equation}\label{eq:gluingdiagram}
\begin{tikzcd}
\bigodot_v \M_v\ar[r, "g"]\ar[d,"\prod_e \ev_e"] &\prod_v \M_v\ar[d,"\prod_{h} \ev_{h}"]\\
E^{\sfE(\Gamma)} \ar[r,"\Delta_u"] & E^{\sfH(\Gamma)}
\end{tikzcd}
\end{equation}

\begin{itemize}[leftmargin=0.6cm]
    \item We denoted by $\sfH(\Gamma)$ the set of \emph{half edges} of $\Gamma;$ each edge $e$ correspond to two half-edges $(h,h')$ each one rooted at one of the two end points of $e$.
    \item The map $\Delta_u$ is the diagonal inclusion, with a subscript $u$ to remember that $A_0$ has a prescribed monodromy datum $u$.
    \item Furthermore we denoted by 
 \[\M_v:=\M_{g_v, m_v}(Y_{\phi(v)}|D; \beta_v;(w_h)_{h\vdash v})\]
 the moduli space of log stable map to the component $Y_{\phi(v)}$
of the central fiber $A_0$, with log structure given by the divisor $D=E_{\phi(v)}^++E_{\phi(v)}^-,$ and numerical data determined by the fixed rigid tropical curve; we denoted by $n_v$ the number of relative markings (nodes to be), i.e. $n_v=|\{ h\in \sfH(\Gamma), h\vdash v\}|$ is the valency of $v$ and $m_v$ the number of markings attached to $v$.
\end{itemize}

In order to ease notation, in what follows we will simply write $\ev$ for the evaluation map. As argued in \cite[Section~9]{kim2018degeneration} we have an \'etale map $\varphi_{\Gamma}\colon\M_{\widetilde{\Gamma}}\to \bigodot_v \M_v$ of degree $\left(\prod_{e\in \sfE(\Gamma)} w_e\right)/l_{\Gamma}$. Preimages by $\varphi_\Gamma$ correspond to lifts of the maps in $\bigodot_v \M_v$ to a log stable map to $A_0$. 

\medskip

Standard results on relative/ log stable maps \cite{li2001stable,gross2013logarithmic} tell us that
$\M_v$ has a perfect obstruction theory relative to the stack $\mathfrak M^{\log}_{g_v,n_v+m_v}$  of log smooth curves with $n_v$ logarithmic and $m_v$ schematic markings given by
by $(R\pi_{v,*}F_v^* T^{\log}_{Y_{\phi(v)}})^{\vee}$ with $\pi_v$ and $F_v$ denoting the projection from the universal curve and the universal log stable map.

\medskip

The main point of \cite[Section~9]{kim2018degeneration} is to prove that 
\[[\M_{\widetilde{\Gamma}}]^{\text{vir}}=\varphi_{\Gamma}^*(\Delta_u^!\prod_v [\M_v]^{\text{vir}}),\]
where $[\M_{\widetilde{\Gamma}}]^{\text{vir}}$ is the virtual class described above, defined from the the pull-back of the usual relative obstruction theory $(R\pi_*F^*T^{\log} A_0)^\vee$ along $\mu_{\widetilde{\Gamma}}$ and 
$ [\M_v]^{\text{vir}}$ is the one just recalled. 

For dimensional reasons, this formula cannot not be applied to get a decomposition of the reduced class; it would in fact simply yield the usual virtual class which produces zero invariants. It can however be adapted to obtain a splitting of the \emph{reduced class} $[\M_{\widetilde{\Gamma}}]^{\text{red}}$ defined before. This is the content of Proposition \ref{prop:splitting-vertices}. 


Let us denote by $\Delta_u^{\text{red}}\colon E^{\sfE(\Gamma)}\to \mathbf{H}^u_{\Gamma} $ the inclusion of the diagonal into the evaluation hyperplane; we have that $\Delta_u=\iota_{\mathbf{H}}\circ\Delta_u^{\text{red}}$ where $\iota_{\mathbf{H}}$ is the embedding of the relevant hyperplane, and thus we have a short exact sequence of bundles on $E^{\sfE(\Gamma)}$, which are actually here all trivial:
    \begin{equation}\label{eq:normalcones}
    0\to N_{\Delta_u^{\text{red}}}\to N_{\Delta_u}\to N_{\iota_{\mathbf{H}}}\rvert_{E^{\sfE(\Gamma)}}\to 0.
    \end{equation}

\begin{prop}
\label{prop:splitting-vertices}
    We have an equality of cycles in $A_*(\M_{\widetilde{\Gamma}}):$
    \[\red{\M_{\widetilde{\Gamma}}}=\varphi_{\Gamma}^*(\Delta_u^{\mathrm{red},!}\prod_v \vir{\M_v}).\]
\end{prop}

\begin{proof}
Reasoning as in \cite[Section~9]{kim2018degeneration} but replacing $\Delta_u$ with $\Delta_u^{\text{red}}$ in the cartesian product \eqref{eq:gluingdiagram} we obtain a commutating diagram of exact triangles, parallel to \cite[Diagram~9.10]{kim2018degeneration} 

\bcd
 \ev^*\mathbb L_{\Delta_u^{\text{red}}}[-1]\ar[r]\ar[d] & g^*\boxtimes_v (R\pi_{v,*}f_v^* T^{\log}_{Y_{\phi(v)}})^{\vee}\ar[d]\ar[r] &\mathbf{E}^{\mathrm{red}}_{\bigodot}\ar[d]\\
\mathbb L_g[-1]\ar[r] & g^*\mathbb L_{\prod_v \M_v/\mathfrak M_v}\ar[r] & \mathbb L_{\bigodot_v \M_v /\prod_v\mathfrak M_v}.
\ecd

Since $\Delta_u^{\text{red}}$ is still a regular embedding, the left vertical arrow in the diagram induces a surjection on the $h^0$ and thus (precisely as also argued in \cite{kim2018degeneration}) since the central arrow is a perfect obstruction theory,  using the two Four Lemmas, also the right arrow gives a perfect obstruction theory.

Furthermore, by the functoriality properties of the virtual class (see for example \cite[Proposition~5.10]{behrend1997intrinsic}) the virtual class on $\bigodot_v\M_v$ induced by $\mathbf{E}^{\mathrm{red}}_{\bigodot},$ wich we denote by $\red{\bigodot_v\M_v}$ satisfies:\[\red{\bigodot_v\M_v}=\Delta_u^{\text{red},!}\prod_v [\M_v]^{\text{vir}}.\]

Moreover, by standard results on diagrams of exact triangles in a triangulated category (see for example \cite[tag 05QN]{stacks-project}), one sees that the standard obstruction theory $\mathbf{E}_{\bigodot},$ given by the cone over 
\begin{equation}\label{eq:ebsfiberprod}
\ev^*\mathbb L_{\Delta_u}[-1]\to g^*\boxtimes_v (R\pi_{v,*}f_v^* T^{\log}_{Y_{\phi(v)}})^{\vee}
\end{equation}
and the reduced one differ by the trivial factor $\ev^*\mathbb L_{\iota_{\mathbf{H}}}=\ev^* N^{\vee}_{\mathbf{H}_{\Gamma}/\prod_h E_h}$, i.e. we have an exact triangle:

\[\ev^*\mathbb L_{\iota_{\mathbf{H}}}\cong\mathcal O_{\bigodot \M_v}[1]\to \mathbf{E}_{\bigodot}\to \mathbf{E}_{\bigodot}^{\mathrm{red}} .\]

Pulling this back along $\varphi_{\Gamma}$ and using the existence of the diagrams of exact triangles \cite[Diagram~(9.12)]{kim2018degeneration} we obtain the following:
\bcd
\varphi_{\Gamma}^*\mathbf{E}_{\bigodot}^{\mathrm{red}}\ar[r] & \mu_{\widetilde{\Gamma}}^*\mathbf{E}^{\mathrm{red}}_{\text{rel},0}\ar[r] &  \varphi_{\Gamma}^*(\bigoplus_{e\in \sfE(\Gamma)} \ev_e^* N^{\log}_{E\slash Y_0})^\vee[1]\\
\varphi_{\Gamma}^*\mathbf{E}_{\bigodot}\ar[r]\ar[d]\ar[u] & \mu_{\widetilde{\Gamma}}^*\mathbf{E}_{\text{rel},0}\ar[r]\ar[d]\ar[u] & \varphi_{\Gamma}^*(\bigoplus_{e\in \sfE(\Gamma)} \ev_e^* N^{\log}_{E\slash Y_0})^\vee[1]\ar[d]\ar[u,"="],\\
\mathbb L_{\bigodot_v \M_v /\prod_v\mathfrak M_v}\ar[r] & \mathbb L_{\M_{\widetilde{\Gamma}} /\mathfrak M_{\widetilde{\Gamma}}}\ar[r] & p_{\Gamma}^*\mathbb L_{\mathfrak M_{\widetilde{\Gamma}} /\prod_v\mathfrak M_v}
\ecd
where $p_{\Gamma}\colon \M_{\widetilde{\Gamma}}\to \mathfrak M_{\widetilde{\Gamma}}$ is the natural morphism forgetting the log map.
The exact triangle in the first line comes from the fact that we have a quasi-isomorphism 

\[\varphi^*_{\Gamma}\ev^*\mathbb L_{\iota_{\mathbf{H}}}[1]\cong \mu_{\widetilde{\Gamma}}^*(R^1\pi_{\Gamma,*}\omega_{\mathcal C\slash\M})^\vee[1]\]
commuting with the given maps to $\varphi_{\Gamma}^*\mathbf{E}_{\bigodot}$ and $\mu_{\widetilde{\Gamma}}^*\mathbf{E}_{\text{rel},0}.$

Once we verify this claim, which we do below, since we have already proved that the cones $\mathfrak C_{\M_{\widetilde{\Gamma}}\slash \mathfrak M_{\widetilde{\Gamma}}}$ and $\mathfrak C_{\bigodot_v\M_v\slash\prod_v \mathfrak M_v}$  are contained in the respective reduced vector bundle stacks, the functoriality of the virtual pull-back construction \cite[Section~4]{manolache2012virtual}
gives the statement.

We can check the isomorphism after dualizing and passing to closed points. Then we have to check that given $[f\colon C\to A_0]\in\M_{\widetilde{\Gamma}}$ inducing (abusing notation) $[f\colon C\to A_0]\in\bigodot\M_v$ we have an isomorphism
\[H^1(C,\omega_C)\to (N_{\mathbf{ H}_{\Gamma}/E^{\sfH(\Gamma)}})_{\ev(f)}\]
commuting with the maps from $H^1(f^*T^{\log}_{A_0/0})$ and $H^1(\mathbf{E}_{\bigodot}^\vee\rvert_{[f]}).$

We have a natural surjection: 
\begin{align}\label{eq:surj}
    T_{\ev(f)}(E^{\sfE(\Gamma)})&\twoheadrightarrow  (N_{\mathbf{ H}_{\Gamma}/E^{\sfH(\Gamma)}})_{\ev(f)}
\end{align}
 coming 
from the pull-back along $f$ of the short exact sequence \eqref{eq:normalcones}. By the explicit description of the evaluation hyperplane, the kernel of this surjection of vector space  is simply the hyperplane in $ T_{\ev(f)}(E^{\sfE(\Gamma)})$ cut out by the equation $\sum w_e\partial x_e=0.$

It follows from the definition of $\mathbf{E}_{\bigodot}$ as the cone over the morphism in \eqref{eq:ebsfiberprod} that \eqref{eq:surj} factors trough  $H^1(\mathbf{E}_{\bigodot}^\vee\rvert_{[f]}).$
Furthermore, tensoring the partial normalization sequence
\[0\to\mathcal O_C\to\nu_*\left(\bigoplus_v \mathcal O_{C_v}\right)\to\bigoplus_e\mathcal O_{q_e}\to 0\]
by $f^*T^{\log}_{A_0/0}$ we see that there is a natural map
\[T_{\ev(f)}(E^{\sfE(\Gamma)})=\bigoplus_e T_{f(q_e)}E_e\to \bigoplus_e f^*T^{\log}_{A_0/0}\rvert_{q_e}\to H^1(f^*T^{\log}_{A_0/0}).\]
The latter induces via the isomorphism $H^1(f^*T^{\log}_{A_0/0})\cong H^1(f^*\Omega^{\log}_{A_0/0})$ a surjective map to $H^1(C,\omega_C).$
The map can also be factored as:
\begin{equation}\label{eq:surj2}
    T_{\ev(f)}(E^{\sfE(\Gamma)})\to \bigoplus_e f^*T^{\log}_{A_0/0}\rvert_{q_e}\cong  \bigoplus_e f^*\Omega^{\log}_{A_0/0}\rvert_{q_e}\to \bigoplus_e \omega_C \rvert_{q_e}\twoheadrightarrow H^1(\omega_C).
\end{equation} 
We can now explicitly describe the map in local coordinates. 
Locally around a point of its singular locus, $A_0$ is cut out by an equation of the form $ts=0$ inside $k[x_e,t,s]$ where  $x_e$  is the coordinate on $E=\left\{t=s=0\right\}$ and $t,s$ are the coordinate on the $\mathbb P^1$-fibers. Then $f^*T^{\log}_{A_0/0}$ in a neighborhood of  $q_e$ by $\partial x_e,t\partial t$ and the $\partial x_e$ is in the image of the first map.
The explicit local expression of the isomorphism $f^*T^{\log}_{A_0/0}\rvert_{q_e}\cong f^*\Omega^{\log}_{A_0/0}\rvert_{q_e}$ via the logarithmic form $d^{\log}t\wedge d_{x_e}$ shows that this sill map via the composition to $w_e d^{\log}u$ where $d^{\log}u$ is the generator of $\omega_C\rvert_{q_e}.$
Thus the
kernels of the  maps \eqref{eq:surj} and \eqref{eq:surj2} are naturally  identified, proving the desired isomorphism.
\end{proof}

\section{Guideline, the MCF for point insertions}

Before proving the general statement, in order to illustrate the strategy exposed in the introduction, we show in this section how to obtain the multiple cover formula for points and $\lambda$-classes insertions using the reduced degeneration formula and the properties of the correlated invariants for $E\times\PP^1.$ For these special insertion we know the local contribution explicitly and the proof simplifies substantially.
\label{sec-MCF-point-case}

\subsection{Local setting and computation}

We consider the correlated GW-invariants of the log-CY $E\times\PP^1$ relative to its divisor $D=E\times\{0,\infty\}$. Let $\bfw$ be a ramification profile, i.e. $\sum_1^n w_i=0$, with $\gcd(w_i)=|\bfw|$, and $a\geqslant 1$. We denote by $\RRR$ (as in ``Ramification") the $\NN$-module of these choices:
$$\RRR = \NN\times\{(w_1,\cdots,w_n)\in\ZZ^n\backslash\{0\} \text{ s.t. }\sum w_i=0\}\subset\NN\times\ZZ^n.$$
The $\NN$-action scales each coordinate, and the norm is defined as the gcd of the $\bfw$ coordinates:
$$|(a,\bfw)| = \gcd(w_1,\dots,w_n).$$
To keep the notation light, we suppress the $a$ and simply write $|\bfw|$ for $|(a,\bfw)|.$

\medskip

Each $(a,\bfw)\in\RRR$ determines an effective homology class, namely $a[E]+(\sum_{w_i>0}w_i)[\PP^1]$ in $H_2(E\times\PP^1,\ZZ).$ We consider the correlated relative GW-invariant
$$\gen{\gen{\pt_0,1_{w_1},\pt_{w_2},\dots,\pt_{w_n}}}^{|\bfw|}_{1,a,\bfw} ,$$
where we fix point insertions for every adjacent point except one.

Suitably rephrasing the computations in \cite[Section~5]{blomme2024correlated}, one sees that the correlated GW-invariants are functions from the $\NN$-module $\RRR$ to $\G$ of diagonal type, namely, for fixed $(a,\bfw),$ these refined invariants take values  in the finite dimensional sub-algebra $\G_{|\bfw|}$. We have an explicit expression provided by \cite[Theorem 5.3]{blomme2024correlated}.

\begin{prop}\label{prop:vertex-contrib-EP1}
    We have the following equality in $\G_{|\bfw|}$:
    $$\gen{\gen{\pt_0,1_{w_1},\pt_{w_2},\dots,\pt_{w_n}}}^{|\bfw|}_{1,a,\bfw} =  w_1^2\cdot a^{n-1}\cdot \bsigma^{|\bfw|}(a).$$
\end{prop}

In particular, by Lemma \ref{lem:monomial-rule}, it satisfies the $(n+1)$-MCF.

\begin{rem}\label{rem:correlatedforE}
The expression provided in an earlier version of \cite{blomme2024correlated} uses an a priori different refinement of the sum of divisors function. It can be proved that they agree by showing that they satisfy the unrefinement relation and have the same $(0)$-coefficient.

Indeed, as proven in \cite[Lemma~5.15]{blomme2024correlated} for the case of $E\times\mathbb P^1$, the correlated GW-invariants at a given level depend on the correlator $\theta$ only through its order as a torsion element. Therefore, using the unrefinement relations, the full correlated invariant at level $\delta$ is determined inductively on $\delta$ by only knowing the coefficient for the $(0)$-correlator and the non-correlated invariant.
\end{rem}

\subsection{Multiple cover formula for correlated invariants of $E\times\PP^1$}
\label{sec:multiplicity-complex-case}

    \subsubsection{Statement and idea of proof}
We now consider the following correlated GW-invariants of $E\times\PP^1$:
$$\gen{\gen{\pt^{n+g-1}}}^{|\bfw|}_{g,a,\bfw},$$
generalizing the example above. This was computed in \cite[Section~6]{blomme2024correlated} applying the correlated refinement of the degeneration formula; we briefly recall the computation below. In particular we have that $\gen{\gen{\pt^{n+g-1}}}^{|\bfw|}_{g,a,\bfw}\in \G_{|\bfw|}.$
We can then consider the $|\bfw|$-primitive coefficient, given by
$$\Prim_{|\bfw|}\gen{\gen{\pt^{n+g-1}}}^{|\bfw|}_{g,a,\bfw} = |\bfw|^2\gen{\pt^{n+g-1}}^{{|\bfw|},\prim}_{g,a,\bfw}$$
where  $\prim$ means any correlator of order $|\bfw|$; this definition is well posed since the correlated invariants only depend on the order of the correlator $\theta\in\ E[|\bfw|]$ (see \cite[Section~5.4.1]{blomme2024correlated}). 
The ${|\bfw|}^2$ factor simply comes from $T_{|\bfw|}=\frac{1}{|\bfw|^2}\sum(\theta)$.

\begin{theo}\label{theo-MCF-EP1-points}
    The function $(a,\bfw)\mapsto\gen{\gen{\pt^{n+g-1}}}^{|\bfw|}_{g,a,\bfw}$ satisfies the $(2n+4g-4)$-MCF, or equivalently we have the following MCF, for any correlator $\theta$:
    $$\gen{\pt^{n+g-1}}^{|\bfw|,\theta}_{g,a,\bfw} = \sum_{k|a,\bfw,\theta} k^{2n+4g-4} \gen{\pt^{n+g-1}}^{|\bfw/k|,\prim}_{g,a/k,\bfw/k},$$
    where $k|\theta$ means that $\theta=k\theta'$ for some $\theta'\in E[|\bfw|]$, which means that $\theta\in E[|\bfw/k|]$.
\end{theo}

In \cite[Section~6]{blomme2024correlated}, we proved that the correlated invariants can be explicitly computed by counting so-called \textit{floor diagrams} with a correlated multiplicity, which we recall below. Then, we proceed as follows:
    \begin{itemize}[leftmargin=0.6cm]
        \item Up to taking the $(\theta)$-coefficients and summing over the correlators, the two proposed versions of the MCF in Theorem \ref{theo-MCF-EP1-points} are indeed equivalent.
        \item We then show that $\gen{\gen{\pt^{n+g-1}}}_g\colon\RRR\to\G$ is actually the push-forward of a multiplicity function $\bfm\colon\FFF\to\G$ on the $\NN$-module $\FFF$ of diagrams, along the natural \textit{degree} morphism $\deg\colon\FFF\to\RRR$. By Lemma \ref{prop:push-forward} it is sufficient to prove the MCF for $\bfm$.
        \item Last, we prove MCF for $\bfm$. Restricting to the orbit of a single \textit{primitive} floor diagram $\widetilde{\Ffk}$, $\bfm(\delta\widetilde{\Ffk})$ is actually the product of a monomial by a product indexed by vertices. By Lemma \ref{lem:monomial-rule}\ and ref{lem:product-rule}, it suffices to show the MCF for the vertex multiplicities, which follow from their explicit expression in Proposition \ref{prop:vertex-contrib-EP1}.
    \end{itemize}

\subsubsection{Floor diagrams}

\begin{defi}\label{defi:floor-diagrams}
    A floor diagram $\Ffk$ is the data of a weighted oriented graph with the following properties:
    \begin{enumerate}[leftmargin=0.6cm]
        \item The graph has three kind of vertices:
            \begin{itemize}
                \item sinks and sources which are univalent vertices, referred as \textit{infinite vertices},
                \item \textit{flat vertices} (set $V_m$) which are bivalent with one ingoing and one outgoing edge,
                \item \textit{floors} (set $V_f$) which carry a label $a_v\in\NN$.
            \end{itemize}
        \item The set of flat vertices and floors carries a total order compatible with the orientation.
        \item The edges have a positive weight $w_e$, such that the weighting makes floors and flat vertices balanced.
        \item The complement of all flat vertices is without cycle and each connected component of the complement contains a unique infinite vertex.
    \end{enumerate}
    Edges adjacent to an infinite vertex are called \textit{ends} (set $E_\infty$), other edges are called \textit{bounded} (set $E_b$). The genus of a floor diagram is $b_1(\Gamma)+|V_f|$. Its degree is $(a,\bfw)$ where $a=\sum a_v$ and $\bfw$ is the signed collection of weights of ends.
\end{defi}
Degeneration graphs that contribute with non-zero multiplicity in the degeneration formula correspond to floor diagrams. 
We refer to \cite[Section~6]{blomme2024correlated} for details on how to see this.

\medskip

Let $\FFF$ be the set of floor diagrams. We define a $\NN$-action on $\FFF$ by scaling edge weights $w_e$ and vertex labels $a_v$. This makes $\FFF$ into a $\NN$-module where the norm is induced by the morphism of $\NN$-modules
$$\deg\colon\Ffk\longmapsto \left(\sum a_v,\bfw(\Ffk)\right)\in\RRR.$$
In other words, $|\Ffk|=|\bfw(\Ffk)|$ is the g.c.d. of the weights of the ends. We also denote by $\ell(\Ffk):=\gcd(w_e,a_v)_{e\in E(\Ffk),v\in V(\Ffk)}$ the  gcd of all edges weights and vertex labels and refer to it as the gcd of $\Ffk$.

\medskip

An application of the correlated decomposition formula \cite[Corollary 4.12]{blomme2024correlated} tells us how to define the correlated  multiplicity $\bfm(\Ffk)$ so that the sum over floor diagrams of degree $(a,\bfw)$ yields $\gen{\gen{\pt^{n+g-1}}}^{|\bfw|}_{g,a,\bfw}$.

\begin{prop}\cite[Definition 6.4, Theorem 6.5]{blomme2024correlated}
\label{prop-mult-EP1-points}
    We have
    $$\gen{\gen{\pt^{n+g-1}}}^{|\bfw|}_{g,a,\bfw} = \sum_{\Ffk}\bfm(\Ffk),$$ where the sum is over genus $g$ degree $(a,\bfw)$ floor diagrams, counted with the multiplicity
    $$\bfm(\Ffk) = \sideset{}{'}\prod w_e^2\prod w_e\cdot\prod_v a_v^{n_v-1}\bsigma^{|\Ffk|}(a_v)\cdot T_{|\Ffk|/\ell(\Ffk)},$$
    where:  the product  $``\sideset{}{'}\prod"$ is  over edges (bounded or not) not adjacent to a flat vertex; the product $\prod$ is over the set of bounded edges; the product $\prod_v$ is over floors.
\end{prop}

Equivalently, $\gen{\gen{\pt^{n+g-1}}}_g=\deg_*\bfm$, where $\deg\colon\FFF\to\RRR$. Furthermore, notice that $|\Ffk|/\ell(\Ffk)$ is actually the norm of the primitive diagram whose orbit contains $\Ffk$.

\begin{proof}[Sketch of proof, details in \cite{blomme2024correlated}]
Let $\Ffk$ be a decorated diagram. Let $\delta_v$ be the gcd of edge weights adjacent to the vertex $v$, and $\delta_{\Ffk}$ the gcd of all edge weights.
By definition of floor diagrams, flat vertices on $\Ffk$ are disposed in such a way  that their complement is a union of trees, each one containing a unique end. This gives  an orientation of the edges of $\Ffk$. Therefore, there is a unique assignment of classes coming from the K\"unneth decomposition of the diagonal yielding non-zero contribution. The multiplicity comes from two terms:
\begin{itemize}[label=$\circ$]
 \item The vertex contributions $a_v^{n_v-1}\cdot w_{v_\mathrm{out}}^2\cdot\bsigma^{\delta_v}(a_v)$, computed in Proposition \ref{prop:vertex-contrib-EP1}, where $w_{v_\mathrm{out}}$ is the weight of the unique edge going out of $v$ (and thus not adjacent to a flat vertex).
 \item The product over bounded edges $\prod w_e$ corresponding to the choices of log-structures needed to lift the map obtained by gluing the log maps at the vertices to a log map to the central fiber.
\end{itemize}
A floor diagram $\Ffk$ contributes in the degeneration formula to a correlator $\theta$ if and only if
$$\frac{|\bfw|}{\delta_{\Ffk}}\theta = \sum_v\frac{\delta_v}{\delta_{\Ffk}}\theta_v,$$
and the contribution spreads  uniformly among the possibilities. So the multiplicity is
$$\bfm(\Ffk) = \sideset{}{'}\prod w_e^2\cdot\prod w_e\cdot \d{\frac{1}{|\bfw|/\delta_{\Ffk}}}\left(\prod_v a_v^{n_v-1}\m{\delta_v/\delta_{\Ffk}}\bsigma^{\delta_v}(a_v)\right),$$
where the first product comes from the edges coming out of non-flat vertices, which are precisely edges not adjacent to a flat vertex.
Using properties of the refined divisor functions, we have that
$$\d{\frac{1}{|\bfw|/\delta_{\Ffk}}}\m{\delta_v/\delta_{\Ffk}}\bsigma^{\delta_v}(a_v) = \d{\frac{1}{|\bfw|/\delta_{\Ffk}}}\bsigma^{\delta_\Ffk}(a_v) = \bsigma^{|\bfw|}(a_v)T_{|\bfw|/\delta_{\Ffk}} = \bsigma^{|\bfw|}(a_v)T_{|\bfw|/\gcd(\delta_\Ffk,a_v)}.$$
An elementary computation yields the expression from the theorem.
\end{proof}
\begin{rem}
    The expression given here for $\bfm(\Ffk) $ looks slightly different than the one in \cite[Definition~6.4]{blomme2024correlated}. It can be obtained from the latter simply applying the arithmetic properties of Lemma~\ref{lem:arithmeticproperties} precisely as done at the end of the proof.
\end{rem}

\subsubsection{Multiple cover formula in $E\times\PP^1$ for point insertions}

\begin{lem}\label{lem:homogeneity-cpx-case}
    The multiplicity $\bfm\colon\FFF\to\G$ satisfies the $(2n+4g-4)$-MCF.
\end{lem}

\begin{proof}
   By Lemma~\ref{lem:reductiontoN}, to prove the MCF  we can restrict to the orbit of a given primitive floor diagram $\tFfk$. We theb have
    $$\bfm(\delta\tFfk) = \sideset{}{'}\prod (\delta w_e)^2\prod \delta w_e\cdot\prod_v (\delta a_v)^{n_v-1}\bsigma^{\delta|\tFfk|}(\delta a_v)\cdot T_{|\tFfk|}.$$
    Since the orbit $\NN\cdot\widetilde{\Ffk}$ is isomorphic to $(\NN,|\cdot|_{|\widetilde{\Ffk}|})$, we can apply Lemma \ref{lem:MCF-for-N-k}. The multiplicity is invariant by $|\tFfk|$-torsion. We then apply $\m{|\tFfk|}$ and reduce to MCF for a diagonal sequence in $\G^{\NN}$, with the standard norm on $\NN$. By Lemma~\ref{lem:arithmeticproperties},
    $$\m{|\tFfk|}\bfm(\delta\tFfk) = \sideset{}{'}\prod (\delta w_e)^2\prod \delta w_e\cdot\prod_v (\delta a_v)^{n_v-1}\bsigma^{\delta}(\delta a_v),$$
    which is the product of a monomial in $\delta$ by a product of functions satisfying the $0$-MCF by Proposition \ref{prop-divisor-fct}. By Lemma~\ref{lem:monomial-rule},
    we only need to compute the exponent of the monomial, which we do in the next Lemma.   
\end{proof}

\begin{lem}
    The exponent of $\delta$ in the MCF above is
    $$2|E'|+|E_b|+\sum_{V_f}(n_v-1) = 2n+4g-4.$$
\end{lem}

\begin{proof}
    Let $\Ffk$ be a floor diagram. Let $E'\subset E_b\sqcup E_\infty$ be the set of edges not adjacent to a flat vertex. We have the following relations:
    \begin{itemize}
        \item $2|V_m| = |E_b|+|E_\infty|-|E'|$, obtained counting the flags adjacent to flat vertices,
        \item $2|V_m|+\sum_{V_f} n_v = 2|E_b|+|E_\infty|$, obtained counting all flags of the graph,
        \item $|V_f|+|V_m|-|E_b| = 1-(g-|V_f|)$, by the genus condition,
        \item $|V_f|+2|V_m|-|E_b| = n$, since for each diagram with non-zero multiplicity the complement of all flat vertices does not  have cycles and each connected components contains a unique end (see for Example \cite[Lemma~6.6,6.7]{blomme2024correlated}).
    \end{itemize}
    The number of ends is equal to $n$. Therefore taking linear combination of the identities above we have:
    $$\left\{ \begin{array}{rl}
        |V_m| = & |E_b|+1-g, \\
        |V_f| = & n+2g-2-|E_b|\\
        |E'| = & n+2g-2-|E_b|,\\
        \sum_{V_f}n_v = & n+2g-2.
    \end{array}\right.$$
   from which we get the expected exponent.
\end{proof}

\begin{proof}[Proof of Theorem \ref{theo-MCF-EP1-points}]
    We know that $\bfm$ satisfies the $(2n+4g-4)$-MCF. Therefore, by Proposition \ref{prop:push-forward}, we also have the MCF with same exponent for
    $$\gen{\gen{\pt^{n+g-1}}}_g\colon (a,\bfw)\mapsto\gen{\gen{\pt^{n+g-1}}}^{|\bfw|}_{g,a,\bfw}\in\G_{|\bfw|}.$$
    Since $\Prim_{|\bfw|}\gen{\gen{\pt^{n+g-1}}}^{|\bfw|}_{g,a,\bfw} = |\bfw|^2\gen{\pt^{n+g-1}}^{{|\bfw|},\prim}_{g,a,\bfw}$, we get
    \begin{align*}
        \gen{\gen{\pt^{n+g-1}}}^{|\bfw|}_{g,a,\bfw} = &\sum_{k|a,\bfw} k^{2n+4g-4} \Prim_{|\bfw|/k}\gen{\gen{\pt^{n+g-1}}}^{|\bfw|/k}_{g,\frac{a}{k},|\bfw|/k} T_{|\bfw|/k}\\
         = & \sum_{k|a,\bfw} k^{2n+4g-4}\gen{\pt^{n+g-1}}^{|\bfw/k|,\prim}_{g,a/k,\bfw/k} |\bfw/k|^2 T_{|\bfw|/k}\\
         =&\sum_{k|a,\bfw} k^{2n+4g-4}\gen{\pt^{n+g-1}}^{|\bfw/k|,\prim}_{g,a/k,\bfw/k} \sum_{|\bfw/k|\theta=0}(\theta)\\
         \sum_{|\bfw|\theta\equiv 0} \gen{\pt^{n+g-1}}^{|\bfw|,\theta}_{g,a,\bfw}\cdot (\theta) = & \sum_{|\bfw|\theta\equiv 0} \sum_{k|a,\bfw,\theta}k^{2n+4g-4}\gen{\pt^{n+g-1}}^{|\bfw/k|,\prim}_{g,a/k,\bfw/k}(\theta). \\
          \end{align*}
    We deduce the equivalent formulation for $\gen{\pt^{n+g-1}}^{|\bfw|,\theta}_{g,a,\bfw}$ looking at the $(\theta)$-coefficient. Going in the reverse order of the computation, we get that the formulations are in fact equivalent.
\end{proof}

\subsection{Multiple cover formula for abelian surfaces}
\label{sec:MCF-abelian-points}

Let $A$ be an abelian surface and  $\beta\in H_2(A,\ZZ)\cong U^{\oplus 3}$ (here $U$ denotes the hyperbolic lattice)  an effective realizable class (i.e. $\beta$ is the class of an algebraic curve $C$ in A), we consider the \emph{reduced GW-invariant}
$$N_{g,\beta} = \gen{\pt^g}^A_{g,\beta} = \int_{\red{\M_{g,g}(A,\beta)}} \prod_1^g\ev_i^*(\pt).$$

Deformation invariance implies that $N_{g,\beta}$ only depends on $\beta$ through its self-intersection $\beta^2$ and its divisibility $\ell(\beta)$. A class is \emph{primitive} if $\ell(\beta)=1.$

\begin{theo}\label{theo-MCF-abelian-points}
    The reduced GW-invariants satisfy the multiple cover formula
    $$N_{g,\beta} = \sum_{k|\beta} k^{4g-3}N_{g,\widetilde{\beta/k}},$$
    where $\widetilde{\beta/k}$ is a primitive class having the same self-intersection as $\beta/k$.
\end{theo}

Though longer to write since the computations have not yet being carried out in \cite{blomme2024correlated} due to the lack of reduced decomposition formula, the proof goes along the exact same steps as the $E\times\PP^1$ case:
    \begin{itemize}[leftmargin=0.6cm]
        \item Considering several elliptically fibered abelian surfaces, rephrase Theorem \ref{theo-MCF-abelian-points} as an MCF for a function $\bfN$ with values in $\G$.
        \item Express $\bfN$ as the push-forward of a multiplicity function $\bfm$ on a $\NN$-module of diagrams $\DDD$ via (reduced) decomposition formula.
        \item Prove the MCF for $\bfm$ by expressing it as a product over vertices and edges and conclude by Lemma \ref{lem:product-rule} and Proposition \ref{prop:push-forward}.
    \end{itemize}

We denote by $\BBB=\NN^2$ the $\NN$-module of \textit{diagram degrees}, endowed with its natural $\NN$-action and norm given by the projection onto the first coordinate. The $\NN$-module $\BBB$ will record, for elliptically fibered abelian surfaces, the intersection of the curve class with a fiber and with a section respectively.

\subsubsection{Decomposition formula and diagrams}
\label{sec:pearl-diagrams}

Let $A_t(u)\to T$ be a one parameter degeneration of abelian surfaces  constructed as in Section~\ref{sec-construction-abelian-families}, with monodromy data $u$, and whose central fiber has $g$ irreducible components Let $\beta(u)$ be the effective realizable class given in Lemma \ref{lem:polarization-u}. 
We want to apply the reduced degeneration formula (Theorem~\ref{theo:degene-formula-reduced}, Corollary~\ref{coro:degen-gp-alg}) distributing the $g$ point constraints  one on each irreducible component of the central fiber $A_0(u)$. 
The reduced degeneration formula, phrased as in Section~\ref{sec:correlatorsindegeneration}, asserts that $N_{g,\beta(u)}$ is expressed as the $(u)$-coefficient of a sum over \textit{pearl diagrams} (alias the abelian surfaces version of \emph{floors diagrams}, already defined in \cite{blomme2022abelian3} and recalled below) counted with a \emph{correlated multiplicity}. Pearl diagrams are actually the rigid tropical curves contributing with a nn-zero multiplicity.

\begin{defi}\label{defi:pearl-diagrams}
    A (pearl) diagram $\Dfk$ is the data of a weighted oriented graph with the following properties:
    \begin{enumerate}[leftmargin=0.6cm]
        \item The graph has two kind of vertices:
            \begin{itemize}
                \item \textit{flat vertices} (set $V_m$) which are bivalent with one ingoing and one outgoing edge,
                \item \textit{pearls} (set $V_f$) which carry a label $a_v\in\NN$.
            \end{itemize}
        \item There is a map $\Dfk\to \Sigma_0 $ for $\Sigma_0$ the oriented cycle graph with $|V_f|+|V_m|$ vertices; the map is 1-to-1 on vertices and compatible with the orientation.
        \item The edges are decorated with weight $w_e$,  such that the signed (with respect to the orientation on $\Sigma_0$)  sum of weights at floors and flat vertices  is zero, namely $\Dfk\to\Sigma_0 $ is balanced and thus a \textit{tropical cover}.
        \item The complement of all flat vertices is a tree (i.e. connected and without cycles).
    \end{enumerate}
    The genus of a pearl diagram is $b_1(\Gamma)+|V_f|$. Its degree is  a \emph{diagram degree}  $B=(|B|,a)\in\BBB$ where $a=\sum a_v$ and $|B|$ is the degree of the tropical cover.
\end{defi}

We denote by $\DDD$ the set of pearl diagrams. It is endowed with a $\NN$-action by scaling edge weights and vertex labels. It thus becomes a $\NN$-module whose norm is induced by the degree map $\deg\colon\DDD\to\BBB$, namely 
$|\Dfk|=|B|$. When $\Dfk$ is fixed, we denote simply by $\delta$ (rather than $\delta_{\Dfk}$) the $\gcd$ of all edge weights $(w_e)$, and by $\ell(\Dfk)=\gcd((a_v),\delta).$

\medskip

If $h_e^\pm$ are the two half-edges corresponding to an edge $e$ (orientation going from $h_e^-$ to $h_e^+$), we set $w_{h_e^\pm}=\mp w_e$, so that the weight for a half-edge going out of a vertex is positive, and negative for a half-edge going into a vertex. In particular, balancing condition asserts that for each vertex $v$, $\sum_{h\vdash v}w_h=0$.

\medskip

As for the case of floor diagrams recalled before, pearl diagrams correspond to rigid tropical curves (Definition~\ref{defi:rigidtropmap}, also called degeneration graphs) with non-zero multiplicity in the degeneration formula for $A_0(u).$ In particular:
\begin{itemize}
    \item flat vertices correspond to genus zero components carrying a marked point;
    \item  pearls correspond to genus $1$ components, also carrying a marked point,  mapping to a component $E\times\PP^1\subset \widetilde{A}_0$ with class $a_v[E]+(\sum_{h\vdash v,w_h>0}w_h)[\PP^1];$ here $\widetilde{A}_0$ is the normalization of of $A_0(u)$, i.e. a disjoint union of $E\times\PP^1$.
\end{itemize}

One can show, with a proof analogous to the floor diagram case (see for example \cite[Proposition~6.6,6.7,6.8]{blomme2024correlated}, \cite[Section~4.1]{blomme2022abelian3},\cite[Section~4.2]{blomme2024bielliptic}), that if a degeneration graph gives a non-zero multiplicity then it is a pearl diagram.
We give a sketch of the argument in the following Lemma.

\begin{lem}
\label{lem:position-points-abelian-case}
Let $\phi\colon\Gamma\to\Sigma_0$ be a rigid tropical curve with non-zero multiplicity in the reduced degeneration formula, then:
\begin{enumerate}[leftmargin=0.6cm]
    \item the vertices have genus $0$ or $1$, and if the genus in $1$ then they must carry a marked point (these are the pearls), and the genus $0$ vertices are bi-valent in $\Gamma;$
   \item the complement of the genus $0$ marked vertices, i.e. the flat vertices, is a tree;
    \item choosing one edge and using the expression of the reduced diagonal determined by this edge, the reduced Poincar\'e insertion are uniquely determined:
        \begin{itemize}
            \item for a flat vertex (marked genus 0) they are $1,1\in H^*(F,\mathbb Z)$;
            \item for a non marked  genus 0 vertex $1,\pt\in H^*(F,\mathbb Z)$;
            \item for a pearl $1, \pt,\dots,\pt\in H^*(F,\mathbb Z)$, with $1$ for the unique edge going out of the pearl (once we cut at flat vertices).
        \end{itemize} 
\end{enumerate}
\end{lem}

\begin{proof}
\begin{enumerate}[leftmargin=0.6cm]
    \item The proof of the first statement is identical to \cite[Proposition~6.6]{blomme2024correlated}, resulting from a dimension count.
    
   \item[(2,3)] Let $\hat{\Gamma}$ be the graph obtained from $\Gamma$ cutting at the flat vertices. We prove (2) and (3) at the same time, showing $\hat{\Gamma}$ is a tree and determining the insertions of the diagonal.
   
   Let $V_f$ be the set of pearls, $V_0$ the unmarked genus zero (bivalent) vertices, $2V_m$ the set of univalent vertices arising from the cut and $E_b$ the set of bounded edges.

   An Euler characteristic computations shows that
   $$\chi(\hat{\Gamma}) = 1-b_1(\Gamma)+|V_m|.$$
   Since the number of marked points matches the genus of the pearl diagram (pearls being of genus 1), we have furthermore $b_1(\Gamma)+|V_f|=g=|V_f|+|V_m|$. Therefore,
   $$\chi(\hat{\Gamma})=1.$$
   In particular, since it is positive, at least one connected component of $\hat{\Gamma}$ needs to be a tree.

   Let $e_0$ be any edge of $\Gamma$. We use the expression of the reduced diagonal $\Delta_u^{\text{red}}$ coming from the K\"unneth type decomposition discussed in Section~\ref{sec:numericaldec}. Dimensional constraints force the following facts on the insertions coming from the diagonal:
    \begin{itemize}
        \item for half edges adjacent to flat vertices, insertions are $1$;
        \item if $e=(h,h')$ is an edge different from $e_0$ and the insertion at $h$ is a $1$, the insertion at $h'$ is a $\pt$;
        \item if one of the two insertions at a vertex of $V_0$ is a $\pt$, the other is a $1$;
        \item if insertions for all half edge adjacent to a pearl except one are $\pt$, then the remaining one is a $1$.
    \end{itemize}
    Using the above facts, we can inductively prune the branches of $\hat{\Gamma}$ determining the insertions. This process terminates if $\hat{\Gamma}$ has no cycles.

    If $\hat{\Gamma}$ is disconnected, we can take the edge $e_0$ (not appearing in the reduced diagonal)  in the complement of tree connected component. We then have a contradiction when finding the insertions at the half edges of this tree component: at some point during the pruning, an edge will get two $1$ insertions at its half edges or a vertex will only get point insertions at its adjacent flags. Hence, $\hat{\Gamma}$ is connected, and thus a tree since it has Euler characteristic $1$. In particular, all insertions have been determined by the above rules.
\end{enumerate}

\end{proof}

\begin{rem}
    The bivalent vertices of genus $0$ without marked point do not contribute to the multiplicity and can thus be deleted from the graph and merging the adjacent edges into a single edge.
\end{rem}

    The following lemma completes the description of the pearl diagrams by expressing  all the combinatorial invariants of the diagram in terms of the number of bounded edges and $g$.

\begin{lem}\label{lem:number-edges-vertex-pearl-diag}
Let $\Dfk$ be a genus $g$ pearl diagram, $E_b$ the set of bounded edges, $E'$ the set of bounded edges not adjacent to a flat vertex, $V_m$ the set of flat vertices, and $V_f$ the set of pearls. We have the following identities, and it has in particular $g$ vertices:
    $$\left\{\begin{array}{rl}
        |V_m| = & 1-g+|E_b| \\
        |V_f| = & 2g-1-|E_b| \\
        |E'| = & 2g-2-|E_b| \\
        \sum_{V_f} n_v = & 2g-2.
    \end{array}\right.$$
\end{lem}

\begin{proof}
    The statement stems from the following relations:
    \begin{itemize}
        \item $2|V_m|=|E_b|-|E'|$, obtained counting edges adjacent to bivalent marked vertices,
        \item $2|V_m|+\sum_{V_f}n_v = 2|E_b|$, obtained counting flags of the graphs in two ways,
        \item $|V_f|+|V_m|-|E_b|=1-(g-|V_f|)$ by the genus condition,
        \item $|V_f|+2|V_m|-|E_b|=1$ since the complement of bivalent marked vertices is a tree.
    \end{itemize}
\end{proof}

We now compute the correlated multiplicity associated to pearl diagrams applying  the reduced decomposition formula.

\begin{prop}
    We define the correlated multiplicity of a pearl diagram $\Dfk$ by
    $$\bfm(\Dfk) = |\Dfk|^2\sideset{}{'}\prod_e w_e^2\prod_e w_e\cdot\prod_v a_v^{n_v-1}\bsigma^{|\Dfk|}(a_v)\cdot T_{|\Dfk|/\ell(\Dfk)}$$
    where we recall that $\sideset{}{'}\prod$ is the product over the edges not adjacent to flat vertices.
    The reduced Gromov-Witten invariant $N_{g,\beta(u)}$ is the sum over genus $g$ degree $B$ pearl diagrams $(u)$-coefficient of $\bfm(\Dfk)$:
    $$N_{g,\beta(u)} = \sum_\Dfk \operatorname{Coeff}_{(u)}\bfm(\Dfk).$$
\end{prop}

\begin{proof}
The reduced decomposition formula in the formulation of Corollary~\ref{coro:decompo-num-vers} gives the following recipe: choose an edge $e_0$ of $\Gamma$ and use the K\"unneth expression of the diagonal forgetting about the edge $e_0$. As computed in Section~\ref{sec:numericaldec} this differs from the expression of the reduced diagonal by a factor $\frac{w_{e_0}^2}{\delta^2}$. By Lemma \ref{lem:position-points-abelian-case}, all the insertions are uniquely determined, and the multiplicity involves the following quantities:
    \begin{itemize}
        \item the product $\prod w_e$ over bounded edges coming from choices of the log structure needed to lift the map obtained from the gluing to a log map to $A_0$
        \item the vertex contributions:
            $$\prod_v a_v^{n_v-1}w_{v_\mathrm{out}}^2 \bsigma^{\delta_v}(a_v),$$
            where $w_{v_\mathrm{out}}$ is the weight of the unique edge  with trivial constraint, determined via the pruning procedure once we choose the edge $e_0.$
        \item $\frac{\delta^2}{w_{e_0}^2}$, the factor appearing in the reduced diagonal class and its K\"unneth decomposition obtained  forgetting the edge $e_0$.
    \end{itemize}
The compatibility condition \eqref{eq:contributingtheta} on correlators contributing non trivially to the reduced decomposition formula for $A_t(u)\to T$  reads:
$$\sum_v\frac{\delta_v}{\delta}\theta_v = \frac{|\Dfk|}{\delta}u,$$
where $\delta=\gcd(w_e)$. Exploiting the group algebra product we can make the left-hand-side naturally appear when writing the invariant as a sum of correlated multiplicities. To do so, we apply the operator $\m{\delta_v/\delta}$ to the vertex multiplicity $a_v^{n_v-1}w_{v_\mathrm{out}}\bsigma^{\delta_v}(a_v)$. This has the effect of replacing each $(\theta_v)$ in the sum by $(\frac{\delta_v}{\delta}\theta_v)$. We then select the $(\frac{|\Dfk|}{\delta}u)$-coefficient in the product.

It is convenient to apply the operator $(|\Dfk|/\delta)^2\d{\frac{1}{|\Dfk|/\delta}}$, which maps a generator $(t)$ to the sum of its $\frac{|\Dfk|}{\delta}$-roots, and then instead  select the $(u)$-coefficient. This way we recognize in the $(u)$-coefficient the contribution for the abelian surface with monodromy $u$. In the end, the multiplicity is therefore
    $$\bfm(\Dfk)=\prod_e w_e \cdot \frac{\delta^2}{w_{e_0}^2}\cdot (|\Dfk|/\delta)^2\d{\frac{1}{|\Dfk|/\delta}}\left(\prod_v a_v^{n_v-1}w_{v_\mathrm{out}}^2\m{\delta_v/\delta}\bsigma^{\delta_v}(a_v)\right).$$
    Notice that if the $e_0$ we forget is adjacent to a flat vertex, the $w_{e_0}^2$ cancels with the corresponding $w_{v_\mathrm{out}}^2$ for the adjacent pearl. If $e_0$ is adjacent to two pearls, the $w_{e_0}^2$ cancels with one of the two $w_{v_\mathrm{out}}^2$ factors coming from the adjacent pearls. In both cases, we get
    $$\bfm(\Dfk) = |\Dfk|^2\sideset{}{'}\prod w_e^2\prod w_e\d{\frac{1}{|\Dfk|/\delta}}\left(\prod_v a_v^{n_v-1}\m{\delta_v/\delta}\bsigma^{\delta_v}(a_v)\right),$$
    which does not depend on the choice of $e_0$.
    We simplify the expression as in Proposition \ref{prop-mult-EP1-points} to get the advertised form for $\bfm(\Dfk)$:
    $$\d{\frac{1}{|\Dfk|/\delta}}\m{\delta_v/\delta}\bsigma^{\delta_v}(a_v)=\d{\frac{1}{|\Dfk|/\delta}} \bsigma^{\delta}(a_v)=\bsigma^{|\Dfk|}(a_v)T_{|\Dfk|/\delta}=\bsigma^{|\Dfk|}(a_v)T_{|\Dfk|/\gcd(\delta,a_v)}.$$
\end{proof}

\begin{rem}
    Notice that for different choices of monodromies $u$ and $u'$, determining different abelian surfaces $A(u), A(u')$ (which in particular carry polarizations with same self-intersection but potentially different divisibility) we could have $\frac{|\Dfk|}{\delta}u=\frac{|\Dfk|}{\delta}u'.$ So, in applying the reduced decomposition formula and considering the contribution of a diagram $\Dfk$, the same correlators $\theta_v$ in the vertex contributions appear.
    Applying the operator $(|\Dfk|/\delta)^2\d{\frac{1}{|\Dfk|/\delta}}$ has the effect of duplicating the coefficient of $\frac{|\Dfk|}{\delta}u$ in such a way that it appears as coefficient for both $u$ and $u'.$
   \end{rem}

\subsubsection{Multiple cover formula in abelian surfaces for point insertions}

\begin{lem}\label{lem-MCF-abelian-points}
    The multiplicity $\bfm\colon\DDD\to\G$ satisfies the $(4g-3)$-MCF.
\end{lem}

\begin{proof}
    We proceed as in the $E\times\PP^1$ case. We can restrict to the orbit a chosen primitive diagram $\tDfk$. We have that
    $$\bfm(\delta\tDfk) = |\delta\tDfk|^2\sideset{}{'}\prod_e (\delta w_e)^2\prod_e \delta w_e \cdot\prod_v (\delta a_v)^{n_v-1}\bsigma^{\delta|\tDfk|}(\delta a_v)\cdot T_{|\tDfk|}.$$
    The orbit $\NN\cdot\tDfk$ is isomorphic to $(\NN,|\cdot|_{|\tDfk|})$, so that we may apply Lemma \ref{lem:MCF-for-N-k}. The multiplicity is invariant by multiplication by $T_{|\tDfk|}$ so that we may apply $\m{|\tDfk|}$ and prove instead the usual MCF in $\NN$. Applying $\m{|\tDfk|}$ has the effect of replacing $\bsigma^{|\tDfk|\delta}$ by $\bsigma^\delta$ and removing $T_{|\tDfk|}$. The multiplicity thus appears as the product of a monomial with a product of $\bsigma^\delta(\delta a_v)$, precisely as in the case of $E\times\PP^1$. Therefore, we conclude it satisfies the MCF, for which we merely need to compute the exponent. Using Lemma \ref{lem:number-edges-vertex-pearl-diag},
    \begin{align*}
        & 2+2|E'|+|E_b|+\sum (n_v-1) \\
    = & 2+2(2g-2-|E_b|)+|E_b|+(2g-2)-(2g-1-|E_b|) = 4g-3.
    \end{align*}
\end{proof}

\begin{proof}[Proof of Theorem \ref{theo-MCF-abelian-points}]
We proceed as in the $E\times \mathbb P^1$ case, with the only difference coming from the geometric interpretation of the various steps.
Here the 
 function $\BBB\to\G$ with values in the group algebra $\G=\QQ[E]$
 we use to encode the invariants, mixes the Gromov-Witten $N_{g,\beta(u)}$ invariants of different elliptically fibered abelian surfaces $A(u)$, for which the classes $\beta(u)$ have same square but different divisibility. This should already sound promising, since the statement of the multiple cover formula for abelian surfaces conjectured in \cite{oberdieck2016curve} does indeed involves different abelian surfaces.
 
\medskip

Let $B$ be a diagram degree. For any $u$ element of $|B|$-torsion in the elliptic curve $E$, we have defined a family of abelian surfaces with $A_t(u)\to T$ with $A_0(u)$ having monodromy data $u$, and endowed with a polarization $\beta(u)$. We consider the group algebra valued function $\bfN_g\colon\BBB\to\G$ defined by
$$\bfN_g(B) = \sum_{|B|u\equiv 0} N_{g,\beta(u)}\cdot(u).$$
Applying the reduced degeneration formula, $\bfN_g(B)$ can be expressed as a sum over diagrams of genus $g$ and degree $B$. In other words, $\bfN_g=\deg_*\bfm$. Since $\bfm$ satisfies that $(4g-3)$-MCF, so does $\bfN_g$ and we have
\[\bfN_g(B)=\sum_{k| B} k^{4g-3}\Prim_{|B|/k}(\bfN_g(B/k)) T_{|B/k|}.\]

\medskip

Last, we prove that the group algebra formulation is actually equivalent to the one stated in the Theorem. We saw that the divisibility of the class $\beta(u)$ considered in Lemma~\ref{lem:polarization-u}, is $\gcd(B,k,l)$  where $B=(|B|,a)$  is the diagram degree and $u=\frac{k+l\tau'}{|B|}$. This can be rewritten as $\gcd(B,\frac{|B|}{\mathrm{ord}(u)})$. Therefore, $k|\beta(u)$ is equivalent to $k|B,u$, where $k|u$ means that $u=ku'$ for $u'$ some other element of $|B|$-torsion. On the other hand
\[\Prim_{|B|/k}(\bfN_g(\frac{B}{k}))=|B/k|^2N^{\prim}_{g,B/k},\] where $N^{\prim}_{g,B/k}$ is the reduced GW invariant for the class constructed in Lemma~\ref{lem:polarization-u} with diagram degree $B/k$ and a primitive choice of monodromy $u$, i.e. any $u$ of order $|B|/k$. Notice that this class has divisibility $1$ and its self-intersection is $(B/k)^2.$
The $|B/k|^2$ factor simply comes from the number of $|B/k|$-torsion elements. So we have
\begin{align*}
\bfN_g(B)&=\sum_{k| B} k^{4g-3}\Prim_{|B|/k}(\bfN_g(B/k)) T_{|B/k|}\\
&= \sum_{k| B} k^{4g-3} N^{\prim}_{g,B/k} |B/k|^2T_{|B/k|}\\
&=\sum_{k| B} k^{4g-3} N^{\prim}_{g,B/k} \sum_{|B/k|u=0}(u)\\
&=\sum_{|B|u\equiv 0} \left(\sum_{k|B,u} k^{4g-3}N^\prim_{g,B/k}\right) \cdot (u)
\end{align*}
which, taking the $(u)$-coefficient, gives us the MCF for reduced GW-invariants for abelian surfaces. Reading the equalities the other way around tells us that the two formulations are equivalent.

\end{proof}

\subsection{Extension to refined invariants}
\label{sec-MCF-lambda}

With little work (essentially the computation of the vertex multiplicities), the paradigma of the previous section can be adapted to also prove that reduced GW invariants with points and a $\lambda$-class insertion, i.e. a Chern classes of the Hodge bundle, satisfy the MCF. In various works, starting from \cite{bousseau2019tropical}, these GW-invariants have been related  to the so-called tropical \emph{refined invariants}, considered also in \cite{blomme2022floor} and \cite{blomme2022abelian1} in the case of $E\times\PP^1$ and an abelian surface respectively.

The key is that in this situation as well we have an explicit expression for the local contributions (it was computed in \cite{blommecarocci2025DR}), which also satisfy the $0$-MCF. 

\subsubsection{Refined local contribution}
\label{sec:local-contrib-lambda-case}

We go back to $E\times\PP^1$ and consider the following correlated GW invariants:
$$\gen{\gen{\lambda_{g-1};\pt_0,1_{w_1},\pt_{w_2},\dots,\pt_{w_n}}}^{|\bfw|}_{g,a,\bfw} \in\G_{|\bfw|},$$
To state the result in a compact way we introduce the homogeneous polynomial in $\bfw$:
$$W_g(\bfw) = \frac{1}{\prod_1^n w_i}\left(\sum_{S\subset[\![1;n]\!]} (-1)^{|S|}w_S^{n+2g-2}\right)\frac{(-1)^{n+g-1}}{(n+2g-2)!},$$
where $w_S=\sum_{i\in S}w_i$. This is in fact a polynomial since  $\prod_1^n w_i$ divides the expression in the brackets. It is homogeneous of degree $2g-2$.

\begin{prop}
    We have
    $$\gen{\gen{\lambda_{g-1};\pt_0,1_{w_1},\pt_{w_2},\dots,\pt_{w_n}}}^{|\bfw|}_{g,a,\bfw} = w_1^2\cdot W_g(\bfw)\cdot a^{n-1}\bsigma^{|\bfw|}_{2g-1}(a).$$
\end{prop}

\begin{proof}
    In \cite[Theorem 7.7]{blommecarocci2025DR}, we provide a formula for the $(0)$-coefficient of
    $$\gen{\gen{\lambda_{g-1};\pt_0,1_{w_1},\pt_{w_2},\dots,\pt_{w_n}}}^{|\bfw|}_{g,a,\bfw} \in\G_{|\bfw|},$$
    which one can see (using and the usual arithmetic properties of the sum of divisors functions proved in Lemma~\ref{lem:arithmeticproperties}) coincides the $(0)$-coefficient of the right-hand side above. 
    As both  the correlated invariants and the refined sum of divisors function satisfy the unrefinement  relations (see \cite[Lemma~3.12]{blomme2024correlated} and Lemma~\ref{lem:arithmeticproperties} respectively), these two functions from $\RRR\to\G$ with the $\NN$-module structure on $\RR$ already considered before coincide. They already have the same total degree, which is the usual total relative GW-invariant with the given insertions (see the computation in \cite[Lemma~7.6]{blommecarocci2025DR} and references therein).
\end{proof}

\subsubsection{Refined multiple cover formula for correlated invariants of $E\times\PP^1$}

We now care about the reduced invariant
$$\gen{\gen{\lambda_{g-g_0};\pt^{n+g_0-1}}}^{|\bfw|}_{g,a,\bfw}\in\G_{|\bfw|}.$$

It is still possible to apply the correlated decomposition formula from \cite{blomme2024correlated}. We sketch how to adapt the computation.

\medskip

First, we find the floor diagrams, namely those degeneration graphs which contribute with non trivial multiplicity in the degeneration formula. This is performed as in \cite{bousseau2021floor} and has nothing to do with correlation. Using \cite[Proposition~3.1]{bousseau2021floor}, for every degeneration graph $\Gamma$ the $\lambda$-class $\lambda_{g-g_0}$ actually splits as a product over vertices of $\sfV(\Gamma)$, and the graphs with non-zero multiplicity in the degeneration  formula must satisfy $b_1(\Gamma)\leqslant g_0$. Furthermore, computations from \cite{bousseau2021floor} also show that vertices with degree $0$ have to be of genus $0$ and thus do not carry a $\lambda$-class.

\medskip

Let $V_f$ denote the floor vertices i.e. the vertices with $a_v\geqslant 1$ and genus $g_v\geqslant 1$, carrying one interior marking. We again call flats the genus zero marked vertices.
Then $\M_v$ has virtual dimension $n_v+g_v$ where $n_v$ is the number of points with boundary evaluation. In order to have a non-zero contribution we need the degree of the constraints to match the dimension and this imposes further conditions on the degeneration graphs, with a decoration corresponding to the insertions at a given vertex that are relevant.

\medskip

First we notice that over $\M_v=\M_{g_v,m_v}(E\times\PP^1|E_0+ E_{\infty},a_v,\bfw_{e\vdash v})$ the Hodge bundle $\mathbb E_{g_v}$ is equipped with a nowhere vanishing section coming from $f^*\Omega_E\to\Omega_{\C_v}\to\omega_{\pi_v}$ so that the class $\lambda_{g_v}$ vanishes. We therefore refer to $\lambda_{g_v-1}$ as the top $\lambda$-class in this setting. 

\medskip

The degree of the geometric constraints is at most $n_v+1,$ a degree $2$ constraint coming from the interior marked point, 
 and at most a degree $n_v-1$ boundary constraint (the $n_v$ boundary points being either points evaluating at a joining divisor or at the boundary divisor which, they either way satisfy the Menelaus relation, thus imposing at most $n_v-1$ point constraints). Thus the only possibility to get a non-zero invariant is precisely with the top $\lambda$-class insertion.

\medskip

In the end, the floor diagrams, alias the degeneration graphs appearing with non-zero multiplicity, are the same as in the case of $n+g_0-1$ point insertions. One only needs to remember that floors come with a genus $g_v\geqslant 1$ and a $\lambda_{g_v-1}$ insertion coming from the decomposition.


\medskip

In particular, the same combinatorial computation performed in \cite[Lemma~6.7]{blomme2024correlated} replacing the genus $1$ floor vertices with the genus $g_v$ but now capped with a $\lambda_{g_v-1}$ class shows that the complement of the flat vertices (which also in this context are genus $0$ bivalent with an interior marking) has connected components which are trees and contain exactly one infinite vertex/unbounded edge. The pruning algorithm determines the Poincar\'e insertions, precisely as done in the point case.

\medskip

The application of the degeneration formula is then completely analogous with the only difference that in  the computation of the multiplicity, the vertex contributions are replaced by the one of Section \ref{sec:local-contrib-lambda-case}, yielding the following refined correlated multiplicity.

\begin{defi}
    The multiplicity of a genus $g$ floor diagram $\Ffk$ is
    $$\bfm(\Ffk) = \sideset{}{'}\prod_e w_e^2\prod_e w_e \cdot \prod_v a_v^{n_v-1}W_{g_v}(\bfw_v)\bsigma^{|\Ffk|}_{2g_v-1}(a_v)\cdot T_{|\Ffk|/\ell(\Ffk)},$$
    where $\bfw_v$ is the collection of edge weights adjacent to the vertex $v$.
\end{defi}

We then proceed as in the case without $\lambda$-class. Since $\bsigma^\delta_{2g-1}(\delta a)$ also satisfies the $0$-MCF, the only difference lies in the exponent which differs by the degree of $\prod_v W_{g_v}(\bfw_v)$, equal to
$$\sum 2g_v-2 = 2g-2g_0,$$
yielding the following results.

\begin{lem}
    The multiplicity function $\bfm$ satisfies the $(2n+2g_0+2g-4)$-MCF.
\end{lem}

\begin{theo}\label{theo-MCF-EP1-lambda}
    The total correlated invariant $(a,\bfw)\mapsto\gen{\gen{\lambda_{g-g_0};\pt^{n+g_0-1}}}^{|\bfw|}_{g,a,\bfw}$ satisfies the $(2n+2g_0+2g-4)$-MCF. In particular, we have
    $$\gen{\lambda_{g-g_0};\pt^{n+g_0-1}}^{|\bfw|,\theta}_{g,a,\bfw} = \sum_{k|a,\bfw,\theta} k^{2n+2g_0+2g-4}\gen{\lambda_{g-g_0};\pt^{n+g_0-1}}^{|\bfw/k|,\prim}_{g,a/k,\bfw/k}.$$
\end{theo}

\subsubsection{Refined multiple cover formula for abelian surfaces}

We consider the following reduced GW-invariants of an abelian surface $A$:
$$\gen{\lambda_{g-g_0};\pt^{g_0}}^A_{g,\beta}.$$
We can still apply the reduced degeneration formula, and the $\lambda$-class still splits as a product of $\lambda$-classes associated to vertices. Here again, the only difference with pearl diagrams from Section \ref{sec:pearl-diagrams} is that pearls (non-flat vertices) acquire a genus $g_v\geqslant 1$. The vertex multiplicity is yet again replaced by the one from \ref{sec:local-contrib-lambda-case}.

\begin{defi}
    The refined correlated multiplicity of a diagram is given by
    $$\bfm(\Dfk) = |\Dfk|^2 \sideset{}{'}\prod_e w_e^2\prod_e w_e\cdot\prod_v W_{g_v}(\bfw_v)a_v^{n_v-1}\bsigma_{2g_v-1}^{|\Dfk|}(a_v)\cdot T_{|\Dfk|/\ell(\Dfk)}.$$
\end{defi}

The functions $\bsigma^\delta_{2g_v-1}(\delta a_v)$ still satisfy the $0$-MCF, so that the only difference with the point case for abelian surfaces is the exponent that comes out of the monomial part. The difference is the homogeneous part $\prod_v W_{g_v}(\bfw_v)$, which has degree $2g-2g_0$. This yields the following results.

\begin{lem}
The multiplicity satisfies the $(2g+2g_0-3)$-MCF.
\end{lem}

\begin{theo}\label{theo-MCF-abelian-lambda}
    The reduced GW-invariants satisfy the following MCF:
    $$\gen{\lambda_{g-g_0};\pt^{g_0}}_{g,\beta} = \sum_{k|\beta}k^{2g+2g_0-3}\gen{\lambda_{g-g_0};\pt^{g_0}}_{g,\widetilde{\beta/k}},$$
    where $\widetilde{\beta/k}$ is a primitive class with same self-intersection as $\beta/k$.
\end{theo}
\section{Logarithmic intermezzo}\label{sec:logintermezzo}
We have already made a fairly extensively use of logarithmic geometry language in Section~\ref{sec-reduced-decomposition}, and already mentioned (with ad hoc or no explanation) \emph{tropicalization} of a logarithmic scheme and \emph{Artin fans}. Nonetheless these notions where so far only relevant for the proof of the (reduced) decomposition formula and not necessary to understand the statements.

In next three sections, which lead to the proof of the MCF for abelian surfaces in its most complete form, the use of techniques coming from the tropical geometry of logarithmic schemes becomes even more prominent, and we thus recall a few notions. A more extensive introduction, adapted to our needs is given in \cite[Section~1]{MaulikRanganathan2025} and references therein.

\subsection{Tropical models}
\label{sec-tropical-models}


\begin{defi}
    An Artin fan $\mathcal A$ is an algebraic stack with logarithmic structure that admits a strict \'etale cover by Artin cones, namely algebraic stacks of the form $[U_{\sigma}/T_{\sigma}] $ for $U_{\sigma}$ an affine toric variety and $T_{\sigma}$ its dense torus.
\end{defi} 

In particular, Artin fans are logarithmically smooth over the trivial log point and equidimensional of dimension zero.
 We refer for example to \cite{abramovich2014comparison,abramovich2016skeletons, cavalieri2020moduli} for an extensive introduction to Artin fans and their first applications in log enumerative geometry. 

 \medskip
 
 Artin fans are combinatorial objects in the sense that by \cite{cavalieri2020moduli}, to each Artin fan $\mathcal A$ corresponds bijectively a cone stack $\Sigma(\mathcal A)$, namely a diagram of rational polyhedral cones $(\sigma, N)$ with face maps $\tau\hookrightarrow\sigma$. The difference between a standard cone complex and a cone stack is in the type of face maps allowed in the diagram: in a cone stack we allow more than one face maps between two cones and self maps different from the identity.

\medskip
 
As explained in \cite[Section~1.3]{MaulikRanganathan2025}, for any algebraic stack $\mathcal X$ with fs (fine and saturated) log-structure, there exists a morphism to some Artin fan $\mathcal X\to\mathcal A$ such that the log-structure on $X$ is pulled back from the one on $\mathcal A$. The previous map, seen as a map of log-algebraic stacks is tautologically strict.

\begin{rem}
Some reader might be more familiar with the notion of \emph{the} Artin fan $\mathcal A_{\mathcal X}$ of an algebraic stack with log-structure as introduced by Abramovich-Wise:
this is the \emph{initial factorization} of the natural morphism $\mathcal X\to\mathcal Log$, where $\mathcal Log$  is Olsson's stack parametrizing log structures \cite{olsson2003logarithmic}. The latter is itself an Artin fan. In this case the cone stack $\Sigma(A_{\mathcal X})$ is usually refer to as the \emph{tropicalization} of the logarithmic stack.
As we closely follow the work of Maulik-Ranganathan, this notion is not relevant for us. 
The reason is that by \cite[Section~5]{abramovich2016skeletons}, allowing more flexibility in the choice of strict morphism $\mathcal X_1\to\mathcal A_1$, we can say that for any log-morphism $\mathcal X_2\to\mathcal X_1$ there is an induces log morphism of  Artin fans $\mathcal A_2\to\mathcal{A}_1$ commuting with the strict morphism $\mathcal{X}_2\to\mathcal{X}_1$. 
\end{rem}


\begin{defi}
    Let $\mathcal X$ an algebraic stack with log-structure. We say that $\mathcal X^{\diamond}\to\mathcal X$ is a \emph{tropical model} or a \emph{logarithmic modification} if $\mathcal X^{\diamond}=\mathcal X\times_{\mathcal A}\mathcal A^{\diamond},$ where $\mathcal X\to\mathcal A$ is strict and $\mathcal A^{\diamond}\to\mathcal A$ is \emph{proper birational} map.
\end{defi}
\'Etale locally, $\mathcal A^{\diamond}$ is obtained by performing torus equivariant modifications, corresponding to subdivisions and changes of lattices along torus invariant subvarieties.

\begin{rem}
    As in \cite[Warning~1.4.1]{MaulikRanganathan2025}, we point out that $\mathcal X^{\diamond}\to\mathcal X$ is in general not birational. This is true for example when the map $\mathcal X\to\mathcal A$ is smooth, which is unfortunately never the case for $\mathcal X$ a moduli space of log stable maps in higher genus.
\end{rem}


Via the categorical equivalence  between Artin fans and cones proved in \cite{cavalieri2020moduli},
the proper birational maps $\mathcal A^{\diamond}\to\mathcal A$ correspond to subdivisions
$\Sigma(\mathcal A^\diamond)\to\Sigma. $ We thus sometimes also refer to  $\mathcal X^{\diamond}\to\mathcal X$ as a subdivision of the logarithmic stack.

\medskip

The tropical models of a logarithmic stacks form an inverse system: indeed given $\mathcal A_1\to\mathcal A$, $\mathcal A_2\to\mathcal A$ logarithmic modifications, there always exists an Artin fan $\mathcal A_3$ mapping to both; this is immediate from the subdivision point of view as given two subdivisions of a cone complex or cone stack is always possible to construct a common refinement.

\medskip

Moreover, the category of Artin fans is closed under fiber products. 
This allows us to talk about logarithmic modifications $\mathcal Y^{\diamond}\to\mathcal X^{\diamond}$ of logarithmic morphisms of algebraic stacks.

\subsection{Logarithmic cohomology ring}

\begin{defi}[Definition~2.4]\cite{Holmes2022Logarithmic}
    Let $\mathcal X$ be an algebraic stack with log structure, then the logarithmic (operational) Chow ring of $\mathcal X$, is the colimit of the operational Chow over all logarithmic modifications/tropical models.
    \[\mathrm{logCH}^*(\mathcal X)=\mathrm{colim}_{\mathcal X^\diamond\to\mathcal X}\mathrm{CH}^*(\mathcal X^\diamond).\] 
    Logarithmic cohomology, logarithmic Chow homology and Borel-Moore homology are defined the same way.
\end{defi}

There exists a natural injective ring homomorphism $ \mathrm{CH}^*_{\mathrm{op}}( \mathcal X)\hookrightarrow \mathrm{logCH}(\mathcal X)$ as well as a \emph{group homomorphism}
$\mathrm{logCH}(\mathcal X)\to \mathrm{CH}^*_{\mathrm{op}}( \mathcal X)$ defined via proper push-forward.

\medskip

A \emph{piece-wise polynomial} $f\in\mathrm{PP(\Sigma)}$ on a cone stack $\Sigma$ is a continuous function $f\colon |\Sigma|\to\mathbb R$ (where $|\Sigma|$ denotes the support of $\Sigma$) such that there exists a subdivision $\widetilde{\Sigma}\to\Sigma$ with $f$ polynomial on each cone of $\widetilde \Sigma.$ If no subdivision is necessary we say that $f$ is strict and write $f\in \mathrm{sPP(\Sigma)}.$

As explained in \cite{Holmes2022Logarithmic,molcho2023hodge}, if $\Sigma(\mathcal A)$ is the cone complex associated to an Artin fan $\mathcal A$ and the logarithmic stack $\mathcal X$ comes equipped with a strict morphism $\X\to\A$, we have a ring homomorphism
\[\Phi\colon \mathrm{PP(\Sigma(\mathcal A))}\to \mathrm{logCH}(\mathcal X).\]
We will thus talk about piece-wise polynomials as elements in the logarithmic cohomology of a log-stack.

\medskip

A more detailed description of  $\mathrm{logCH}(\mathcal X)$ for $\mathcal X=\bar{\M}_{g,n}$ endowed with its natural divisorial log structure is given in \cite{holmes2025logDR,pandharipande2024logarithmic}.

\subsection{Virtual birational models}\label{sec:virtual_birational model}

Let $\mathcal X$ a Deligne-Mumford stack with log-structure either defined via a strict map to an Artin fan $\mathcal X\xrightarrow{g}\mathcal A$, or a strict map to a log-smooth algebraic stack $\mathcal X\xrightarrow{g}\mathcal Y$. Assume furthermore that $g$ is endowed with a perfect obstruction theory $\mathbf{E}_g$ giving a virtual class 
$\vir{\mathcal X}=g^![\mathcal A]$ (or $\vir{\mathcal X}=g^![\mathcal Y]$).

\begin{lem}\label{lem:virtualmodels}
    The class $\vir{\mathcal X}$ lifts to a class in $\mathrm{logCH}_*(\mathcal X)$. In other words, any tropical model $\mathcal X^\diamond$ is equipped with a virtual class $\vir{\mathcal X^\diamond}$ pushing forward to $\vir{\mathcal X}$ along the proper map $\mathcal X^\diamond\to \mathcal X.$
\end{lem}

\begin{proof}
    It is an immediate consequence of the definition of tropical model and Costello's virtual push-forward formula \cite{manolache2012virtual,herr2023costello} (see also \cite[Section~3]{ranganathanGWexpansion}).
\end{proof}

\begin{rem}
 Application of this Lemma for specific geometric situations appeared all over the logGW literature. Indeed, allowing oneself to work with different tropical models has been essential to: prove birational invariance and comparisons theorems for logGW
 \cite{abramovich2014comparison, abramovich2018birational}, product formulas in logGW theory, relate invariants of toric varieties with DR-cycles \cite{ranganathan2024logarithmic}, prove the decomposition formula for logGW and logDT \cite{ranganathanGWexpansion,MaulikRanganathan2025}, handle the vanishing cohomology in the degeneration formula \cite{maulik2025gromov}.
\end{rem}

We are interested in virtual birational models and lifts of the virtual class to logarithmic Chow/homology for the following geometries:

\begin{itemize}[leftmargin=0.6cm]
    \item The moduli space  $\M(a):=\M_{g,n+m}(E,a)$ of stable maps to the elliptic curve $E$. It comes with a strict map $g$ to the stack of pre-stable curves $\mathfrak M_{g,n+m}$ and admits the usual perfect obstruction theory \cite{behrend1997intrinsic}. The stack $\mathfrak M_{g,n+m}$ endowed with its divisorial log structure is log-smooth.

    \item The moduli spaces of log-stable maps $\M(A)=\M_{g,m}(A/T,\beta)$. By the discussion in Section \ref{sec-tropical-models}, the  log-structure on this moduli space \cite{gross2013logarithmic} comes by pull-back along a strict map to some Artin fan $\A_{\M}\to\mathcal Log_T.$
    As recalled in \cite[Section~1.3]{MaulikRanganathan2025}, these morphisms are \'etale and representable.
    This implies that the obstruction theory

\[ \mathbf{E}_{\M(A)}:=\left (R\pi_*R\mathcal Hom([F^*\Omega ^{\text{log}}_{A\slash B}\to \Omega^{\text{log}}_{\mathcal C\slash \mathcal M(A)}(x_1+\dots +x_m)] ,\mathcal O)\right)^\vee\]
can be taken to be the obstruction theories relative to $\A_{\M}.$ Moreover, still by \'etaleness and representability, the intrinsic cones for $\M(A)$ are identified $\mathfrak C_{\M(A)/ \mathcal Log _T}\cong\mathfrak C_{\M(A)/\A_{\M}}$. In particular, the latter is embedded in the reduced obstruction bundle (see Lemma~\ref {lem:relativeconered}) and $[\M(A)]^{\mathrm{red}}$ is obtained as virtual pull-back of the fundamental cycle $[\mathcal A_{\M}]$ for the reduced obstruction bundle $\mathbf{E}_{\mathrm{red}}$ of $\M(A)$ over $\mathcal A_{\M}.$ Then Lemma~\ref{lem:virtualmodels} holds for the reduced class.

\item The moduli space $\M_{\Lambda}(X|D)$ of log maps to a normal crossing pair $(X|D)$ with $\lambda$ the necessary discrete data. We more specifically care about $\M_{g,m+n}(E\times\PP^1|E_++E_-; a,\bfw)$.  As before, the map to $\mathcal Log$, which admits a perfect obstruction theory, factors through some Artin fan $\mathcal A_{\M_{\Lambda}(X|D)}$ and thus the same argument as before applies.

\item The moduli space $R\mathcal  M_{g,m+n}(E\times\PP^1|E_++E_-; a,\bfw)$ of rubber maps to   $E\times\PP^1|E_++E_-.$ We refer for example to \cite{janda2020double,marcuswiselog},\cite[Section~6]{bae2023pixton},  at references therein for the definition of the log structure and the virtual class on this space. 
\end{itemize}

\subsection{Exotic insertions}

Given a moduli space of log maps $\M_{\Lambda}(X|D)$ (always keeping in mind $\M_{g,m+n}(E\times\PP^1|E_++E_-; a,\bfw)$), one has a (logarithmic) evaluation morphism to an \textit{evaluation space} $\mathrm{Ev}_{\Lambda}(X|D)$:

\[\ev\colon \M_{\Lambda}(X|D)\to \prod_{i=1}^n W_i=:\mathrm{Ev}_{\Lambda}(X|D),\]

where  $W_i$ is the stratum to which the $i$th-marked point is mapped, endowed with the strata logarithmic structure (i.e. the log structure induced by taking the intersection with those divisor \emph{not} containing $W_i$). 

\begin{expl}
    For $\M(A)$, the evaluation space is simply the $m$-th fiber product of $A_B^m$.
    For $\M(E\times\mathbb P^1|E^++E^-)$, denoting by $m$ the number of interior markings and by $n$ the number of boundary markings, the evaluation space is $(E\times\mathbb P^1)^m\times E^n.$
\end{expl}

Let $\mathrm{Ev}_{\Lambda}(X|D)^{\diamond}$ a tropical model of the evaluation space.\footnote{In this case where the logarithmic structure is simply the one of a product of normal crossing pairs, this is nothing but a strata blow-up.} Consider the fiber product in the category of fs log-schemes (see for example \cite[Section~1.7]{MaulikRanganathan2025} and references therein):

\begin{equation}
\begin{tikzcd}
\M_{\Lambda}(X|D)^{\diamond}\ar[r]\ar[d] &\M_{\Lambda}(X|D)\ar[d]\\
\mathrm{Ev}_{\Lambda}(X|D)^{\diamond} \ar[r] & \mathrm{Ev}_{\Lambda}(X|D).
\end{tikzcd}
\end{equation}
Then $\M_{\Lambda}(X|D)^{\diamond}$ is a tropical model of $\M_\Lambda(X|D)$ and in particular a virtually birational model  of the moduli space of log maps. As explained above, it is endowed with a virtual class pushing forward to the usual one along the top horizontal morphism. 

Crucially, the fiber product in the category of fs log-schemes usually does \emph{not} produce a cartesian diagram in the underlying category of schemes (see for example \cite[Section~1.7]{MaulikRanganathan2025}). 

\begin{defi}\cite{maulik2025gromov}
    An \emph{exotic insertion} is a cohomology class $\gamma^{\diamond}$ on some tropical model $\mathrm{Ev}_{\Lambda}(X|D)^{\diamond}$ of the evaluation space which is \emph{not pulled-back} from $\mathrm{Ev}_{\Lambda}(X|D).$

    An exotic Gromov-Witten invariants is then simply defined as integrating an exotic insertion $\gamma^\diamond$ and a class in $H^*(\bar{\M}_{g,n},\QQ)$ against the virtual class in this model:
    \[\deg \left( \vir{\M_{\Lambda}(X|D)^{\diamond}}\cap (\ev^*\gamma^{\diamond}\cup\mathrm{ft}^*\alpha)\right).\]
\end{defi}
Because the diagram in not cartesian of underlying algebraic stack, exotic invariants cannot be recovered by usual invariants via push-pull.

\section{Multiple cover formula for the correlated (log) DR-cycle}\label{sec:MCFlogcorrelatedDR}

In this section we recall the Definition of the correlated DR-cycle introduced in \cite{blommecarocci2025DR}, focusing on the specific case where the target is an elliptic curve $E$. We denote the usual DR-cycle by and $\DR(a,\bfw)=\epsilon_* \vir{R\M(a,\bfw)}$, where $R\M(a,\bfw)$ is the moduli space of rubber stable maps to $E\times\PP^1|E_0+E_{\infty}$ of  
 degree $a$ on $E$  (i.e. $\pi_*f_*[C]=a[E]$ for $\pi\colon E\times\PP^1\to E$) and  with ramification profile $\bfw=(w_1,\dots,w_n)$. We conclude the section with the proof of the MCF for the correlated (log) GW-cycles obtained from the logDR-cycle, i.e. Theorem \ref{theo:MCF-for-correlated-DR}.

\subsection{Recollections in the non-correlated DR-cycle formula}

\subsubsection{Pixton's formula} The usual DR-cycle with target variety $\DR(a,\bfw)$ (see \cite{janda2020double}) has an expression as a sum of decorated strata classes, known as Pixton's formula:
\begin{equation}\label{eq:DRtarget}
  \DR(a,\bfw) = \sum_\Gamma \vir{\M_\Gamma(\bfa)}\cap\operatorname{Dec}_\Gamma(\bfw).  
\end{equation}

In the notation we forget about the target, genus and markings; we record instead the curve class $a$ as well as the tangency profile. 
\begin{itemize}[leftmargin=0.6cm]
    \item The sum is indexed by so-called $E$-stable graphs (see \cite[Section 0.3]{janda2020double}), which are graphs whose vertices are decorated by a curve class $a_v$ and a genus $g_v$ such that the total genus $b_1(\Gamma)+\sum g_v=g$. We denote by $\bfa=(a_v)$ the collection of curve classes indexed by vertices.
    \item The stratum $\M_\Gamma(\bfa)$ is the fiber product over the diagonal $\Delta\colon E^{\sfE(\Gamma)}\to E^{\sfH(\Gamma)}$ of  the product moduli spaces $\prod_v\M_v(a_v):=\M_{g_v,n_v+m_v}(E,a_v)$ mapping to $E^{\sfH(\Gamma)} $ via the evaluation maps. It has a virtual class $\Delta^!(\prod_v[\M_v(a_v)]^{\text{vir}})$ induced by Gysin pull-back.
    We denoted by $n_v$ the valency of the vertex $v$, i.e. the flags of $\Gamma$ based at $v;$ $m_v$ is just the number of interior markings adjacent to $v$ and do not play a role in the gluing. 
    \item The decoration $\Dec_\Gamma(\bfw)$ is given explicitly in Pixton's formula \cite{janda2020double} and is a polynomial in the entries $w_i$ whose coefficients are cohomology classes pulled-back from the moduli space of curves $\prod\overline{\M}_{g_v,n_v+m_v}$. In the case at hand, since we are looking at the trivial bundle, these decorations only come from products of $\psi$-classes. We refer to \cite{janda2020double} for the explicit formula as we do not need the precise shape of the decoration. 
\end{itemize}

\subsubsection{Logarithmic double ramification cycle with target variety}

Before correlating the Pixton formula, we need an extension of the DR-cycle with target variety to the log-Chow ring to be later able to deal with the exotic insertions. We thus introduce this lift.

\medskip

After the seminal work in \cite{bae2023pixton}, we can see Pixton's formula for DR with target varieties\footnote{We remark that \cite{janda2020double} is an ingredient for the proof of the main result in \cite{bae2023pixton} not merely a special case.} as an application of the universal DR formula, indeed given $\mathcal L$ a line bundle on $X$ one obtains a morphism
  \[\M_{g,n}(X,\beta)\to\mathfrak Pic_{g,n,d}\]
  defined by $(f\colon C\to X)\mapsto f^*L.$ Here we denoted by $d:=c_1(\mathcal L)\cdot\beta.$

  For the special case of $\L$ being the trivial line bundle, and the map $\M_{g,n}(X,\beta)\xrightarrow{\mathrm{ft}}\bar{\M}_{g,n}$ does not contract any component of the source curve,  which is the only case of interest for us, the morphism to $\mathfrak Pic_{g,n,0}$ actually factors trough $\bar{\M}_{g,n}$ and thus \eqref{eq:DRtarget} simply follows  rewriting 

  \[ \vir{\M(a)}\cap \mathrm{ft}^*  \DR^{\mathrm{op}}_{\bar{\M}_{g,n}}(a,\bfw)\]

  using the splitting axiom for the virtual class, where $\DR^{\mathrm{op}}_{\bar{\M}_{g,n}}(a,\bfw)$ denotes the pull-back to $\bar{\M}_{g,n}$ of the universal DR-operational class defined and computed in \cite{bae2023pixton}. 
  
  In this situation we simply define:

\begin{defi}\label{defi:triviallogDR}
    Let $\mathrm{logDR}_{\bar{\M}_{g,n}}(a,\bfw)\in\mathrm{logCH}^g(\bar{\M}_{g,n})$ the logarithmic double ramification cycle defined in \cite{holmes2021extending,Holmes2022Logarithmic} and computed in \cite{holmes2025logDR}. We define

    \[ \mathrm{logDR}(a,\bfw):=\vir{\M(a)}\cap\mathrm{ft}^*\mathrm{logDR}_{\bar{\M}_{g,n}}(a,\bfw)\in \mathrm{logCH}(\M(a)) \]
    The product makes sense since $\vir{\M(a)}$ lifts to a class in $ \mathrm{logCH}(\M(a))$ by Lemma~\ref{lem:virtualmodels}.
\end{defi}

\begin{rem}
    The above definition is the correct one only in the restrictive hypothesis that $\mathcal L$ is trivial and there is no stabilization in the forgetful map. More in general, a representative of the log DR class with target varieties (at least under the assumption that there is a marking) should come applying the methods of \cite[Section~4,5]{holmes2025logDR} to the moduli space of stable maps endowed with the log structure pulled-back from pre-stable curves.
\end{rem}

 By \cite{holmes2025logDR}, there exists an explicit logarithmic modification $\bar{\M}^\diamond\to\bar{\M}_{g,n}$  and an explicit twist $\mathcal L_{\bfw}$ of $\mathcal O_{\mathcal C^\diamond}(-\sum w_ix_i)$  by a strict-piecewise linear function (in the sense of \cite[Section~2.5]{holmes2025logDR} and references therein), both depending on  $\bfw$ such that $\mathrm{logDR}(a,\bfw)$ is obtained by applying the universal DR-formula \cite{bae2023pixton} to

 \[\bar{\M}^\diamond\xrightarrow{\varphi_{\mathcal L_{\bfw}}} \mathfrak Pic_{g,n,0}. \]

 The strata of the logarithmic modification $\bar{\M}^\diamond$ are indexed by certain discrete data described in the remark below; we simply denote by $\Gamma^\diamond$ the collection of these discrete data.

\begin{rem}
The modification and the twist are constructed using stability conditions for line bundle on quasi-stable curves as explained in \cite[Section~4]{holmes2025logarithmic}. The logarithmic strata of $\bar{\M}^\diamond$ are indexed by $\Gamma^\diamond=(\Gamma_{\bfw}^{\sigma},D_{\bfw}^\sigma,\alpha^{\sigma})$ where $\Gamma_{\bfw}^{\sigma}$ is a quasi-stable model for $\Gamma$ (i.e. obtained adding at most one bi-valent vertex on each edge); $D_{\bfw}^\sigma$ a $\sigma$-stable tropical divisor equivalent to $\underline{\deg}(\mathcal O(-\bfw))$ and $\alpha^{\sigma}$ the piecewise linear function realizing the equivalence between the two divisors. The exact description of these strata is not necessary.
\end{rem}

It then follows that 

\[ \mathrm{logDR}(a,\bfw)=\vir{\M(a)^\diamond}\cap \mathrm{ft}^* \varphi_{\mathcal L_{\bfw}}^*\DR^{\mathrm{op}} .\]

In particular, it admits a tautological expression as in \eqref{eq:DRtarget}. In the expression,  $\vir{\M_{\Gamma}(a)}$ is replaced by $\vir{\M^\diamond_{\Gamma^\diamond}(a)}$, where $\bar{\M}_{\Gamma^\diamond}$ is (an \'etale torsor over the normalization of the closure of) a logarithmic stratum in $\bar{\M}^\diamond$ (see  \cite[Section~6]{holmes2025logDR} for details).  

\medskip

The virtual class $\vir{\M^\diamond_{\Gamma^\diamond}(a)}$ has two possible descriptions. Each strata $\bar{\M}_{\Gamma^\diamond} $ comes with a morphism (not necessarely birational, it could have positive dimensional fibers) $\bar{\M}_{\Gamma^\diamond} \xrightarrow{p_{\Gamma^\diamond}} \bar{\M}_{\Gamma}$, which is l.c.i and thus admits a Gysin pull-back. Then, the virtual class may be obtained as Gysin pull-back along either of the cartesian diagrams

\begin{equation*}
\begin{tikzcd}
\M^\diamond_{\Gamma^\diamond}(a)\ar[r]\ar[d] &\M^\diamond(a)\ar[d]\\
\bar{\M}_{\Gamma^\diamond} \ar[r] & \bar{\M}^\diamond
\end{tikzcd} ,\quad \text{ or }
\begin{tikzcd}
\M^\diamond_{\Gamma^\diamond}(a)\ar[r]\ar[d] &\M_{\Gamma}(a)\ar[d]\\
\bar{\M}_{\Gamma^\diamond} \ar[r, "p_{\Gamma^\diamond}"] & \bar{\M}_{\Gamma}
\end{tikzcd},
\end{equation*}
the first one endowing the virtual strata of the tropical model $\M^\diamond(a)$ with its natural class $\vir{\M^\diamond_{\Gamma^\diamond}(a)}.$
By functoriality of the virtual pull-back, and an easy diagram chase, we see that  \begin{equation} \label{eq:gysinpull-back}
    p_{\Gamma^\diamond}^!\vir{\M_{\Gamma}(a)}=\vir{\M^\diamond_{\Gamma^\diamond}(a)}
\end{equation}
The equation allow us to relate the virtual class of the strata on the tropical model with  $\vir{\M_{\Gamma}(a)}$, which we know possess a description in terms of $\prod\M_v(a_v)$ capped with the diagonal.

\medskip

Furthermore, these classes are always compatible with splittings of the space of maps in open and closed components as those considered in \cite[Section~4]{blommecarocci2025DR} and recalled below.

\medskip

This allows to transfer in a straithforward manner all the computations of the next section to any tropical model of $\M(\bfa)$ and yielding a lift of the correlated DR formula and of the MCF formula for the latter in log-homology.

\subsubsection{Goal}

Given a geometric insertion $\gamma\in H^*(E^{n+m},\QQ)$, we can cap with the log DR-cycle and push to the log-homology of $\bar{\M}_{g,n+m}$, i.e. consider
$$\gen{\gamma}_{a,\bfw} = \ft_*( \mathrm{log}\DR(a,\bfw)\cap \ev^*\gamma) \in \mathrm{log}H^*(\bar{\M}_{g,m},\QQ),$$
where $\ev$ and $\mathrm{ft}$ are respectively the evaluation map and the forgetful morphism to the moduli space of curve. For fixed insertions, this defines a function on $\RRR$, the $\NN$-module of ramification profiles.

\medskip

Defining a correlated (log) DR-cycle $\mathbf{logDR}(a,\bfw)$ and using the correlated Pixton's formula, we get a correlated refinement $\gen{\gen{\gamma}}_{a,\bfw}$, which can be seen as a diagonal type function $\RRR\to \mathrm{log}H^*(\bar{\M}_{g,n+m},\QQ)\otimes\G$. We then prove that this refinement satisfies the MCF.

\medskip

The correlated Pixton formula involves classes of boundary strata of a refined stratification of $\M(a)=\M_{g,n}(E,a),$ which we recall  below for the special case of target the elliptic curve.


\subsection{Ingredients for the correlated refinement}

In this section we recall the ingredients needed to prove the correlated refinement of the DR with target variety formula obtained in \cite{blommecarocci2025DR} and explain the (straightforward) lift to a formula in the logarithmic Chow ring.





\subsubsection{Refined stratification of $\M_{g,n}(E,a)$}
\label{sec-refined-stratification}


\begin{defi}\cite[Definition 4.23]{blommecarocci2025DR}
A \emph{monodromy graph}  $(\Gamma,K,\varphi)$ is the data of: $\Gamma$ a $E$-valued stable graph, $K$  a subgroup of $\Pic^0(E)[|\bfw|]$ and a morphism $\varphi\colon K\to H^1(\Gamma,\mu_{|\bfw|})$, with $\mu_{|\bfw|}$ the $|\bfw|$-roots of unity.  
\end{defi}

The stratum indexed by $(\Gamma,K,\varphi)$ is the substack of $\M_\Gamma(\bfa)$ parametrizing stable maps $f\colon C\to E$ such that
\begin{enumerate}[label=(\roman*)]
    \item the dual graph of $C$ decorated with the degrees $a_v$ of the restrictions of $f\colon C\to E$ to the components $C_v$  and genus at $g_v$ at the vertices is  the $E$-valued stable graph $\Gamma$;
    \item the subgroup of $\Pic^0(E)[|\bfw|]$ of componentwise trivial line bundles ($\L$ such that $f\rvert_{C_v}^*\L\cong\O_{C_v}$ for each $v\in\sfV(\Gamma)$) is $K$;
    \item for a line bundle $\L\in K$, its pull-back $f^*\L\in H^1(\Gamma,\mu_{|\bfw|})$ is thus determined by its monodromy along loops of $\Gamma$, and prescribed by $\varphi$: for each a loop $\gamma\in H_1(\Gamma,\ZZ)$, the monodromy of $\L$ along $\gamma$ is $\varphi(\L)(\gamma)$.

\end{enumerate}

We denote by $\varphi^\intercal\colon H_1(\Gamma,\ZZ_{|\bfw|})\to \widehat{K}$ the adjoint of $\varphi$, where $\widehat{K}=\Hom(K,\mu_{|\bfw|})$ is the Pontrjagin dual of $K$.

\begin{rem}
    In particular, $\operatorname{Ker}(\varphi)\subseteq K$ is the subgroup of globally trivial line bundles: $f^*\L\cong\O_C$.
    This kernel is called the \emph{core of the stratum}; when needed, we denote it by $C$. In \cite{blommecarocci2025DR}, the notation for $(K,C)$ is $(\tK,K)$. As the core does not play any role in the present paper, we drop the tilde from the notation of the kernel.
\end{rem}

As proved in \cite{blommecarocci2025DR}, the stratum $\M_\Gamma(\bfa)$ splits as a disjoint union of open and closed sub-stacks indexed by the possible $(K,\varphi)$. In particular, the smooth stratum of $\M(a)$ splits as a disjoint union of open and closed sub-stacks indexed by subgroups $K\subset\Pic^0(E)[|\bfw|]$, yielding a  natural decomposition of its virtual class:
$$\M(a) = \bigsqcup_K \M_K(a), \quad \vir{\M(a)} = \sum_K \vir{\M_K(a)}.$$

\subsubsection{Description of the smooth stratum}
We recall the description of $\M_K(a)$ which allows us to compute its virtual class.

\medskip

\textbf{$\bullet$ Coverings of $E$.} To each subgroup $H\subset\Pic^0(E)[|\bfw|]$ is associated a degree $|H|$ normal covering
$$\pi_H\colon E_H\to E,$$
 The deck-transformation group of $\pi_H$ is the Pontrjagin dual of $H$, i.e. $\widehat{H}=\Hom(H,\ZZ_{|\bfw|})$.
  See \cite[Section 4.1]{blommecarocci2025DR} for an explicit construction.

\medskip

\textbf{$\bullet$ Map between moduli spaces.} 
A map $f\colon C\to E$ lifts to $E_H$ if and only if the pull-back along $f$ of any line bundle in $H$ is trivial. A family of maps $F\colon \C/S\to E$ such that, for each $\L\in H, F^*\L$ is trivial fiberwise, lifts to a family of maps to $E_H$ up to base change along a $|H|$-covering $S'\to S.$ If a map lifts, the set of lifts is a torsor under $\hat{H}.$

 Using this, it is proven in \cite[Section 4]{blommecarocci2025DR},  that composition with the covering defines a degree $|H|$ map
$$\pi_H\colon \M(E_H,a/|H|)\longrightarrow \bigsqcup_{K\supset H}\M_K(E,a)=\M^H(E,a)\subset\M(E,a).$$
The equality on the right is the definition of $\M^H(E,a)$.
The moduli space $\M(E_H,a/|H|)$ is understood to be empty if $|H|$ does not divide $a$.

In particular, the virtual class of the components $\M_K(a)$ are related, up to inclusion-exclusion relations, to the virtual classes of moduli spaces of coverings of $E$. 
Using the M\"obius function on the lattice of subgroups of $\Pic^0(E)[|\bfw|]\simeq\ZZ_{|\bfw|}^2$, we can give an explicit expression.

\begin{prop}
\label{prop:relation-to-covering-smooth-case}
\begin{enumerate}[label=(\roman*)]
    \item If $H\subset\Pic^0(E)[|\bfw|]$, we have
    $$\sum_{K\supset H}\vir{\M_K(E,a)} = \vir{\M^H(E,a)} = \frac{1}{|H|}\pi_{H\ast}\vir{\M(E_H,a/|H|)}.$$
    \item If $K\subset\Pic^0(E)[|\bfw|]$, we have
    $$\vir{\M_K(E,a)} = \sum_{H\supset K}\mu(H/K)\frac{1}{|H|}\pi_{H\ast}\vir{\M(E_H,a/|H|)},$$
    where $\mu$ is the M\"obius function on the lattice of subgroups of $\ZZ_{|\bfw|}^2$.
    \item If $\gamma\in H^*(E,\QQ)^{\otimes n}$, for $K\subset\Pic^0(E)[|\bfw|]$, we have
    $$\vir{\M_K(E,a)}\cap\ev^*\gamma = \sum_{H\supset K}\mu(H/K)\frac{1}{|H|}\pi_{H\ast}\left( \vir{\M(E_H,a/|H|)}\cap\ev^*\pi_H^*\gamma\right).$$
\end{enumerate}
\end{prop}

\begin{proof}
    The first point comes from the degree $|H|$ morphism and is \cite[Lemma 4.13]{blommecarocci2025DR}. The second is obtained performing a M\"obius inversion on the lattice of subgroups. The last point comes from push-pull-formula.
\end{proof}

\medskip

\textbf{$\bullet$ Transferring.} 
To compute the integrals that appear in the Proposition above, we need to describe the action of $\pi_H^*$ on cohomology, using that all $E_H$ are actually elliptic curves.

\medskip

If $\pt$ is the Poincar\'e dual class to a point, we have $\pi_H^*\pt=|H|\pt$. In order to describe the action on $1$-cycles, we use the following lemma.

\begin{lem}\label{lem:pull-back-covering}
    Consider a multilinear function $\varphi\colon(\QQ^2)^{\otimes 2r}\to\QQ$ invariant by the diagonal action of $SL_2(\ZZ)$ on the source. Then, if $M$ is any matrix with $\QQ$-coefficients, we have
    $$\varphi(Mx_1,\dots,Mx_{2r}) = (\det M)^r\varphi(x_1,\dots,x_{2r}).$$
\end{lem}

\begin{proof}
    The formula is true for $SL_2(\ZZ)$ by assumption. Therefore, it is also true for $SL_2(\QQ)$, which is the Zariski closure of $SL_2(\ZZ)$ in $GL(2,\mathbb Q)$. The formula is also true for scalar matrices by multilinearity. These generate the subgroup of $GL_2(\QQ)$ whose determinant is a square in $\QQ$, which is also Zariski dense in $\M_2(\QQ)$. The result follows.
\end{proof}

Let $E_H\to E$ be a covering  and $\gamma\in H^*(E,\QQ)^{\otimes n}$. As both $E$ and $E_H$ are elliptic curves, $\gamma$ makes sense both as insertion for $E$ and $E_H$ in the sense that we always have an isomorphism $s\colon H^*(E,\QQ)\to H^*(E_H,\QQ)$ with $s\in SL_2(\ZZ)$. the precise isomorphism does not matter due to invariance. By abuse of notation, we write $\gamma$ for the cohomology class in $E_H$ obtained as image of $\gamma$ along the isomorphism $s$; i.e. on the right-hand-side of the following corollary, $\gamma$ means the class obtained by deformation from $E$ to $E_H$ and not the pull-back.

\begin{coro}\label{coro:pull-back-covering}
    Let $\pi_H\colon E_H\to E$ be a degree $|H|$ covering of elliptic curves, and take
    $\gamma\in H^1(E,\QQ)^{\otimes 2r}\otimes H^2(E,\QQ)^{\otimes s}$.
    We then have in $H^*(\bar{\M}_{g,n},\QQ)$
    $$\ft_*\left( \vir{\M(E_H,a/|H|)}\cap\ev^*(\pi_H^*\gamma) \right) = |H|^{r+s} \cdot \ft_*\left( \vir{\M(E_H,a/|H|)}\cap\ev^*\gamma \right).$$
    In particular, dropping $E$ from the notation and for $K\subset\Pic^0(E)[|\bfw|]$, we have
    $$\ft_*(\vir{\M_K(a)} \cap \ev^*\gamma) = \sum_{H\supset K} \mu(H/K)|H|^{r+s-1} \cdot \ft_*(\vir{\M(a/|H|)}\cap\ev^*\gamma).$$
\end{coro}

\begin{proof}
    After applying $\ft_*$, deformation invariance ensures that seen as cycle-valued functions on $H^*(E,\QQ)^{\otimes n}$, they are $SL_2(\ZZ)$-invariant. Applying Lemma \ref{lem:pull-back-covering}, the pull-back $\pi_H^*$ acts on $H^1(E,\ZZ)$ via a linear map of determinant the degree of the covering $|H|$, thus yielding the first identity.

    For the second part, we use Proposition \ref{prop:relation-to-covering-smooth-case}(iii). We can drop the elliptic curves from the notation again since we have no $\pi_H$ appearing anymore.
\end{proof}

\begin{rem}
    In particular, we see that $\ft_*(\vir{\M_K(a)}\cap\ev^*\gamma)$ is invariant when applying the action of $SL_2(\ZZ)$ on $\gamma$, which is far from obvious since when deforming $E$, one should also act on $K$.
\end{rem}

\subsubsection{Description of the boundary strata}

We recalled above that the strata $\M_{\Gamma}(\bfa)$ are described as fiber product of  $\prod_v\M_v(a_v)$ with the diagonal $E^{\sfE(\Gamma)}\xrightarrow{\Delta}E^{\sfH(\Gamma)}$ over the evaluation space $E^{\sfH(\Gamma)}.$ We provide an analogous description for the strata of $\M_K(a).$

\medskip

\textbf{$\bullet$ Refined evaluation map.} 
Since $\M^H(a)$ parametrizes maps that lift to $E_H$ (up to a base change) and the set of lifts is then a $\widehat{H}$-torsor. We have an evaluation map
$$\ev\colon\M^H(E,a)\to E_H^n/\widehat{H},$$
where $\widehat{H}$ acts diagonally by deck-transformation. 
This is a  \emph{rubber evaluation space}, in a sense made precise by \cite[Lemma~4.1]{blommecarocci2025DR}.
Notice that $E_H^n/\widehat{H}\to E^n$ is a covering of degree $|H|^{n-1}$.

\medskip

\textbf{$\bullet$ Twisted diagonals.} Let $\Gamma$ be a $E$-stable graph with set of edges (resp. half-edges) $\sfE(\Gamma)$ (resp. $\sfH(\Gamma)$), and $K$ a subgroup of $\Pic^0(E)[|\bfw|]$.

We consider the pull-back $\operatorname{\Delta}_K$ of the diagonal 
$E^{\sfE(\Gamma)}\to E^{\sfH(\Gamma)}$ along the covering map $E_K^{\sfH(\Gamma)}/\widehat{K}^{\sfV(\Gamma)}\to E^{\sfH(\Gamma)}$.
Following \cite[Lemma 4.20]{blommecarocci2025DR}, $\operatorname{\Delta}_K$ is disconnected;  its $|K|^{b_1(\Gamma)}$ connected components, called \textit{twisted diagonals}, are indexed by a choice of $\varphi\in H^1(\Gamma,\widehat{K})$,  equivalently morphisms $\varphi\colon K\to H^1(\Gamma,\ZZ_{|\bfw|})$. We denote the twisted diagonal associated to $\varphi$ by $\Delta_{K,\varphi}$.

We notice that twisted diagonals are complex subtori  that differ by the action of deck-automorphism groups, which act by translation. In particular, all twisted diagonals are cobordant and realize the same homology class.

\medskip

\textbf{$\bullet$ Fiber product.} Let $\M_v(a_v)=\M_{g_v,n_v+m_v}(E,a_v)$ be the moduli space associated to a vertex of $\Gamma$. Given the splitting in open and closed sub-stacks $\M_v(a_v)=\bigsqcup_{K_v}\M_{v,K_v}(a_v)$, the product over vertices splits as follows:
$$\prod_v\M_v(a_v) = \bigsqcup_{K}\left(\bigsqcup_{\cap K_v=K}\prod_v\M_{v,K_v}(a_v)\right),$$
where the second disjoint union is over the collection of subgroups $(K_v)$ of $\Pic^0(E)[|\bfw|]$ satisfying $\bigcap_v K_v=K$. To make the notation easier, we set:

$$\P_\Gamma(\bfa) = \prod_v\M_v(a_v),\quad \text{ and }\P_{\Gamma,K}(\bfa)=\bigsqcup_{\cap K_v=K}\prod_v\M_{v,K_v}(a_v),$$
with $\P$ as in ``product'' 
so that we can write 
\[\P_\Gamma(\bfa)=\bigsqcup_K\P_{\Gamma,K}(\bfa).\]

These spaces are all endowed with their natural virtual classes, defined by taking product and sums of $\vir{\M_{v,K_v}(a_v)}$ simply obtained restricting the usual obstruction theory to this open (and closed) substack.

\medskip

We can then consider the fiber product

\bcd
\bigsqcup_{\varphi}\M_{\Gamma,K,\varphi}(\bfa) \ar[r]\ar[d] & \P_{\Gamma,K}(\bfa) \ar[d,"\ev"] \\
\bigsqcup_\varphi\Delta_{K,\varphi} \ar[r] & E_K^{\sfH(\Gamma)}/\widehat{K}^{\sfV(\Gamma)}. \\
\ecd

The bottom arrow is the (regular) embedding  of the union of twisted diagonals indexed by $\varphi\colon K\to H^1(\Gamma,\ZZ_{|\bfw|})$. 
Via the Gysin pull-back $\Delta_{K,\varphi}^!,$ each component
 $\M_{\Gamma,K,\varphi}(\bfa)$ is endowed with a virtual fundamental 
 $\Delta_{K,\varphi}^!(\vir{\P_{\Gamma,K}(\bfa)}).$
\medskip

\textbf{$\bullet$ Relation to moduli of stable maps to covering.} Let $\pi_H\colon E_H\to E$ be a covering and $\Gamma$ be a $E$-valued graph.

The pull-back of the diagonal $E^{\sfE(\Gamma)}\to E^{\sfH(\Gamma)}$ along the covering $E_H^{\sfH(\Gamma)}\to E^{\sfH(\Gamma)}$ is also disconnected, indexed by elements $\psi\in C^1(\Gamma,\widehat{H}),$ i.e. a choice of deck automorphism for each edge of the graph.  (See 
\cite[Lemma 4.20]{blommecarocci2025DR} for further details).
We also call these components \textit{twisted diagonals}, and denote them by $\widetilde{\Delta}_{H,\psi}$. 

We can then consider the fiber product
\bcd
\bigsqcup_\psi \M_{\Gamma,\psi}(E_H,\bfa/|H|) \ar[r]\ar[d] & \P_\Gamma(E_H,\bfa/|H|), \ar[d,"\ev"] \\
\bigsqcup_\psi\widetilde{\Delta}_{H,\psi} \ar[r] & E_H^{\sfH(\Gamma)} \\
\ecd
where the disjoint union on the left is over $\psi\in C^1(\Gamma,\widehat{H})$. Choosing $\psi=0$, we recover the usual diagonal and thus  the usual stratum $\M_\Gamma(E_H,\bfa/|H|)$ of the moduli space $\M(E_H,a)$. In this situation as well, for the same argument as before, the different twisted diagonals are all cobordant.

\medskip

For different choices of $\psi$, $\M_{\Gamma,\psi}(E_H,\bfa/|H|)$ should be thought as parametrizing maps that glue up to a choice of deck-automorphisms given by $\psi$.


\medskip

Composition by $E_H\to E$ yields the following commutative diagram, where vertical maps have degree $|H|^{|\sfV(\Gamma)|}$:
\bcd
\bigsqcup_\psi \M_{\Gamma,\psi}(E_H,\bfa/|H|) \ar[r]\ar[d] & \P_\Gamma(E_H,\bfa/|H|) \ar[d] \\
\bigsqcup_{K\supset H}\bigsqcup_\varphi \M_{\Gamma,K,\varphi}(E,\bfa) \ar[r] & \bigsqcup_{K\supset H}\P_{\Gamma,K}(E,\bfa).
\ecd

The image of  $\M_{\Gamma,\psi}(E_H,\bfa/|H|)$ along the left vertical map is contained in the locus where $\varphi|_H=[\psi]\in H^1(\Gamma,\widehat{H}).$ 

\begin{prop}
\label{prop:relation-to-covering-nodal-case}
The virtual classes enjoy the following properties.
    \begin{enumerate}[leftmargin=0.6cm]
        \item Enumerative invariants arising from cap product against the  classes
        $$\vir{\P_\Gamma(E_H,\bfa/|H|)}\cap\ev^*\widetilde{\Delta}_{H,\psi} \quad (\text{resp.} \; \P_{\Gamma,K}(E,\bfa)\cap\ev^*\Delta_{K,\varphi}),$$
        do not depend on $\psi$ (resp. $\varphi$). 
        We thus replace the index $\psi$ (resp. $\varphi$) with $0$ in what follows.
        \item We have the following relation in the homology with $\QQ$-coefficients of $\P_{\Gamma}(E,\bfa)$:
        $$|H|^{b_1(\Gamma)-1}\pi_{H\ast}\left(\vir{\P_\Gamma(E_H,\bfa/|H|}\cap\ev^*\widetilde{\Delta}_{H,0}\right) = \sum_{K\supset H}|K|^{b_1(\Gamma)}\vir{\P_{\Gamma,K}(E,\bfa)}\cap\ev^*\Delta_{K,0},$$
    \end{enumerate}
\end{prop}

\begin{proof}
    \begin{enumerate}[leftmargin=0.6cm]
        \item By definition, the push-forward of $\widetilde{\Delta}_{H,\psi}^!\vir{\P_\Gamma(E_H,\bfa/|H|}$ in $\P_\Gamma(E_H,\bfa/|H|)$ is
        $$\vir{\P_\Gamma(E_H,\bfa/|H|}\cap\ev^*\widetilde{\Delta}_{H,\psi},$$
        still writing $\widetilde{\Delta}_{H,\psi}$ for the Poincar\'e dual class of the diagonal. As twisted diagonals are all cobordant, we conclude.

        \item As the classes we get in $\P_\Gamma(E_H,\bfa/|H|)$ after capping with the diagonal do not depend on $\psi$, we have that the the push-forward of $\sum_\psi\widetilde{\Delta}_{H,\psi}^!\vir{\P_\Gamma(E_H,\bfa/|H|})$ is simply
        $$|H|^{|\sfE(\Gamma)|}\vir{\P_\Gamma(E_H,\bfa/|H|}\cap\ev^*\widetilde{\Delta}_{H,0},$$
        where $|H|^{|\sfE(\Gamma)|}$ is the number of possible $\psi\in C^1(\Gamma,\widehat{H})$. 
        
        Reasoning in a similar way on $\P_{\Gamma,K}(E,\bfa)$, for each $K$ there are $|K|^{b_1(\Gamma)}$ possible choices of $\varphi$. As the classes do not depend on the choices we get:
        \begin{align*}
         & |H|^{|\sfE(\Gamma)|}\pi_{H\ast}\left(\vir{\P_\Gamma(E_H,\bfa/|H|}\cap\ev^*\widetilde{\Delta}_{H,0}\right) \\
         = & |H|^{|\sfV(\Gamma)|}\sum_{K\supset H}|K|^{b_1(\Gamma)}\vir{\P_{\Gamma,K}(E,\bfa)}\cap\ev^*\Delta_{K,0},
        \end{align*}
        where the $|H|^{|\sfV(\Gamma)|}$ factor is due to the degree of the map being $|H|^{|\sfV(\Gamma)|}$. 
    \end{enumerate}
\end{proof}

\begin{rem}
Alternatively we can restrict to one component $\M_{\Gamma,\psi}(E_H,\bfa/|H|)$ mapping to its image. As there are $|H|^{|\sfV(\Gamma)|-1}$ elements $\psi$ in the same class, the degree is only $|H|$. However, given $[\psi]|_H$, there are $|K/H|^{b_1(\Gamma)}$ possible extensions for each $K\supset H$.
\end{rem}

Stress that the push-forwards of the classes in $\M_\Gamma(\bfa)$ are not equal since they belong to distinct components $\M_C(a)$, where $C$ is the core of the monodromy graph.

\medskip

\textbf{$\bullet$ Transferring.} As in the smooth case, we can compute against the virtual class of $\M_{\Gamma,K,\varphi}(E,\bfa)$ 
passing to moduli spaces of maps to $E_H.$

\begin{prop}
    Let $K$ be a subgroup 
    and $\gamma\in H^1(E,\QQ)^{\otimes 2r}\otimes H^2(E,\QQ)^{\otimes s}$. 
    We have the following identity of cycles in $H^*(\bar{\M}_{g,n},\QQ)$:
    
    \begin{align*}
     & |K|^{b_1(\Gamma)} \cdot \ft_*\left( \vir{\P_{\Gamma,K}(\bfa)} \cap\ev^*(\Delta_{K,0}\cup\gamma) \right) \\
     = & \sum_{H\supset K} \mu(H/K)|H|^{b_1(\Gamma)+r+s-1} \cdot \ft_*\left( \vir{\P_{\Gamma}(\bfa/|H|)}\cap\ev^*( \widetilde{\Delta}_{H,0}\cup\gamma) \right).
    \end{align*}
\end{prop}

\begin{proof}
    We start from Proposition \ref{prop:relation-to-covering-nodal-case}(2). We then perform a M\"obius inversion. Using push-pull formula we then proceed as in the smooth case to replace $\pi_H^*\gamma$ by $\gamma$, where we do the same abuse of notation considering $\gamma$ as an insertion both for $E$ and $E_H$.
\end{proof}

\subsection{Correlated DR-cycle formula}
\label{sec-correlated-DR}

We denote by  $\M(a)(\frac{1}{|\bfw|})$ the moduli of \textit{spin stable maps} parametrising a stable map together with a (log-)root of the trivial bundle of order $|\bfw|$. 
Up to taking a root stack, obtained as explained in \cite[Section~5]{blommecarocci2025DR} and references therein, the universal log root $\T$ on the moduli space of spin stable maps can be represented by an honest line bundle $\O_{\C}(D)$ for $D$ a divisor supported on the boundary, whose tropicalization is a tropical torsion divisor $\ttD$ on the dual graph.

\medskip

 The space $\M(a)(\frac{1}{|\bfw|})$  also splits as a union of open and closed substack and we thus have a refined stratification of its boundary.
The refined stratification essentially coincides with the stratification of $\M(a)$ described above, up to the additional information of the tropical torsion divisor $\ttD$ on the dual graph. Strata are hence indexed by $(\Gamma,K,\varphi,\ttD)$ with  $(\Gamma,K,\varphi)$ as above and $\ttD$ a tropical $|\bfw|$-torsion divisor.

\medskip

The correlated DR-cycle formula in \cite{blommecarocci2025DR} is obtained applying the universal DR-formula \cite{bae2023pixton} on 
the \emph{correlated components} of the moduli space of spin maps and then pushing forward to $\M(a)$, which requires a degree computation.


\subsubsection{Logarithmic Weil pairing}
In order to talk about the correlated components of  $\M(a)(\frac{1}{|\bfw|})$, we recall the definition of the logarithmic Weil pairing and of \emph{correlator} for a torsion line bundle. For this part we assume the reader to be familiar with the Logarithmic Picard group of Molcho-Wise \cite{molchowiselog}, and more in general with the theory of (logarithmic) abelian varieties.

\medskip

Let $\C\to \M(a)(\frac{1}{|\bfw|})$ be the universal curve; this is endowed with a natural vertical logarithmic structure coming from the moduli space of spin curves (we refer to \cite{holmes2023root} and references therein; a description cal also be found in \cite[Section~5]{blommecarocci2025DR}).

\medskip

We denote by $\LogPic^0( \C/\M(a)(\frac{1}{|\bfw|}))$ the logarithmic Jacobian introduced in \cite{molchowiselog} and by $\LogPic^0( \C/\M(a)(\frac{1}{|\bfw|}))[|\bfw|]$ the subgroup of $|\bfw|$-torsion log line bundle. After a base change along the the root stack map necessary to have a line bundle representative of the universal root, it is shown to be an \'etale finite group scheme over $\M(a)(\frac{1}{|\bfw|})$ \cite[Theorem~4.3]{holmes2023root}. We collect the facts we need in the following Proposition; the proofs of the various statements can be found in \cite[Sections~2-3]{blommecarocci2025DR}.

\begin{prop}
    \begin{itemize}
        \item $\LogPic^0( \C/\M(a)(\frac{1}{|\bfw|}))$  is endowed with a (log) Abel-Jacobi section inducing a  self duality isomorphism as a logarithmic abelian variety in the sense of \cite{kajiwara2008logabelian2,molchowiselog};
        \item $\LogPic^0( \C/\M(a)(\frac{1}{|\bfw|}))$ is initial for morphism from log smooth curves to (log) abelian varieties; in particular we have a homomorphism of log abelian varieties (dual to the pull-back map)
        \[F_*\colon \LogPic^0( \C/\M(a)(\frac{1}{|\bfw|}))\to\Alb(E)\times \M(a)(\frac{1}{|\bfw|}),\]
        commuting with the log Abel-Jacobi and the Albanese map respectively;
        \item There is a non-degenerate pairing extending the usual Weil pairing for abelian varieties:
        \[W_{\C}\colon\LogPic^0( \C/\M(a)(\frac{1}{|\bfw|}))[|\bfw|]\times \LogPic^0( \C/\M(a)(\frac{1}{|\bfw|}))[|\bfw|]\to\mu_{|\bfw|}\]
        \item The pairing is compatible with $F_*$ in the sense that
        \[W_{\C}(\T,F^*\L)=W_E(F_*\T,\L)\]
    \end{itemize}
\end{prop}

\medskip

\subsubsection{Correlators via Weil-Pairing}

The \emph{correlator function}
\[\theta_{(\cdot,\cdot)}\colon\M(a)(\frac{1}{|\bfw|})\to \Alb(E)[|\bfw|]\cong\Hom(\Pic^0(E)[|\bfw|],\mu_{|\bfw|}) \]
is defined by
\[(f,\T)\mapsto \left( \theta_{(f,\T)}\colon \L\mapsto W_E(f_*\T,\L)=W_{\C}(\T,f^*\L)\right).\]
As the pairing is non-degenerate, the data is equivalent to the data of $f_*(\mathcal T)\in\Alb(E)[|\bfw|]$. 

The correlator function is of course locally constant, so $\M(a)(\frac{1}{|\bfw|})$
is a disjoint union of open and closed sub-stacks indexed by
 $\theta\in\Alb(E)[|\bfw|]$ which we denote by $\M^\theta(a)(\frac{1}{|\bfw|})$.
 It is proved in \cite[Section~6]{blommecarocci2025DR} that this definition of correlator coincides with the one given in \cite{blomme2024correlated} on the locus of rubber maps.

\subsubsection{Correlators on a given refined stratum}

Let us take $[f\colon C\to E]\in \M_{\Gamma,K,\varphi}(\bfa)$ and $\T$ a torsion (log-)line bundle on $C$. We want to understand the relation among:
the correlator of $(f,\T),$ the tropicalization of $\T$ (i.e. the associated torsion tropical divisor on $\Gamma$) and the fixed combinatorial data $(K,\varphi)$.
To do so, let us look at the following commutative diagram. Here, all subgroups are to be understood as $|\bfw|$-torsion, and we dropped the $[|\bfw|]$ to keep a ``light'' diagram:
\bcd
K \ar[rd,"f^*"]\ar[d] &  & & & \\
  \Pic^0(E) \ar[rd,"f^*"] & H^1(\Gamma,\CC^*) \ar[d,hook] & & & \\
  0 \ar[r] & \Pic^{[0]}(C) \ar[r]\ar[d,two heads] & \LogPic^0(C) \ar[r,"\rm{trop}"]\ar[d,"f_*"] & \operatorname{TroPic}^0(C) \ar[d] \ar[r] & 0. \\
  & \bigoplus\Pic^0(C_v) & \Alb(E) \ar[r,"q"] & \Alb(E)/K^\perp & \\
\ecd
\begin{itemize}
    \item The middle short exact sequence is obtained taking  $|\bfw|$-torsion from the short exact sequence relating the Jacobian, logarithmic Jacobian and tropical Jacobian (\cite[Section~4]{molchowiselog},\cite{kajiwara2008logabelian2}). The right morphism is called \emph{tropicalization}.
    \item The vertical short exact sequence in the second column is the usual short exact sequence for the group of multidegree $0$ line bundles on a nodal curve; the second map is the restriction to the components of the normalization.
    \item The diagonal maps on the left are simply the pull-back from $E$, with $K$ the subgroup of  line bundles restricting to the trivial one componentwise.
    \item The middle vertical arrow is the push-forward map $f_*$ mapping a log-line bundle to its correlator as recalled above, and the last map on the right is explained in the next Lemma.
\end{itemize}

 Looking at the diagram above, we can understand the image of the correlator map restricted to torsion log-line bundles with fixed tropicalization, i.e. on the $\Pic^{[0]}(C)$-torsor $\mathrm{trop}^{-1}(\ttD)$ for $\ttD$ the linear equivalence class of a torsion tropical divisor on $\Gamma$.
\begin{lem}\label{lem-description-relation-correlator-tropicalization}
We have the following facts:
\begin{enumerate}
    \item The image of $\Pic^{[0]}(C)$ by $f_*$ is $K^\perp$, where $K^\perp$ is the orthogonal with respect to the Weil pairing.
    \item There is a well-defined map $\TroPic(C)\to\Alb(E)/K^\perp$ obtained by choosing a lift $\U\in \LogPic(C)$ of $\ttD\in\TroPic(C)$ and applying $f_*$. The image is well-defined up to $K^\perp$.
    \item The image of $f_*$ restricted to the torsor of lifts of $\ttD$ is a $K^\perp$-torsor
    \item The map in (2) is actually to  $\varphi^\intercal\colon H_1(\Gamma,\ZZ_{|\bfw|})\to \widehat{K}\simeq \Alb(E)/K^\perp$ 
\end{enumerate}
\end{lem}

\begin{proof}
\begin{enumerate}
    \item To compute $f_*(\Pic^{[0]}(C))$, we compute instead its Weil orthogonal. Since
    $$W_E(\L,f_*\U) = W_C(f^*\L,\U),$$
     looking at $\U\in\Pic^{[0]}(C)$, we have that $\L\in (f_*\Pic^{[0]})^\perp$ if and only if $f^*\L\in\Pic^{[0]}(C)^\perp = H^1(\Gamma,\CC^*)$ (See \cite[Proposition~2.13]{blommecarocci2025DR}). The orthogonal thus consists of those line bundles that restrict to the trivial line bundle on each component $C_v,$ which is the definition of $K$.

    \item Since $\rm{trop}$ is surjective, we can always find a lift of  a given tropical line bundle, only well-defined up to translation by  $\Pic^{[0]}(C)$. By (1), the image by $f_*$ is  thus well-defined up to translation by an element of $K^\perp$.
\item Follows immediately from (2).

    \item Let $\T$ be a torsion log-line bundle with tropicalization $\mathtt{T}$. The correlator is the data of the homomorphism $W_C(\T,f^*(-)) $.
    
Via adjunction, looking at the correlator
modulo $K^\perp$ is the same as looking at $W_C(\T,f^*(\cdot)) $ restricted to $\L\in K$.

Now, for  $\L\in K$, $W(\T,f^*\L)=\widetilde{W}(\mathtt{T},f^*\L)$, where $\widetilde{W}$ is the induced pairing between $H^1(\Gamma,\ZZ_{|\bfw|})$ and $\TroPic(C)[|\bfw|]$ (See again \cite[Proposition~2.13]{blommecarocci2025DR}).
This is precisely $\varphi(\L)(\mathtt{T})$, where $\mathtt{T}$ is interpreted as an element in $H_1(\Gamma,\ZZ_{|\bfw|})$.
\end{enumerate}
\end{proof}
\begin{rem}
 An element in $H_1(\Gamma,\ZZ_{|\bfw|})$ does indeed correspond to a torsion tropical line bundle given that
  $\TroPic^0(\Gamma)=\Hom(H_1(\Gamma,\ZZ),\bar{M}^{\mathrm{gp}})/H_1(\Gamma,\ZZ).$ Here $C$ is a log curve with respect to the minimal log structure on the moduli space of spin stable maps, and $\bar{M}$ is just the base monoid.
\end{rem}

Given the identification of  $\TroPic^0(C)\to \Alb(E)/K^\perp$ with $\varphi^\intercal$ provided in the previous Lemma, we have that the image of the correlator map restricted to the log-line bundles with fixed tropicalization  $\ttD$ is the
 $K^\perp$-torsor of lifts of $\varphi^\intercal(\ttD)$. We set $T_{\varphi^\intercal(\ttD)}$ to be the group algebra average of elements in this class:
$$T_{\varphi^\intercal(\ttD)} = \frac{1}{|K^\perp|}\sum_{\theta\in q^{-1}(\varphi^\intercal(\ttD))}(\theta).$$
This is a element in the group algebra of $\Alb(E)[|\bfw|]$ of total degree $1$. 

\medskip

Using \cite[Lemma 5.14]{blommecarocci2025DR}, we have that for a given $[f\colon C\to E]\in\M_{\Gamma,K,\varphi,\ttD}^{\theta}(\bfa)(\frac{1}{|\bfw|})$,
the cardinality of  $\mathrm{trop}^{-1}(\ttD)\cap f_*^{-1}(\theta)$
 is $|\bfw|^{2g-2-b_1(\Gamma)}|K|$: the $|\bfw|^{2g-b_1(\Gamma)}$ split evenly among the $|K^\perp|=\frac{|\bfw|^2}{|K|}$ correlators. In particular this is saying that the restriction of the (finite) forgetful morphism $\M_{\Gamma,K,\varphi}(\bfa)(\frac{1}{|\bfw|})\to \M_{\Gamma,K,\varphi}(\bfa)$ to the component $\M_{\Gamma,K,\varphi,\ttD}^{\theta}(\bfa)(\frac{1}{|\bfw|})$ has constant degree equal to $|\bfw|^{2g-2-b_1(\Gamma)}|K|$.


\begin{theo}\cite[Theorem 6.7]{blommecarocci2025DR}
\label{theo:correlated-DR-formula}
The total correlated DR-formula at level $|\bfw|=\gcd(w_i)$ for genus $g$ degree $a$ stable maps to an elliptic curve $E$ has the following form:
$$\mathbf{DR}(a,\bfw) = |\bfw|^{2g}  \sum_{\Gamma,\ttD} \frac{1}{|\bfw|^{b_1(\Gamma)}}\sum_{K,\varphi} \vir{\M_{\Gamma,K,\varphi}(\bfa)}\cap\Dec_{\Gamma,\ttD}\left(\frac{\bfw}{|\bfw|}\right) \cdot T_{\varphi^\intercal(\ttD)}  .$$
The first sum is over $E$-valued stable graphs of genus $g$ and degree $a$ together with a choice of equivalence class $\ttD$ of $|\bfw|$-torsion tropical divisor; the second sum is over subgroups $K\subset\Pic^0(E)[|\bfw|]$ and monodromies $\varphi\colon K\to H^1(\Gamma,\ZZ_{|\bfw|})$. 
\end{theo}

We refer to \cite{holmes2023root} (see also\cite[Section~6]{blommecarocci2025DR}) for the explicit expression of the cohomology class $\Dec_{\Gamma,\ttD}\left(\frac{\bfw}{|\bfw|}\right)$, not needed here. The only important thing for us to know about the decoration $\Dec_{\Gamma,\ttD}\left(\frac{\bfw}{|\bfw|}\right)$, is that it depends only \emph{polynomially} in $\frac{\bfw}{|\bfw|}$.

\begin{proof}[Proof of the shape]
    The formula from \cite{blommecarocci2025DR} appears as a sum over the strata $\M_{\Gamma,K,\varphi}(\bfa)$ in the language of piecewise polynomial functions, each stratum coming with a decoration
    $$|\bfw|^{2g-2q(E)}\frac{|K|}{|\bfw|^{b_1(\Gamma)}}\sum_{\ttD}\Dec_{\Gamma,\ttD}(\frac{\bfw}{|\bfw|}),$$
    where the sum is over $\ttD$ which are orthogonal to the image of $\varphi$. To get the advertised expression, the formula is modified as follows:
    \begin{itemize}[leftmargin=0.6cm]
        \item The irregularity $q(E)$ has value $2$.
        \item Since that $|K||K^\perp|=|\bfw|^2$, we can replace $|\bfw|^{-2}|K|=\frac{1}{|K^\perp|}$,
        \item In addition, Lemma \ref{lem-description-relation-correlator-tropicalization} describes the possible correlators on a given stratum. Thus, instead of summing over the $\ttD$ contributing to the $(0)$-correlator, we assign to each $\ttD$ the group algebra element equal to the sum of its correlators: $|K^\perp|T_{\varphi^\intercal}(\ttD)$. The $|K^\perp|$ cancels with the previous $\frac{1}{|K^\perp|}$.
    \end{itemize}
\end{proof}

\begin{rem}\label{rem:fromDRtologDR}
    
The correlated refinement of $\mathrm{logDR}(a,\bfw)$ has essentially the same shape, with the difference that since the class representing the log lift is obtained applying the universal DR-formula \cite{bae2023pixton} to a twist of the universal root of $\mathcal C\to\M(a)(\frac{1}{|\bfw|})$ on some tropical model $\mathcal C^\diamond\to\M^\diamond(a)(\frac{1}{|\bfw|})$. In the first sum, $\Gamma$ is replaced by the finer data structure $\Gamma^\diamond$, $\ttD$ is replaced by a tropically equivalent divisor $\ttD^{\sigma}$ and $\M_{\Gamma,K,\phi}(\bfa)$ by some strata  $\M_{\Gamma^\diamond,K,\phi}(\bfa)$ of the a tropical model $\M(\bfa)^\diamond.$

We refer to \cite[Section~5]{blommecarocci2025DR}, \cite[Section~3.6]{BCUK2026corrDDR}
 for an explanation on how to suitably choose the torsion divisor $\ttD^{\sigma}$ depending on $\bfw$. 
\end{rem}


\subsection{Multiple cover formula for the correlated DR}

Let $\gamma\in H^*(E,\QQ)^{\otimes n}$ be of the following type:
$$\gamma\in H^0(E,\QQ)^{\otimes t}\otimes H^1(E,\QQ)^{\otimes 2r}\otimes H^2(E,\QQ)^{\otimes s},$$
with $t+2r+s=n$, which we may always assume up to relabeling the marked points. The degree of $\gamma$ is $\deg\gamma=r+s$. We now consider the following group algebra valued cycle in $H^*(\bar{\M}_{g,n},\QQ)$:
$$\ft_*\left( \mathbf{DR}(a,\bfw) \cap \ev^*\gamma \right).$$
To get concrete numbers, one would have to intersect with an homology class $\alpha\in H^*(\bar{\M}_{g,n},\QQ)$ of complementary dimension, i.e.
$$\deg\alpha + \deg\gamma +g= 2g-2+n$$

\begin{theo}\label{theo:MCF-for-correlated-DR}
    The correlated cycle $\ft_*(\mathbf{DR}(a,\bfw)\cap\ev^*\gamma)$ where $\gamma\in H^*(E,\QQ)^{\otimes n}$ satisfies the $(2g-1+\deg\gamma)$-multiple cover formula.
\end{theo}

\begin{rem}
    In terms of a complementary insertion $\alpha$, the exponent is
    $$2g-1+\deg\gamma = 3g-3+n-\deg\alpha.$$
\end{rem}

To prove the MCF formula we show that the formula is already true on each stratum $(\Gamma,\mathtt{D})$, where it further reduces to the MCF for refined group counting functions from Section \ref{sec:refined-groups-counting}.

\begin{proof}[Proof of Theorem \ref{theo:MCF-for-correlated-DR}]
The idea is to express the function as a sum indexed by strata $(\Gamma,\ttD)$, and prove the MCF for the summand. Along the way, we emphasize the role of the smooth stratum where $\Gamma$ is the trivial graph, for which the situation is slightly easier to understand.

\medskip

\textbf{$\bullet$ Step 1: expression of the contribution.} The correlated DR-cycle for a choice of $(a,\bfw)$ and insertion $\gamma$ appears as a sum indexed by $\Gamma$ a $E$-valued stable graph of degree $a$ and $\ttD$ be a tropical $|\bfw|$-torsion divisor. Our first step is to rewrite the expression as a sum.

Since the degree $a$ also gets scaled in the process, we separate the degree information from the graph: we rather sum over $E$-valued stable graphs devoided of the degree information $\bfa$, and later sum over the possible degree assignments $\bfa$ such that $\sum a_v=a$.

For the tropical torsion divisor, we sum over all possible torsion divisors, of any order, with a contribution understood to be $0$ if $\ttD$ is not $|\bfw|$-torsion.

With the above modification, we get
$$\ft_*(\mathbf{DR}(a,\bfw)\cap\ev^*\gamma) = \sum_{\Gamma,\ttD} \sum_{\bfa\vdash a}\C_{\Gamma,\ttD}(\bfa,\bfw),$$
where the first sum is over graphs and tropical torsion divisors, the second sum is over $\bfa=(a_v)$ such that $\sum a_v=a$, abbreviated $\bfa\vdash a$. The expression of $\C_{\Gamma,\ttD}(\bfa,\bfw)$ is provided by Theorem \ref{theo:correlated-DR-formula}:
$$\C_{\Gamma,\ttD}(\bfa,\bfw) = |\bfw|^{2g-b_1(\Gamma)}\sum_{K,\varphi} \ft_*\left( \vir{\M_{\Gamma,K,\varphi}(\bfa)}\cap ( \Dec_{\Gamma,\ttD}(\frac{\bfw}{|\bfw|})\cup\ev^*\gamma) \right)T_{\varphi^\intercal(\ttD)}.$$


The quantity $\frac{\bfw}{|\bfw|}$ does not change when $\bfw$ gets scaled. As already explained in the previous sections, to prove the MCF we may restrict to multiples of a primitive element $(\widetilde{a},\widetilde{\bfw})\in\RRR.$ Therefore, we drop the $\frac{\bfw}{|\bfw|}$ from the Dec to lighten the notation.

As a toy-case, we carry over the computation for the smooth stratum. The contribution for the smooth stratum is obtained taking the trivial graph $\Gamma^\mathrm{sm}$ with a unique vertex, no loops, and $\ttD$ trivial as well. The contribution is then
$$\C^\mathrm{sm}(a,\bfw) = |\bfw|^{2g}\sum_K \ft_*\left( \vir{\M_K(a)}\cap( \Dec_{\Gamma^\mathrm{sm},0}\cup\ev^*\gamma) \right)T_{K^\perp}.$$

The cycle-valued function $\C_{\Gamma,\ttD}$ is defined on the $\NN$-module $\RRR_\Gamma=\ZZ_{\geqslant 0}^{\sfV(\Gamma)}\times\{\bfw:\sum w_i=0\}$ endowed with norm $|(\bfa,\bfw)|=|\bfw|$. Hence,
$$(a,\bfw)\mapsto\sum_{\bfa\vdash a}\C_{\Gamma,\ttD}(\bfa,\bfw),$$
is actually the push-forward of $\C_{\Gamma,\ttD}$ by the morphism of $\NN$-module $((a_v),\bfw)\in\RRR_\Gamma\mapsto(\sum a_v,\bfw)\in\RRR$. Therefore, to prove the MCF, it suffices to prove that each of the $\C_{\Gamma,\ttD}$ satisfies the MCF. We then conclude by Proposition \ref{prop:push-forward} and linearity. The rest of the proof consists in proving the MCF for each $\C_{\Gamma,\ttD}$.

\medskip

\textbf{$\bullet$ Step 2: Reducing to a product over the vertices.} This step is not needed when $\Gamma$ is the trivial graph since there is a unique vertex, and the class is already $\vir{\M_K(a)}$. In case $\Gamma$ is not trivial, we replace $\vir{\M_{\Gamma,K,\varphi}(\bfa)}$ by its expression as a product over vertices:
$$\vir{\M_{\Gamma,K,\varphi}(\bfa)}\cap(\Dec_{\Gamma,\ttD}\cup\ev^*\gamma) =  \vir{\P_{\Gamma,K}(\bfa)}\cap (\ev^*\Delta_{K,\varphi}\cup\Dec_{\Gamma,\ttD}\cup\ev^*\gamma),$$
where $\Delta_{K,\varphi}$ is the twisted diagonal corresponding to $\varphi$. 

We proved in Proposition~\ref{prop:relation-to-covering-nodal-case} that these numbers do not depend on $\varphi$. Recall that $\P_{\Gamma,K}(\bfa)$ is the component of the product  $\prod\M_v(a_v)$ where  $K=\cap K_v$.

\medskip

\textbf{$\bullet$ Step 3: reducing to invariants for coverings of $E$ for $\C^\mathrm{sm}(a,\bfw)$.} We now replace the refined virtual classes for $E$ by virtual classes for coverings of $E$. Therefore, we now exceptionally write again the target spaces in the moduli spaces to avoid any confusion. We use the relation between $\M_K(E,a)$ and $\M(E_H,a/|H|)$ for $H\supset K$ provided by Proposition~\ref{prop:relation-to-covering-smooth-case}:
$$\vir{\M_K(E,a)} \cap (\Dec_{\Gamma^\mathrm{sm},0}\cup\ev^*\gamma) = \sum_{H\supset K}\frac{\mu(H/K)}{|H|} \pi_{H\ast} \left( \vir{\M(E_H,a/|H|)}\cap (\Dec_{\Gamma^\mathrm{sm},0}\cup\ev^*\pi_H^*\gamma)\right).$$
Using Corollary \ref{coro:pull-back-covering}, and recalling that  $E$ and $E_H$ are both elliptic curves (allowing for the already explained abuse of notation for the insertions $\gamma$) we can rewrite the push-forward to $\bar{\M}_{g,n}$ as follows:
$$\ft_*\left( \vir{\M_K(a)}\cap (\Dec_{\Gamma^\mathrm{sm},0}\cup\ev^*\gamma )\right) = \sum_{H\supset K}\mu(H/K) |H|^{r+s-1} \ft_*\left( \vir{\M(a/|H|)} \cap (\Dec_{\Gamma^\mathrm{sm},0}\cup\ev^*\gamma)\right).$$
Now, summing over $K$ and regrouping according to the value of the cardinality of $|H|$, we get
\begin{align*}
    \C^\mathrm{sm}(a,\bfw) = & |\bfw|^{2g}\sum_K\sum_{H\supset K}\mu(H/K)|H|^{r+s-1}\ft_*\left(\vir{\M(a/|H|)}\cap(\Dec_{\Gamma^\mathrm{sm},0}\cup\ev^*\gamma)\right) T_{K^\perp} \\
    = & |\bfw|^{2g} \sum_{\delta||\bfw|^2}\delta^{r+s-1}\ft_*\left(\vir{\M(a/\delta)}\cap(\Dec_{\Gamma^\mathrm{sm},0}\cup\ev^*\gamma\right)\sum_{\substack{H\supset K \\ |H|=\delta}}\mu(H/K)T_{K^\perp}.
\end{align*}

As $E_H$ is an elliptic curve just as $E$, they have the same GW-cycles, thus depending only on the degree since the insertions are always identified and the remaining discrete data agree. We denote by
$$\GW^\mathrm{sm}(a)=\ft_*\left(\vir{\M(a)}\cap (\Dec_{\Gamma^\mathrm{sm},0}\cup\ev^*\gamma)\right).$$
Passing to the orthogonals with respect the Weil pairing, the inner sum is a sum over the $K^\perp\supset H^\perp$ with $H^\perp$ of index $\delta$. Furthermore, $\mu(K^\perp/H^\perp)=\mu(H/K)$ since $K^\perp/H^\perp$ is actually isomorphic to the dual $\widehat{H/K}$, itself isomorphic (non canonically) to $H/K$. We thus recognize the group counting function $\bfF(\delta,|\bfw|)$ from Section \ref{sec:refined-groups-counting}. In the end, we get
$$\C^\mathrm{sm}(a,\bfw) = |\bfw|^{2g} \sum_{\delta||\bfw|^2}\delta^{r+s-1}\GW^\mathrm{sm}(a/\delta)\bfF(\delta,|\bfw|).$$

\textbf{$\bullet$ Step 3bis: reducing to invariants for coverings of $E$ for $\C_{\Gamma,\ttD}(a,\bfw)$.} We now take a general $\Gamma$. Using that the virtual class $\vir{\M_{\Gamma,K,\varphi}(\bfa)}$ is the Gysin pull-back of the virtual class of $\prod_v\M_{v,K_v}(a_v)$ along the twisted diagognal $\Delta_{K,\varphi}$, after push-forward to the product 
and using the fact that the twisted diagonals $\Delta_{K,\varphi}$ are all cobordant we get:
\begin{align*}
    \C_{\Gamma,\ttD}(\bfa,\bfw) = & |\bfw|^{2g-b_1(\Gamma)}\sum_{K,\varphi} \ft_*\left(\vir{\M_{\Gamma,K,\varphi}(E,\bfa)}\cap(\Dec_{\Gamma,\ttD}\cup\ev^*\gamma)\right) T_{\varphi^\intercal(\ttD)} \\
    = & |\bfw|^{2g-b_1(\Gamma)}\sum_{K,\varphi} \ft_*\left(\vir{\P_{\Gamma,K}(E,\bfa)}\cap (\ev^*\Delta_{K,\varphi}\cup\Dec_{\Gamma,\ttD}\cup\ev^*\gamma)\right) T_{\varphi^\intercal(\ttD)} \\
    = & |\bfw|^{2g-b_1(\Gamma)}\sum_{K} \ft_*\left(\vir{\P_{\Gamma,K}(E,\bfa)}\cap(\ev^*\Delta_{K,0}\cup\Dec_{\Gamma,\ttD}\cup\ev^*\gamma)\right) \sum_{\varphi}T_{\varphi^\intercal(\ttD)},
\end{align*}
where in the last row we used that the class $\Delta_{K,\varphi}$ does not depend on $\varphi$, and we can thus factor the integral out of the sum over $\varphi$. We can then use Proposition \ref{prop:relation-to-covering-nodal-case} to express the contribution in terms of invariants for the coverings of $E$. Then, we group them according to the cardinality of $|H|$ as in the smooth case:
\begin{align*}
     & \C_{\Gamma,\ttD}(\bfa,\bfw) \\
    = & |\bfw|^{2g-b_1(\Gamma)}\sum_K\sum_{H\supset K}\mu(H/K)\frac{|H|^{b_1(\Gamma)+r+s-1}}{|K|^{b_1(\Gamma)}}\ft_*\left(\vir{\P(\bfa/|H|)}\cap(\ev^*\widetilde{\Delta}_{H,0}\cup\Dec_{\Gamma,\ttD}\cup\ev^*\gamma)\right)\sum_\varphi T_{\varphi^\intercal(\ttD)} \\
    = & |\bfw|^{2g-b_1(\Gamma)}\sum_{\delta||\bfw|^2}\delta^{b_1(\Gamma)+r+s-1}\GW_{\Gamma,\ttD}(\bfa/\delta) \cdot \frac{1}{|K|^{b_1(\Gamma)}}\sum_{\substack{H\supset K,\varphi \\ |H|=\delta}}\mu(H/K)T_{\varphi^\intercal(\ttD)},
\end{align*}
where we set
$$\GW_{\Gamma,\ttD}(\bfa) = \ft_*\left(\vir{\P_\Gamma(\bfa)}\cap(\ev^*\Delta_0\cup\operatorname{Dec}_{\Gamma,\ttD}\cup\ev^*\gamma)\right),$$
with $\Delta_0$ the usual diagonal.

We now reformulate the group algebra count appearing in the formula.
The sum over $\varphi$ can equivalently be seen as a sum over $\varphi^\intercal\colon H_1(\Gamma,\ZZ_{|\bfw|})\to\widehat{K}\simeq\Alb(E)/K^\perp$. If $\omega$ is the order of $\ttD$, we can write $\ttD=\frac{|\bfw|}{\omega}\mathtt{T}$ for a primitive tropical torsion divisor $\mathtt{T}$ and complete $\mathtt{T}$ to a basis of $H_1(\Gamma,\ZZ_{|\bfw|})$.

Then, though we are summing over choices of morphism $\varphi^\intercal\colon\ZZ_{|\bfw|}^{b_1(\Gamma)}\to\Alb(E)/K^\perp$, only the choice of $\varphi^\intercal(\mathtt{T})$  matters.
Since there are  $|K|$ morphisms $\ZZ_{|\bfw|}\to\Alb(E)/K^\perp$, there are $|K|^{b_1(\Gamma)-1}$ possible choices for the other basis elements. Thus, we get:
$$\frac{1}{|K|^{b_1(\Gamma)}}\sum_{\substack{H\supset K,\varphi \\ |H|=\delta}} \mu(H/K) T_{\varphi^\intercal(\ttD)} = \frac{1}{|K|}\sum_{\substack{H\supset K \\ |H|=\delta}}\mu(H/K) \sum_{\psi\colon\ZZ_{|\bfw|}\to\Alb(E)/K^\perp}T_{q^{-1}(\psi(|\bfw|/\omega))},$$
where $q$ denotes the projection $\Alb(E)[|\bfw|]\to\Alb(E)[|\bfw|]/K^\perp$.  As in the smooth case, the first sum can be viewed as a sum over $K^\perp\supset H^\perp$ with $H^\perp$ of index $\delta.$ 
We thus recognize the second group counting function $\bfG_\omega(\delta,|\bfw|)$ from Section \ref{sec:group-counting-Gomega} with $\omega$ the order of the torsion tropical divisor $\ttD$. 

In the end, since $\deg\gamma=r+s$, we get
$$\C_{\Gamma,\ttD}(\bfa,\bfw) = |\bfw|^{2g-b_1(\Gamma)}\sum_{\delta||\bfw|^2}\delta^{b_1(\Gamma)+\deg\gamma-1}\GW_{\Gamma,\ttD}(\bfa/\delta)\cdot\bfG_\omega(\delta,|\bfw|).$$

\textbf{$\bullet$ Step 4: Reducing to MCF for counts of subgroups.} From the expression found in the previous step, we have
$$\Prim_{|\bfw|}\C_{\Gamma,\ttD}(\bfa,\bfw) = |\bfw|^{2g-b_1(\Gamma)}\sum_{\delta||\bfw|^2}\delta^{b_1(\Gamma)+\deg\gamma-1}\GW_{\Gamma,\ttD}(\bfa/\delta)\Prim_{|\bfw|}\bfG_\omega(\delta,|\bfw|).$$
We can now compute the right-side of the expected MCF:
\begin{align*}
     & \sum_{k|\bfa,\bfw} k^{2g+\deg\gamma-1}\Prim_{|\bfw|/k}\C_{\Gamma,\ttD}(a/k,\bfw/k) \\
    = & \sum_{k|\bfa,\bfw} k^{2g+\deg\gamma-1}\left(\frac{|\bfw|}{k}\right)^{2g-b_1(\Gamma)}\sum_{\delta|(\frac{|\bfw|}{k})^2}\delta^{b_1(\Gamma)+\deg\gamma-1}\GW_{\Gamma,\ttD}(\bfa/k\delta)\Prim_{|\bfw/k|}\bfG_\omega(\delta,\frac{|\bfw|}{k}) \\
    = & |\bfw|^{2g-b_1(\Gamma)}\sum_{k|\bfa,\bfw}\sum_{\delta|(\frac{|\bfw|}{k})^2} (k\delta)^{b_1(\Gamma)+\deg\gamma-1}\GW_{\Gamma,\ttD}(\bfa/k\delta)\Prim_{|\bfw/k|}\bfG_\omega(\delta,\frac{|\bfw|}{k}) \\
    = & |\bfw|^{2g-b_1(\Gamma)}\sum_{d||\bfw|^2}d^{b_1(\Gamma)+\deg\gamma-1}\GW_{\Gamma,\ttD}(\bfa/d)\sum_{k|d,\frac{|\bfw|^2}{d},\bfw}\Prim_{|\bfw/k|}\bfG_\omega(\frac{d}{k},\frac{|\bfw|}{k}),
\end{align*}
where in the last row we set $d=k\delta$. It is subject to the conditions $d||\bfw|^2$. Then $k$ has to satisfy the following: $k|d,\bfa,\bfw,\frac{|\bfw|^2}{d}$. Notice that we need to have $d|\bfa$ as well for $\GW_{\Gamma,\ttD}(\bfa/d)$ to be non-zero so that the condition $k|\bfa$ is redundant. The condition $k|d,\bfw,\frac{|\bfw|^2}{d}$ may be seen as $k|(d,\bfw)$ and still having $\frac{d}{k}|(\frac{|\bfw|}{k})^2$, which is equivalent to $k|\frac{|\bfw|^2}{d}$.

As the function $\bfG_\omega(\delta,n)$ satisfies the $0$-MCF, the inner sum is actually $\bfG_\omega(d,|\bfw|)$, so that we recover the expression of $\C_{\Gamma,\ttD}(\bfa,\bfw)$, finishing the proof.
\end{proof}


\section{Multiple cover formula for $E\times\PP^1$}
\label{SEC-MCFEP1}

We now get into the proof of Theorem \ref{theo-MCF-EP1-general}, i.e. the MCF for (exotic) relative invariants of $E\times\PP^1$. The idea is to adapt the results of \cite{ranganathan2024logarithmic}, which reduce computation of (correlated) GW-cycles for $E\times\PP^1$ to intersection with the (correlated) DR-cycle, for which Theorem \ref{theo:MCF-for-correlated-DR} provides the MCF. 


\subsection{Correlated exotic Gromov-Witten invariants}
\label{sec:MCF-EP1-setting}

In Section~\ref{sec:correlation} we recalled the definition of correlated virtual class and correlated GW invariants for the geometric setting of interest. 
To keep a light notation, we continue denoting our
 moduli space of log-stable maps to $E\times\PP^1$ by 
\[\M(a,\bfw)=\M_{g,m}(E\times\PP^1|E_0+E_{\infty},a,\bfw), \]
recalling only $a,\bfw$ that will play a role in the MCF. 

\medskip

Logarithmic GW-invariants of  $E\times \PP^1$ are obtained capping the virtual fundamental class $\vir{\M(a,\bfw)}$ of degree $n+m+g-1$ with cohomology classes pulled-back along the forgetful and evaluation morphisms:
\begin{align*}
 \ft& \colon\M(a,\bfw)\to\bar{\M}_{g,n+m},  \\ 
 \ev &\colon\M(a,\bfw)\to (E\times\PP^1)^m\times E^n :=\Ev(E\times\PP^1).
\end{align*}

Similarly, we can consider the moduli space of rubber stable maps
\[R\M(a,\bfw)=R\M_{g,m}(E\times\PP^1|E_0+E_{\infty},a,\bfw),\]
associated to the same discrete data, whose virtual class of degree  $n+m+g-2$ pushes forward to the DR-cycle along $\DR(a,\bfw)$ along 
$\epsilon\colon R\M(a,\bfw)\to\M(a).$
We also have evaluation and forgetful maps

\begin{align*}
 \ft& \colon R\M(a,\bfw)\to\M(a)\to\bar{\M}_{g,n+m},  \\ 
\bar{\ev} &\colon R\M(a,\bfw)\to\M(a)\to  E^{m+n},
\end{align*}
though the evaluation map is going to be refined to get value in a bigger space. We define \emph{rubber invariants} capping the DR-cycle with classes pulled-back along $\ft$ and $\bar{\ev}.$







\subsubsection{Correlation}

We have the correlator maps
$$\kappa\colon\M(a,\bfw)\to E[|\bfw|],\quad \kappa\colon R\M(a,\bfw)\to E[|\bfw|],$$
both denoted by $\kappa$ and with values in the $|\bfw|$-torsion of $E$. In particular, we get, as already recalled before a correlated refinement and correlated virtual classes \cite{blomme2024correlated}:

\[\fvir{\M(a,\bfw)}{|\bfw|}=\sum_{|\bfw|\cdot\theta=0}\vir{\M(a,\bfw)^\theta}(\theta),\qquad  \fvir{R\M(a,\bfw)}{|\bfw|}=\sum_{|\bfw|\cdot\theta=0}\vir{R\M(a,\bfw)^\theta}(\theta)\]
The map $\mu\colon \M(a,\bfw)\to  R\M(a,\bfw)$ is compatible with the correlator maps by definition. Pushing forward along $\epsilon$, we obtain the correlated DR- cycle $\mathbf{DR}(a,\bfw) $ considered in Section~\ref{sec-correlated-DR}\footnote{The definition of correlation recalled in Section~\ref{sec-correlated-DR} seems at first sight different from the one recalled here. We refer to \cite{blommecarocci2025DR} for a proof of the fact that the two definition in fact coincide.}.

The \emph{correlated} GW-invariants (resp. correlated rubber invariants) are defined capping the correlated class against cohomology class pulled-back along evaluation/forgetful maps. We adopt similar notations: for $\alpha\in H^*(\bar{M}_{g,n},\QQ)$ and $\gamma$ a cohomology class in the suitable evaluation space,

\begin{align*}
   \gen{\gen{\alpha;\gamma}}_{a,\bfw}  & =  \fvir{\M(a,\bfw)}{|\bfw|}\cap(\ft^*\alpha\cup\ev^*\gamma) \in \G_{|\bfw|},\\
    \gen{\gen{\gamma}}_{a,\bfw} & = \ft_*(\fvir{\M(a,\bfw)}{|\bfw|}\cap\ev^*\gamma) \in H^*(\bar{\M}_{g,n+m},\QQ)\otimes\G_{|\bfw|}\\
\gen{\gen{\alpha;\gamma}}_{a,\bfw} & \mathbf{DR}(a,\bfw) \cap (\ft^*\alpha\cup\bar{\ev}^*\bar{\gamma})\in  \G_{|\bfw|},\\
\gen{\gen{\gamma}}_{a,\bfw} & \ft_*( \mathbf{DR}(a,\bfw) \cap \bar{\ev}^*\bar{\gamma})\in H^*(\bar{\M}_{g,n+m},\QQ)\otimes\G_{|\bfw|}
  \end{align*}

\subsubsection{Exotic insertions}
By Section~\ref{sec:virtual_birational model}, the correlated virtual class and the correlated DR-cycle can be lifted to classes in $\mathrm{logCH}_*(\M(a,\bfw))\otimes \G_{|\bfw|}$ and $\mathrm{logCH}_*(\M(a))\otimes \G_{|\bfw|}$ respectively. The lift of the correlated DR-cycle has already been discussed in Section~\ref{sec:MCFlogcorrelatedDR}.

Indeed, the map $\kappa$ used to show that the moduli space is disconnected and to define the correlated components lifts (simply by post-composition) to any tropical model $\M^\diamond(a,\bfw)$. The projection $\M^\diamond(a,\bfw)\to\M(a,\bfw)$ identifying the respective classes via pushforward is the disjoint union of the projections $\M^{\theta,\diamond}(a,\bfw)\to\M^\theta(a,\bfw)$.

\medskip

For any $\Ev(E\times\PP^1)^\diamond\to \Ev(E\times\mathbb P^1)$ logarithmic modification, and given
$\M^\diamond(a,\bfw)$ a tropical (and virtually birational ) model admitting an evaluation map

\[\M^\diamond(a,\bfw)\to \Ev(E\times\PP^1)^\diamond\]
we define the correlated exotic invariants:
\begin{align*}
     \gen{\gen{\alpha;\gamma^\diamond}}_{a,\bfw}  & =  \fvir{\M^\diamond(a,\bfw)}{|\bfw|}\cap(\ft^*\alpha\cup\ev^*\gamma^\diamond) \in \G_{|\bfw|},\\
    \gen{\gen{\gamma^\diamond}}_{a,\bfw} & = \ft_*(\fvir{\M^\diamond(a,\bfw)}{|\bfw|}\cap\ev^*\gamma^\diamond) \in \mathrm{logH}^*(\bar{\M}_{g,n+m},\QQ)\otimes\G_{|\bfw|}
\end{align*}

Since the logarithmic structure on $\Ev(E\times\PP^1)=E^{n+m}\times (\PP^1)^m$ is pulled back from 
$(\PP^1)^m$, every logarithmic modification $ \Ev(E\times\PP^1)^\diamond$ is going to be of the form $E^{n+m}\times P_{\Sigma}$ for $P_{\Sigma}\to(\PP^1)^m$ some toric blow-up.



\subsection{From rigid to rubber}
\label{sec:MCF-EP1-rigid-to-rubber}

This section is a straightforward adaptation of a special case of \cite{ranganathan2024logarithmic}. 
The main result is the following:

\begin{theo}
\label{theo:rigid-to-rubber}
    Let $\gamma^\diamond = \bar{\gamma}\otimes\gamma_t\in H^*( E^{n+m})\otimes H^*(P_{\Sigma})$ be an exotic insertion. There exists $\Phi_{\bfw}$ a piecewise polynomial on $\Sigma_{\bar{\M}_{g,n}}$ (corresponding to a cohomology class in some tropical model $\bar{\M}^\diamond$) such that 
    \[\ft_*(\vir{\M(a,\bfw)^{\diamond}}  {}\cap\ev^*\gamma^\diamond) = \ft_*(\mathrm{logDR}(a,\bfw)\cap ( \ev^*\bar{\gamma}\cup \Phi_{\bfw}) \in \mathrm{logCH}(\bar{\M}_{g,n+m})).\]
\end{theo}


\begin{rem}
We deal only with the case of $E\times\PP^1$ to try and minimize the length of this already long paper, but we are certain that there is a more general framework to which \cite{ranganathan2024logarithmic} can be adapted, allowing a comparison between log-GW invariants of toric variety bundles in the sense of \cite{CN}, and higher ramification cycles with target varieties.
\end{rem}

\subsubsection{Minimal tropical models and virtual classes}

Since in what follows $a$ and $\bfw$ do not play any major role, we instead want to emphasize the target space: we denote $\M(a,\bfw)$ by $\M(E\times\PP^1)$, and $\M(a)$ by $\M(E)$.


\medskip

For the comparison, it is convenient to work with certain \emph{minimal} tropical model of the moduli spaces of logarithmic resp. rubber maps.

\medskip

In this section, all schemes are logarithmic schemes, but we will simply write $S$ for $(S,M_S).$
Let $\Gmlog, \Gmtrop$ the logarithmic and tropical multiplicative group respectively, whose values at a log-scheme $S$ are respectively $H^0(S,M_S^\gp), H^0(S,\bar{M}_S^\gp)$. 

\medskip

We consider the following moduli stack over logarithmic schemes (see also \cite[Section~2]{blomme2024correlated},\cite[Section~1]{ranganathan2024logarithmic}):

\begin{itemize}[leftmargin=0.6cm]
    \item The stack over log-schemes $\M(E\times\Gmlog)$; an $S$-point of this stack for  a fs log-scheme is the data of:
    $f\colon C/S\to E$ in $\M(E)$ a usual stable map of the underlying schemes, together with a section of $\pi_*M^\gp_C$.
    
    As explained in \cite[Section~1]{ranganathan2024logarithmic}, it is endowed with an evaluation map:
    $$\ev\colon \M(E\times\Gmlog) \to E^{n+m}\times\Gmlog^m := \Ev(E\times\Gmlog).$$

    We remark that  $\M(E\times\Gmlog)$ is \emph{not} representable by an algebraic stack with log structure, but admits tropical models that are.

    \item The stack over log-schemes $R\M(E\times\Gmlog)$ (see \cite[Definition~1.6.1]{ranganathan2024logarithmic}) parametrizing a usual stable maps $f\colon C/S\to E$ in $\M(E)$ of the underlying schemes, together with a section of $\pi_*M^\gp_C/M^{\gp}_S$. It can also be described as the  fiber product \cite{marcuswiselog}:
   \bcd
   R\M(E\times\Gmlog) \ar[r]\ar[d] & \operatorname{Div}_{g,m,\bfw}\ar[d, "\mathtt{aj}"] \\
   \M(E) \ar[r, "\mathcal O"] & \mathfrak{Pic}.
   \ecd
   where $\mathfrak{Pic}$ is the universal Picard stack over the space of pre-stable curves, $\operatorname{Div}_{g,m,\bfw}$ is Marcuse-Wise algebraic stack with log structure (see \cite[Theorem~4.2.4]{marcuswiselog} parametrizing log-smooth curves together with sections of $\pi_*\bar{M}^\gp_C/\bar{M}^\gp_S.$

   In \cite{ranganathan2024logarithmic} the stack $\operatorname{Div}_{g,m,\bfw}$ is denoted by $\mathfrak{M}_{g,m,\bfw}^{\mathrm{rub}}(\Gmtrop),$ or $R_{g,m,\bfw}(\Gmtrop)$ for the open sub-stack where the curve marked curve is stable.

   \medskip

   Unlike $\M(E\times\Gmlog),$ the log stack $   R\M(E\times\Gmlog)$ is representable by an algebraic stack with log structure. Furthermore, by \cite[Section~4.6]{marcuswiselog}, the morphism $\mathtt{aj}$ is endowed with a perfect obstruction theory, we thus get a class
\[\vir{R\M(E\times\Gmlog)}=\mathtt{aj}^!\vir{\M(E)}.\]
    
    \item The stack $R\M(E\times\Gmtrop)$ parametrizing a usual  stable maps $f\colon C/S\to E$ in $\M(E)$ together with a section of $\pi_*\bar{M}^\gp_C/\bar{M}^\gp_S$, or in other words a PL-function on the tropicalization up to translation.  This can also be described as a fiber product:

    \bcd
R\M(E\times\Gmtrop)\ar[d]\ar[r] & \operatorname{Div}_{g,m,\bfw}\ar[d]\\
\M(E)\ar[r] & \bar{\M}_{g,n+m},
\ecd\textbf{}
where $\bar{\M}_{g,n+m}$ is considered as a stack over log-schemes endowing it with its natural boundary log structure.

By \cite[Corollary~4.2.6]{marcuswiselog} the right vertical arrow is in fact a proper and birational logarithmic modification between log smooth (and in particular irreducible) algebraic stack with log structure. It then follows from Section~\ref{sec:virtual_birational model} that $R\M(E\times\Gmtrop)$ is a virtual birational model of $\M(E).$

Furthermore, again using the obstruction theory of the morphism $\mathtt{aj}$ and the functoriality of virtual pull-back we can also write 
    \[\vir{R\M(E\times\Gmlog)}=aj^!\vir{R\M(E\times\Gmtrop)}.\]
\end{itemize}

By \cite[Propositions~1.5.5,1.5.6]{ranganathanGWexpansion} 
and by \cite[Section~5]{marcuswiselog} respectively, 
we have logarithmic modifications 
 \begin{align*}
 \M(E\times\PP^1)&\longrightarrow \M(E\times\Gmlog),\\
  R\M(E\times\PP^1)&\longrightarrow R\M(E\times\Gmlog). 
 \end{align*}
By \cite[Section~4,5]{marcuswiselog} the obstruction theories for the rubber spaces $  R\M(E\times\PP^1)$ and $ R\M(E\times\Gmlog)$ are compatible; in particular push-forward along the horizontal map identify the virtual classes and thus we also have $\epsilon_*\vir{R\M(E\times\Gmlog)}=\DR(a,\bfw).$

\subsubsection{Evaluation spaces}

We have the following commutative diagram of evaluation maps:
\bcd 
\M(E\times\PP^1)\ar[r]\ar[d,"\ev"] &  \M(E\times\Gmlog) \ar[d,"\ev"]\\
E^{m+n}\times (\PP^1)^m\ar[r] & E^{m+n}\times\Gmlog^m,
\ecd
where the top map is induced by the logarithmic subdivision $\PP^1\to\Gmlog$.

\medskip

Following  \cite[Section~2]{ranganathan2024logarithmic}, there is on the rubber side an evaluation map that we now describe. Its description can also be found in \cite[Section 6.2.3]{BCUK2026corrDDR} in the toric case. We have the diagonal action of $\Gmlog$ on $\Gmlog^m$, and thus on $\Ev(E\times\Gmlog)$ by acting on the last components. Let
$$\Ev^\rub(E\times\Gmlog) = \Ev(E\times\Gmlog)/\Gmlog = E^{m+n}\times \Gmlog^m/\Gmlog,$$
where the quotient is taken in the category of abelian groups valued sheaves on logarithmic schemes.

As explained in \cite[Section~2.4]{ranganathan2024logarithmic}, the evaluation map $\ev\colon\M(E\times\Gmlog)\to\Ev(E\times\Gmlog)$ descends to an evaluation map $\ev^{\rub}\colon R\M(E\times\Gmlog)\to \Ev^{\rub}(E\times\Gmlog)$ defined on the rubber space. The idea is that though evaluation in the $\Gmlog$-coordinates does not make sense for maps defined up to translation, the ratio of such coordinates still makes sense. We obtain the following cartesian diagram:
\begin{equation}\label{eq:rigid-torubbercartesian}
    \begin{tikzcd}
        \M(E\times\Gmlog)\ar[d]\ar[r,"\ev"] & \Ev(E\times\Gmlog)\ar[d]\\
R\M(E\times\Gmlog)\ar[r,"\ev^{\rub}"] & \Ev^{\rub}(E\times\Gmlog)
    \end{tikzcd}.
\end{equation}

\subsubsection{}

The stacks over log-schemes appearing in \eqref{eq:rigid-torubbercartesian} (including the evaluation spaces) are not representable by algebraic stacks/schemes with logarithmic structures. Nevertheless, we can always choose logarithmic subdivisions that are representable. This way, it is proved in  \cite[Proposition~2.4.5]{ranganathan2024logarithmic} that \eqref{eq:rigid-torubbercartesian} can be replaced, up to taking logarithmic modifications, by a cartesian square of algebraic stacks (in fact Deligne-Mumford stacks) with fs log structures:
\begin{equation}\label{eq:rigid-torubbercartesiandagger}
    \begin{tikzcd}
        \M(E\times\Gmlog)^\diamond\ar[d]\ar[r,"\ev"] & \Ev(E\times\Gmlog)^\diamond\ar[d]\\
        R\M(E\times\Gmlog)^\diamond\ar[r,"\ev^{\rub}"] & \Ev^{\rub}(E\times\Gmlog)^\diamond
    \end{tikzcd}.
\end{equation}

To do so, we can for example proceed as follow: start with some logarithmic subdivision $\Ev^\diamond(E\times\PP^1)=E^{n+m}\times P_{\Sigma}$ in which the exotic insertion we want to consider lives,  and take the induced subdivision of $\Ev^{\rub}(E\times\Gmlog)$ obtained taking the image $\bar{\Sigma}\subseteq \RR^m/\RR$ of the fan  $\Sigma\subseteq\RR^m$. We then take $\Ev(E\times\Gmlog)^\diamond$ and $\Ev^\rub(E\times\Gmlog)^\diamond$ to be refinements of the subdivisions above needed to ensure that $\Ev(E\times\Gmlog)^\diamond\to \Ev^\rub(E\times\Gmlog)^\diamond$ is flat and that these two spaces are smooth.



\begin{proof}[Proof of Theorem \ref{theo:rigid-to-rubber}]
\begin{description}[leftmargin=0.6cm]
\item[Step I: rigid to rubber] The moduli spaces and the evaluation spaces considered above fit, as shown in \cite[Proposition~2.4.5]{ranganathan2024logarithmic}, into the following commutative diagram, parallel to the diagram appearing in \cite[Proof Thm~B,Step I]{ranganathan2024logarithmic}:

\begin{equation}
    \begin{tikzcd}
        \M(E\times\PP^1)^\diamond\ar[r,"b"]\ar[dr,"\mu "] & F\ar[r,"\ev"]\ar[d,"t"]\ar[rd, phantom, "\square"] & \Ev(E\times\Gmlog)^\diamond\ar[d,"\delta"]\\
       &  R\M(E\times\Gmlog)^\diamond\ar[r,"\ev^{\rub}"] & \Ev^\rub(E\times\Gmlog)^\diamond ,
    \end{tikzcd}
\end{equation}
where $F$ is completing the cartesian square, $t$ and $\delta$ are flat, $b$ is virtually birational (by Costello's pushforward \cite{costello2006higher}, see also \cite[Lemma~3.3.1]{ranganathan2024logarithmic}) and $t$ is a subdivision of
    \[\M(E\times\Gmlog)\xrightarrow{\mu} R\M(E\times\Gmlog).\]

    By \cite[Theorem~3.3.2]{ranganathan2024logarithmic}, $F$ is endowed with a natural virtual class satisfying 
    \[\vir{F}=t^*\vir{R\M(E\times\Gmlog)^\diamond}=b_*\vir{\M(E\times\PP^1)^\diamond}.\]

We conclude applying the push-pull formula for the cap product we obtain that
    \begin{align}\label{eq:rigidtorubber}
        \mu_*(\vir{\M(E\times\PP^1)^\diamond}\cap b^*\ev^*\gamma)&=t_* (b_* \vir{\M(E\times\PP^1)^\diamond}\cap \ev^*\gamma)\\
        &= t_*( t^*\vir{R\M(E\times\Gmlog)^\diamond} \cap \ev^*\gamma) \nonumber\\
        &= \vir{R\M(E\times\Gmlog)^\diamond}\cap \ev^{\rub,*}\delta_*\gamma \nonumber
    \end{align}
    where in the last equality we have used that the square is cartesian with proper and flat vertical maps.

\medskip

    Recall that $\gamma$ is assumed to be split as $\bar{\gamma}\otimes\gamma_t$. The push-forward $\delta_*$ only acts on the second factor. The cohomology class $\delta_*\gamma$ can be thus be written as $\bar{\gamma}\otimes \ell_0\in H^*(E^{n+m})\otimes H^*(P_{\bar{\Sigma}})$. Since the dimensions of the toric targets differ by $1$ and the push-forward prserves the homological rank, the (Chow) cohomological rank of $\ell_0$ is $1$ less as the rank of $\gamma_t$.



    \medskip

    We have thus managed to express the relative invariants in terms of rubber invariants capped with a suitable class pulled back from the rubber evaluation space.

\medskip

\item[Step II: rewriting the rubber insertion ]

  Recall that we have a subdivision
  $$\Ev^\rub(E\times\Gmlog)^\diamond= E^{n+m}\times P_{\bar{\Sigma}},$$
  for some smooth $m-1$-dimensional toric variety. In particular, we have the surjectivity of the pull-back map 
  $$H^*([P_{\bar{\Sigma}}/T_{\bar{\Sigma}}])\to H^*(P_{\bar{\Sigma}}),$$
  where $T_{\bar{\Sigma}}$ is the dense torus of $P_{\bar{\sigma}}$. This remains true when tensoring with $H^*(E^{n+m})$. We denote by $\gamma^\rub:=\bar{\gamma}\otimes \ell$  a lift of $\delta_*\gamma=\bar{\gamma}\otimes \ell_0.$ 

    \medskip 
    
As argued in \cite[Proof Thm~B, Step II]{ranganathan2024logarithmic} via the usual description for the equivariant cohomology of toric varieties (see for example \cite{payne2006equivariant}) $\ell$ can be interpreted as a piecewise polynomial function on the cone complex $\bar{\Sigma}.$ 

\medskip

Composing the rubber evaluation $\ev^{\rm{rub}}$ with the projection to the Artin fan we obtain an evaluation
    \[R\M(E\times\Gmlog)^\diamond\xrightarrow{\underline{\ev}} E^{n+m}\times[P_{\bar{\Sigma}}/T_{\bar{\Sigma}}]\]
    which factors through $R\M(E\times\Gmtrop)^\diamond$. As explained above, $R\M(E\times\Gmtrop)^\diamond$ is simply a tropical model of $\M(E)$, which we can write as fiber product (which is also the fs fiber product) with
    $$\mathbf{R}_{g,m,\bfw}\to\mathbf{Div}_{g,m,\bfw}\to\bar{\M}_{g,m+n}^{\log},$$
    for some logarithmic modification $\mathbf{R}_{g,m+\bfw}$.

    \medskip

    Since the morphism $\underline{\ev}$ is  logarithmic, it descends to a morphism between the Artin fans/tropicalizations. As the log structure on $R\M(E\times\Gmtrop)^\diamond$ is pulled-back from  $\mathbf{R}_{g,m,\bfw},$ this simply becomes a morphism

    \[\Sigma(\mathbf{R}_{g,m,\bfw})\xrightarrow{\mathrm{Trop}(\underline{\ev})}\bar{\Sigma}.\]

    In \cite{BCUK2026corrDDR} the cone complex $\Sigma(\mathbf{R}_{g,m,\bfw})$ is denoted by $\DR^{\mathrm{trop}}_{g,m}(\bfw).$

    \medskip

    Via $\mathrm{Trop}(\underline{\ev}),$ the piece-wise polynomial $\ell$ can be thought as a piece-wise polynomial on $\Sigma(\mathbf{R}_{g,m+\bfw}).$ But now we are exactly in the situation of \cite[Theorem~B,Step III]{ranganathan2024logarithmic} and we can choose an extension of this piece-wise linear polynomial to a piece-wise linear polynomial on $\bar{\M}_{g,n+m},$ which we denote by $\Phi_{\bfw}.$

Using the cone complex/logarithmic maps dictionary, this is the same as a map from some subdivision $\bar{\M}^\diamond_{g,n+m}$ to $[P_{\bar{\Sigma}}/T_{\bar{\Sigma}}].$ As explained again in \cite[Theorem~B,Step III]{ranganathan2024logarithmic}, the subdivision $\bar{\M}_{g,n}$ is such that 
$\Sigma(\mathbf{R}_{g,m+\bfw})\to\Sigma (\bar{\M}^\diamond_{g,n+m})$ is strict.

\medskip

Recapping, we have written the class \[\ev^{\rub,*}\delta_*\gamma =\epsilon^{\diamond,*}(\bar{\ev}^*\bar{\gamma}\otimes\mathrm{ft}^{\diamond,*}\Phi_{\bfw}),\]
where
\[R\M(E\times\Gmlog)^\diamond\xrightarrow{\epsilon^\diamond}R\M(E\times\Gmtrop)^\diamond\xrightarrow{\mathrm{ft}^{\diamond}} \bar{\M}^\diamond_{g,n+m},\]
and 
\[ R\M(E\times\Gmlog)^\diamond\to  R\M(E\times\Gmtrop)^\diamond\xrightarrow{\bar{\ev}^\diamond} E^{n+m}.\]

\item[Step III: projection formula]
Starting from \eqref{eq:rigidtorubber} and applying the projection formula along the proper morphism $\epsilon^\diamond$ and then along the forgetful morphism to 
$\bar{\M}^\diamond_{g,n+m}$ we conclude.
\end{description}
\end{proof}

\begin{coro}\label{cor:rigidtorubbercorrelated}
    Let $\gamma^\diamond = \bar{\gamma}\otimes\gamma_t\in H^*( E^{n+m})\otimes H^*(P_{\Sigma})$ be an exotic insertion. There exists $\Phi_{\bfw}$ a piece-wise polynomial on $\Sigma_{\bar{\M}_{g,n}}$ (corresponding to a cohomology class in some tropical model $\bar{\M}^\diamond$ ) such that 
    \[\ft_*(\fvir{\M(a,\bfw)^{\diamond}}{|\bfw|} \cap\ev^*\gamma^\diamond) = \ft_*(\mathbf{logDR}(a,\bfw)\cap ( \ev^*\bar{\gamma}\cup \Phi_{\bfw}) \in \mathrm{logCH}(\bar{\M}_{g,n+m}))\otimes\G_{\bfw}.\]
\end{coro}

\begin{proof}
    All the steps in the proof of Theorem~\ref{theo:rigid-to-rubber} are compatible with restriction to a fixed correlator $\theta$-component. See also \cite[Section~6]{BCUK2026corrDDR}
\end{proof}

\begin{rem}\label{rem:PLpolscaling}
 Notice that, as explained for example in the proof of \cite[Theorem~6.13]{BCUK2026corrDDR}, scaling the tangency order $\bfw$ by an integer $k$ we get that 

\[\Phi_{k\bfw}=k^{\deg\ell}\Phi_{\bfw}.\]

This comes from the explicit description of $\mathrm{Trop}(\underline{\ev})$ and the fact that there is an embedding of cone complexes $\Sigma(\mathbf{R}_{g,m,\bfw})\hookrightarrow \Sigma(\mathbf{R}_{g,m,k\bfw})$: one can choose as an extension of $\ell$ for $\bfw$ the PP-function $\Phi_{kw}$, but due to the $k$-scaling in the tropical picture, one gets a $k^{\deg\ell}$ factor doing so.       
\end{rem} 

\subsection{Multiple cover formula for $E\times\PP^1$}
\label{sec:MCF-EP1-proof}

We now prove the MCF for $E\times\PP^1$.

\begin{theo}
\label{theo-MCF-EP1-general}
    Let $\gamma^\diamond$ be an exotic insertion. Then, the correlated exotic GW-cycle $\gen{\gen{\gamma^\diamond}}_{a,\bfw}$ satisfies the $(2g-2+\deg\gamma^\diamond)$-MCF.
\end{theo}

\begin{proof}
    By Corollary~\ref{cor:rigidtorubbercorrelated}, we have
    \[\gen{\gen{\gamma^\diamond}}_{a,\bfw} = \ft_*(\mathbf{logDR}(a,\bfw)\cap (\ev^*\bar{\gamma}\cup \Phi_{\bfw}) ).\]
By Theorem \ref{theo:MCF-for-correlated-DR}, $\mathbf{logDR}(a,\bfw)\cap \ev^*\bar{\gamma}$ satisfies the $(2g-1+\deg\bar{\gamma})$-MCF. So

$$\mathbf{logDR}(a,\bfw)\cap (\ev^*\bar{\gamma}\cup \Phi_{\bfw}) ) =\sum_{k|a,\bfw} k^{2g-1+ \deg\bar{\gamma}}  \mathbf{logDR}(a/k,\bfw/k)\cap \ev^*\bar{\gamma}\cup \Phi_{\bfw}).$$
Next, we use Remark~\ref{rem:PLpolscaling} saying that for the $k$-factor, $\Phi_\bfw$ can be replaced by $k^{\deg\ell}\Phi_{\bfw/k}$:
$$\mathbf{logDR}(a,\bfw)\cap (\ev^*\bar{\gamma}\cup \Phi_{\bfw}) )=\sum_{k|a,\bfw} k^{2g-1+ \deg\bar{\gamma}+\deg\ell}  \mathbf{logDR}(a/k,\bfw/k)\cap \ev^*\bar{\gamma}\cup \Phi_{\bfw/k}).$$
Due to the rank dropping by $1$, we have actually $\deg\bar{\gamma}+\deg\ell=\deg \gamma^\diamond-1$, finishing the proof.
\end{proof}

Capping with a class $\alpha$ of complementary dimension, one gets the analog formula for correlated GW-invariants instead, and taking the $(\theta)$-coefficient, one gets the version of the formula free of the group algebre formalism.

\begin{coro}\label{coro:MCF-for-EtimesP}
    The correlated GW-invariants satisfy the $(3g-3+n+m-\deg\alpha)$-MCF:
    $$\gen{\alpha;\gamma}_{g,a,\bfw}^{|\bfw|,\theta} = \sum_{k|a,\bfw,\theta} k^{3g-3+n+m-\deg\alpha}\gen{\alpha;\gamma}_{g,a/k,\bfw/k}^{|\bfw/k|,\prim} ,$$
    where $k|\theta$ means that $\theta$ is of $(|\bfw|/k)$-torsion.
\end{coro}

We also treat separately the easy case where we have no interior marking, where the formula also applies.

\begin{lem}\label{lem:reducetogenpearls}
    Assume that we have no interior marking: $m=0$. Then,
    \begin{enumerate}
        \item We have $g=0$, forcing $a=0$, and $n=2$, i.e. the curve class in $E\times\PP^1$ is a multiple of $[\PP^1]$.
        \item Furthermore in this case, the invariant can be explicitly computed: there are only two possibility for the Poincar\'e insertions, either one is a point class and the other is trivial, or they are both $1$-cycles from $E$. In both cases we have
    \[\gen{\gen{\pt_w,1_{-w}}}^w_{0,0,(w,-w)} =\gen{\gen{\alpha_w,\beta_{-w}}}^w_{0,0,(w,-w)}  = \frac{1}{w}\cdot (0), \]
    which satisfies the $(-1)$-MCF.
    \end{enumerate}
\end{lem}

\begin{proof}
Unless $g=0$ and $n=2$, we can consider the natural map $\M(a,\bfw)\xrightarrow{\pi_1} R\M(a,\bfw)$ to the moduli space of rubber stable maps. Due to the absence of interior markings, all possible insertions are pulled back from $R\M(a,\bfw)$ since both boundary evaluations and the forgetful map to the moduli space of curves factor through the latter.

Since the map $\M(a,\bfw)\to R\M(a,\bfw)$ has one dimensional fibers (and  also the rank of the virtual classes differ by one), we have that
\[\pi_{1,*}\vir{\M(a,\bfw)}=0.\]

Therefore, only remains the case where $g=0$, forcing $a=0$, and $n=2$. The explicit computation can be found for example in \cite[Proposition~3.17]{blomme2024bielliptic}. There is no correlation in this case since a fiber is mapped to a point by the projection.

The $(-1)$-MCF is obvious from the expression, though it also follows from the general MCF.
\end{proof}
\section{Multiple cover formula for abelian surfaces}\label{sec:degenerationformula}

We conclude  with the proof of Theorem \ref{theo-MCF-abelian-general}, the general case of the MCF for abelian surfaces. The proof relies on Theorem \ref{theo-MCF-EP1-general} and the reduced decomposition formula. 
To handle vanishing cohomology in the degeneration formula, we start reviewing the results of \cite[Section 2]{maulik2025gromov} and carry out the straightforward adaptation to the reduced virtual class.

\subsection{Vanishing cohomology in the reduced degeneration formula}\label{sec:vanishingcohomology}

We recall two foundational results due to Deligne \cite{deligne1971theorie} about the singular cohomology of families of algebraic varieties:

\begin{theo}\cite[Theorem~4.2.6, Theorem~4.1.1]{deligne1971theorie} 
    \begin{enumerate}[leftmargin=0.5cm]
        \item Let $X\to T^\circ=T\setminus\left\{t_1,\dots,t_k\right\}$ a smooth projective morphism of algebraic varieties with $T$ a smooth projective connected curve, and $t\in T^\circ$. Then, the monodromy action of $\pi_1(T^\circ)$ on the cohomology $H^*(X_t)$ is semisimple.  In particular, we have a projection operator
        \[H^*(X_t)\to H^*(X_t)^{\pi_1(T^\circ)}.\]
        \item Let $X\to T^\circ$ be as before, and $X\hookrightarrow\bar{X}$ an open immersion into a smooth variety, projective over a proper base $T$. Then 
        \[\mathrm{Im}(H^*(\bar{X})\to H^*(X_t))= H^*(X_t)^{\pi_1(T^\circ)}.\]
    \end{enumerate}
\end{theo}

Let $A\to T^\circ$ be a smooth family of abelian surfaces with polarization $\beta$. We denote by $\M(A_t)=\M_{g,m}(A_t,\beta)$ the moduli space of stable maps, forgetting about the unnecessary information. Point $(2)$ states that the restriction map to the cohomology of a fiber is often not surjective. However, from a GW-theory point of view, deformation invariance ensures that monodromy-invariant classes are enough to recover the full theory. Using point $(1)$ above together with the invariance by algebraic deformations of $A$ of the reduced virtual class $\red{\M(A_t)}$, we get the following Proposition.

\begin{prop}\cite[Proposition~2.2.2]{maulik2025gromov}
    Let $\gamma\in H^*(A^m_t)$ a cohomology class and $\bar{\gamma}\in H^*(A^m_t)^{\pi_1(T^\circ)}$ its monodromy invariant projection. We then have the following equality of GW-cycles:
    \[\ft_*(\red{\M(A_t)}\cap\ev^*\gamma)=\ft_*(\red{\M(A_t)}\cap\ev^*\bar{\gamma})\in H_*(\overline{\M}_{g,m},\mathbb Q).\]   
\end{prop}

Let $A\to T$ be a semi-stable degeneration of abelian varieties over a smooth proper curve, as constructed in Section~\ref{sec-construction-abelian-families}. We consider the $m$th-power $A_T^m$. The problem is that even though the total space of $A$ is smooth, the $m$-fold fiber product $A_T^m$ acquires toroidal singularities. Therefore, one cannot apply directly Deligne's Theorem: it is not true that any monodromy invariant class $\bar{\gamma}\in H^*(A^m_t)^{\pi_1(T^\circ)}$ extends to a class in $H^*(A^m_T)$.

\medskip

To remedy this problem, one has to modify the total space $A_T^m$, and this leads to the appearance of \textit{exotic insertions}. As proved in \cite[Lemma~2.3.1]{maulik2025gromov}, there exists a tropical model $(A_T^m)^\diamond\to T$ such that: the total space is smooth, the map to $T$ is semi-stable and  $(A_T^m)^\diamond$ has no new components. In other words, the map between tropicalizations $\Sigma( (A^m_T)^\diamond)\to\Sigma(A^m_T)$ is a subdivision, but does not contain any new ray.

\medskip

By Deligne's Theorem, every monodromy invariant class $\bar{\gamma}\in H^*(A^m_t)^{\pi_1(T^\circ)}$ extends to a cohomology class on a resolution of the family of evaluation spaces $(A_T^m)^\diamond$. We consider the fine and saturated (not cartesian) fiber product where $\M(A_T)=\M_{g,m}(A_T,\beta)$:

\begin{equation}
\begin{tikzcd}
\M(A_T)^{\diamond}\ar[r]\ar[d] &\M(A_T)\ar[d]\\
(A_T^m)^\diamond \ar[r] & A_T^m.
\end{tikzcd}
\end{equation}

Then $\M(A_T)^\diamond$ is a tropical model for the space of maps. Recall that the map to Olsson's stack $\M(A_T)\to\mathcal Log_T$, and thus the map to an Artin fan $\mathcal A_{\mathcal M}\xrightarrow{p}\mathcal Log_T$ (with $p$ representable and strict \'etale) for the moduli space, factors through the log smooth moduli space $\mathfrak {M}$. In particular, we have a tropical model (still log smooth by definition) $\mathfrak {M}^{\diamond}$ and the relative obstruction theory $\mathbf{E}_{\mathrm{rel}}$ and $ \mathbf{E}^{\mathrm{red}}_{\mathrm{rel}}$ are automatically obstruction theories for 
$\M(A_T)^{\diamond}\to \mathfrak {M}^{\diamond}.$ 

\medskip

Let $A_0$ be one of the singular fibers of the family. By restriction, we obtain $\M(A_0)^\diamond$ a tropical model for $\M(A_0)=\M_{g,m}(A_0,\beta)$ equipped with an evaluation map to the fiber $(A_0^m)^\diamond$.  Given $\gamma_0$ the specialization of a cohomology class $\gamma\in H^*((A_T^m)^\diamond)$, we now seek a decomposition formula adapted to compute invariants with insertion $\gamma_0$.

\subsubsection{Virtual decomposition}

Proposition~\ref{prop:splitting-diagrams} applies also for the tropical modification $\M(A_0)^\diamond$ (see \cite[Section~3]{abramovich2020decomposition}), with the only difference that the \emph{rigid tropical curves} in this case is in the sense of \cite[Definition~3.6]{abramovich2020decomposition}, i.e. these are tropical maps of fixed type whose associated cone is a ray in the cone stack $\Sigma(\mathfrak M^\diamond).$
Notice that since $\Sigma(\mathfrak M^\diamond)\to\Sigma(\mathfrak M)$ is a subdivision, there could be tropical types which are not rigid in $\mathfrak M$ but index rays of the chosen model (see also  \cite[Section~7]{MaulikRanganathan2025}). We get
$$\red{\M(A_0)^\diamond}=\sum_{\Gamma} l_{\Gamma}\cdot \mu_{\Gamma,*} \red{\M^\diamond_{\Gamma}}.$$
    where the sum is over rigid tropical curves indexing rays of $\Sigma(\mathfrak M^\diamond).$

\subsubsection{Splitting across vertices}

Let $\phi\colon\Gamma\to\Sigma_0$ be one of the (new) rigid tropical curves. The evaluation map for the interior markings
 $\M_{\Gamma}\to A_0^m$ factors as
 $$\M_{\Gamma}\to\prod_{v} Y_{\phi(v)}^{m_v}\to A_0^m,$$
 where $\prod_{v} Y_{\phi(v)}^{m_v}$ is the irreducible components of $A_0^m$ determined by $\Gamma.$ This induces a factorization
 \[\M^\diamond_{\Gamma}\to \left(\prod_{v} Y_{\phi(v)}^{m_v}\right)^\diamond\to (A_0^m)^\diamond,\]
where the central term is determined by the blow-up $(A_0^m)^\diamond$. Then, one gets a tropical model $\left(\prod_v \M_v\right)^\diamond$ of $\prod\M_v$ induced by fs pull-back:
 \begin{equation}
\begin{tikzcd}
\left(\prod_v \M_v\right)^{\diamond}\ar[r]\ar[d] & \prod_v \M_v\ar[d]\\
\left(\prod_{v} Y_{\phi(v)}^{m_v}\right)^\diamond \times\prod_v E^{n_v} \ar[r] & \prod_{v} Y_{\phi(v)}^{m_v}\times\prod_v E^{n_v}.
\end{tikzcd}
\end{equation}

It follows immediately by Proposition~\ref{prop:splitting-vertices} and Lemma~\ref{lem:virtualmodels} an equality of cycles in $A_*(\M^{\diamond}_{\widetilde{\Gamma}}):$
    \[\red{\M^\diamond_{\widetilde{\Gamma}}}=\varphi_{\Gamma}^*\Delta_u^{\mathrm{red},!}\left( \vir{(\prod_v \M_v)^\diamond}\right).\]

Therefore, we know that the entire GW-theory of $A_t$ can be reconstructed by invariants of $\vir{(\prod_v\M_v)^\diamond}$ with insertions coming from $\left(\prod_{v} Y_{\phi(v)}^{m_v}\right)^\diamond \times\prod_v E^{n_v}$. Notice that by Lemma~\ref{lem:virtualmodels}, we are free to pass to even finer subdivisions 
    \[\left(\prod_v \M_v\right)^{\spadesuit}\to \left(\prod_{v} Y_{\phi(v)}^{m_v}\right)^\spadesuit\times\prod_v E^{n_v}.\]

\subsubsection{Product structure}

The disadvantage of the above formulation is that we do not (yet) have an actual splitting as a product over the vertices. The situation can be fixed as in \cite[Theorem~2.4.2]{maulik2025gromov}. Actually, the same technique is also used in \cite[Section~12]{MaulikRanganathan2025} to obtain a strong splitting formula for simple normal crossing degenerations.

\medskip

 By \cite[Proposition~1.7.2]{MaulikRanganathan2025} we can choose tropical models $\M_v^\dagger$ and $(Y_{\phi(v)}^{m_v})^\dagger\times E^{n_v}$ such that the induced subdivided evaluation map 
 \[\prod _v \M_v^\dagger\to  \prod_v (Y_{\phi(v)}^{m_v})^\dagger\times E^{n_v},\]
 is \emph{combinatorially flat and with reduced fibers}, as in \cite[Definition~1.7.1]{MaulikRanganathan2025}. Up to further blow-up, one can consider that we have a map from a blow-up $\left(\prod_{v} Y_{\phi(v)}^{m_v}\right)^\spadesuit$ of the previously constructed $\left(\prod_{v} Y_{\phi(v)}^{m_v}\right)^\diamond$. We can consider the fine and saturated fiber product
\begin{equation}
\begin{tikzcd}
\left(\prod_v \M_v\right)^{\spadesuit}\ar[r]\ar[d] & \prod_v \M_v^\dagger\ar[d]\\
\left(\prod_{v} Y_{\phi(v)}^{m_v}\right)^\spadesuit\times\prod_v E^{n_v} \ar[r] & \prod_{v} (Y_{\phi(v)}^{m_v})^\dagger\times\prod_v E^{n_v}.
\end{tikzcd}
\end{equation}
The fact that the right vertical map is combinatorially flat with reduced fibers ensures that the underlying diagram of Deligne-Mumford stacks is cartesian (see \cite[Proposition~1.7.3]{MaulikRanganathan2025}). 

\medskip

Then, for any insertion $\gamma^\spadesuit$ in the cohomology of the evaluation space on the bottom left, denoting by $\gamma^\dagger$ its image via proper push-forward along the bottom map, via a push-pull-argument we get that

\[\ft_*\left (\vir{(\prod_v \M_v)^{\spadesuit}}\cap \ev^*(\gamma^\spadesuit\otimes\Delta_u^{\mathrm{red} })\right)=\ft_*\left(\prod_v \vir{\M^{\dagger}_v}\cap\ev^*(\gamma^\dagger\otimes \Delta_u^{\mathrm{red}})\right).\]

In conclusion, any GW-invariant of $A_t$ can be reconstructed from the exotic invariants of $E\times\PP^1$ via reduced degeneration formula.

\begin{rem}
    Notice that since there are no rational curves on an abelian surface, the $\psi_i$ classes are pulled-back from $\overline{\M}_{g,m}$ so the descendant theory just comes from geometric insertions and insertions pulled back from the moduli of curves.
\end{rem}

\subsection{Multiple cover formula for abelian surfaces}


Let $A$ be an abelian surface with polarization $\beta$, and its moduli space of stable maps $\M(A,\beta)=\M_{g,m}(A,\beta)$. Let $\gamma\in H^*(A^m,\QQ)$ be a cohomological insertion. We consider the following (reduced) GW-cycles:
$$\gen{\gamma}_{g,\beta} = \ft_*(\red{\M(A,\beta)}\cap\ev^*\gamma).$$
When the setting is clear, we also delete $g$ from the notation. To get invariants, one just has to intersect with a cohomology class from $\bar{\M}_{g,m}$:
$$\gen{\alpha;\gamma}_{g,\beta} = \int_{\red{\M(A,\beta)}}\ft^*\alpha\cup\ev^*\gamma.$$

For each $k|\beta$, choose $\varphi_k$ be an isometry of $H^*(A,\QQ)$ mapping $\beta/k$ to an integral primitive class, thus sharing the same self-intersection.

\begin{theo}\label{theo-MCF-abelian-general}
    The reduced GW-cycles satisfy the MCF
    $$\gen{\gamma}_{g,\beta} = \sum_{k|\beta} k^{2g-3+\deg\gamma}
    \gen{\varphi_k(\gamma)}_{g,\varphi_k(\beta/k)}.$$
    Similarly, reduced GW-invariants satisfy the MCF:
    $$\gen{\alpha;\gamma}_{g,\beta} = \sum_{k|\beta} k^{3g-3+n-\deg\alpha}
    \gen{\alpha;\varphi_k(\gamma)}_{g,\varphi_k(\beta/k)}.$$
\end{theo}


The proof follows the same steps as in the case of point insertions:
\begin{itemize}
    \item Reformulate the above MCF in the group algebra language for a group algebra valued function $\BBB\to\G\otimes H^*(\bar{\M}_{g,m},\QQ)$, working with the families of elliptically fibered abelian surface from Section~\ref{sec-construction-abelian-families}.
    \item Use the reduced degeneration formula to write the latter as the push-forward of a (correlated) multiplicity function $\bfm$ defined on the $\NN$-module $\TTT$ of (rigid) tropical curves.
    \item Using Lemma \ref{lem:product-rule} and Lemma \ref{lem:restriction-suborbit}, we show that the multiplicity function satisfies the MCF by restricting to the orbit of a single tropical curve, and conclude by Proposition \ref{prop:push-forward} (push-forward of a function).
\end{itemize}

\medskip

The group algebra valued function is defined as follows. Let $B=(|B|,a)$ be a diagram degree. We have the families of elliptically fibered abelian surfaces $A_T(u)$, with $u\in E[|B|]$ the possible monodromies. For $u=0$, one gets the product $E\times F$. The diagram degree defines the polarization $\beta(0)=|B|[F]+a[E]$ on $E\times F$, which we also denote by $B$. For each other $u$, we provided an isomorphism $\varphi_u\colon H^*(A(0))\to H^*(A(u))$ such that the polarization on $A(u)$ is $\beta(u)=\varphi_u(\beta(0))=\varphi_u(B)$.

Let $\gamma\in H^*(A(0)^m,\QQ)$ be an insertion. Forgetting about the genus $g$ in the notation, we set
$$\gen{\gen{\gamma}}_B = \sum_{u\in E[|B|]} \gen{\varphi_u(\gamma)}_{\beta(u)}\cdot (u) \; \in \; \G\otimes H^*(\bar{\M}_{g,m},\QQ),$$
which is a group algebra element encoding the GW-cycles for all elliptically fibered abelian surfaces constructed at the same time. If $\alpha\in H^*(\bar{\M}_{g,m},\QQ)$, we also have the associated invariants:$$\gen{\gen{\alpha;\gamma}}_B = \sum_{u\in E[|B|]} \gen{\alpha;\varphi_u(\gamma)}_{\beta(u)}\cdot (u).$$

\medskip

We denote by $\TTT$ the set of tropical curves $\phi\colon\Gamma\to\Sigma_0$. It is endowed with a $\NN$-action by scaling edge weights $w_e$ (and lengths) and vertex labels $a_v$. It becomes a $\NN$-module with the norm pulled back from the degree morphism $\deg\colon\TTT\to\BBB$ mapping $\phi\colon\Gamma\to\Sigma_0$ to $(\deg\phi,\sum a_v)$, where $\deg\phi$ is the degree of the harmonic map $\phi$. In particular, $\deg$ becomes a morphism of $\NN$-modules.

\medskip

Given a tropical curve $\phi\colon\Gamma\to\Sigma_0$, we denote by $\bfw$ (resp. $\bfw_v$) the collection of edge weights (resp. adjacent to $v$). We set $\delta_\Gamma=\gcd(\bfw)$ the g.c.d. of edge weights and $\delta_v$ the g.c.d. of edge weights adjacent to $v$. Notice $\delta_\Gamma\mid |\Gamma|$ but they usually differ since $\deg\phi$ is a sum of some edge weights. For instance, the degree could be non-trivial and every edge weight being equal to $1$. Similarly, $\ell(\Gamma)|\delta_\Gamma$ but they usually differ since $\ell(\Gamma)=\gcd(\delta_\Gamma,(a_v))$.

\begin{proof}[Proof of Theorem \ref{theo-MCF-abelian-general}]
\textbf{$\bullet$ Step 0: group algebra formulation.} We start showing that the MCF for abelian surfaces is equivalent to the MCF for $B\mapsto\gen{\gen{\gamma}}_B$. This follows from a computation fairly similar to the $E\times\PP^1$ case taking the $(u)$-coefficient. From the expression of $\gen{\gen{\gamma}}_B$, one gets
$$\Prim_{|B|}\gen{\gen{\gamma}}_B = |B|^2\cdot \gen{\varphi_{u_\prim}(\gamma)}_{\beta(u_\prim)},$$
where $u_\prim$ is a monodromy of order $|B|$. Therefore, assuming the group algebra MCF, and denoting by $u_{k,\prim}$ an order $|B|/k$ monodromy, we get
\begin{align*}
        \gen{\gen{\gamma}}_B = &\sum_{k|B} k^{2g-3+\deg\gamma} \Prim_{|B/k|}\gen{\gen{\gamma}}_{B/k}\cdot T_{|B|/k}\\
         = & \sum_{k|B} k^{2g-3+\deg\gamma}\gen{\varphi_{u_{k,\prim}}(\gamma)}_{\varphi_{u_{k,\prim}}(B/k)} |B/k|^2 T_{|B|/k}\\
         =&\sum_{k|B} k^{2g-3+\deg\gamma}\gen{\varphi_{u_{k,\prim}}(\gamma)}_{\varphi_{u_{k,\prim}}(B/k)}\sum_{|B/k|u=0}(u)\\
         \sum_{|B|u\equiv 0} \gen{\varphi_u(\gamma)}_{\varphi_u(B)}\cdot (u) = & \sum_{|B|u\equiv 0} \sum_{k|B,u}k^{2g-3+\deg\gamma}\gen{\varphi_{u_{k,\prim}}(\gamma)}_{\varphi_{u_{k,\prim}}(B/k)}(u), \\
          \end{align*}
where we used that $k|u$ is equivalent to $u$ being of $(|B|/k)$-torsion. We deduce the MCF for abelian surfaces taking each $(u)$-coefficient. Conversely, it suffices to sum over the possible monodromies. We get all cases since for each self-intersection, it is possible to choose a diagram degree $B$ with maximal divisibility.

\medskip

\textbf{$\bullet$ Step 1: Applying the decomposition formula.} By Corollary \ref{coro:degen-gp-alg}, one has
\[ [\![\M_0]\!]^\mathrm{red}=\sum_{\widetilde{\Gamma}}\frac{l_{\Gamma}}{|\sfE(\Gamma)|!}\mu_{\widetilde{\Gamma},*}\varphi_{\Gamma}^*\Delta^{\mathrm{red},!}\left( (|B|/\delta_\Gamma)^2\d{\frac{1}{|B|/\delta_\Gamma}}\left(\prod_v\m{\delta_v/\delta_\Gamma}\fvir{\M_v}{\delta_v}\right) \right).\]
Let $\gamma$ be class insertion allowed by the usual decomposition formula, i.e. coming from the total space of the degeneration. Capping with $\ev^*\gamma$ and pushing forward to $\bar{\M}_{g,m}$ using that $\gen{\gen{\gamma}}_B = \ft_*([\![\M_0]\!]^\mathrm{red}\cap\ev^*\gamma)$ and $\mu_{\widetilde{\Gamma}}$ has degree $\frac{\prod w_e}{l_\Gamma}$, one gets
\begin{align*}
\gen{\gen{\gamma}}_B = & \sum_{\widetilde{\Gamma}} \frac{\prod_e w_e}{|\sfE(\Gamma)|!} (|B|/\delta_\Gamma)^2\d{\frac{1}{|B|/\delta_\Gamma}}\ft_*\left(\prod_v\m{\delta_v/\delta_\Gamma}\fvir{\M_v}{\delta_v}\right)\cap\ev^*(\gamma\cup\Delta^\mathrm{red}_{e_0}) \\
= & \sum_{\Gamma} \frac{\prod_e w_e}{|\operatorname{Aut}(\Gamma)|} (|B|/\delta_\Gamma)^2\d{\frac{1}{|B|/\delta_\Gamma}}\ft_*\left(\prod_v\m{\delta_v/\delta_\Gamma}\fvir{\M_v}{\delta_v}\right)\cap\ev^*(\gamma\cup\Delta^\mathrm{red}_{e_0}). 
\end{align*}
This associates to each rigid tropical curve a class multiplicity $\bfm(\Gamma)$.

In case $\gamma$ does not come from the total space of the degeneration, Section \ref{sec:vanishingcohomology} states that the equality remains valid up to adding some blow-ups and consider an exotic insertion $\gamma^\diamond$. The remainder of the proof is just a notational game adding some $\diamond$ which we thus avoid to lighten the already heavy notation.

To express the multiplicity in a nicer way, we use the reduced K\"unneth formula for $\Delta^\mathrm{red}_{e_0}$:
$$\Delta_{e_0}^\mathrm{red}=\frac{\delta_\Gamma^2}{w_{e_0}^2}\sum_{P}\bigotimes_v\Delta_v(P),$$
where the sum is over the possible ``Poincar\'e insertions'' $P$ at every half-edge, inherited from the usual K\"unneth formula, and $\Delta_v(P)$ are the insertions adjacent to $v$. The coefficient $\frac{\delta_\Gamma^2}{w_{e_0}^2}$ does not really matter since it is invariant when scaling the tropical curve. The multiplicity of a tropical curve becomes
$$\bfm(\Gamma) = \sum_P \frac{\prod_e w_e}{|\operatorname{Aut}(\Gamma)|} \frac{|B|^2}{w_{e_0}^2}\d{\frac{1}{|B|/\delta_\Gamma}}\prod_v\m{\delta_v/\delta_\Gamma}\operatorname{Mult}_{\gamma_v,P}(a_v,\bfw_v),$$
$$\text{with}\; \operatorname{Mult}_{\gamma_v,P}(a_v,\bfw_v)=\ft_*\left( \fvir{\M_v}{\delta_v}\cap\ev^*(\gamma_v\cup\Delta_v(P)) \right),$$
and where $\gamma_v$ is the part of $\gamma$ associated to the markings adjacent to $v$, and the sum is over the possible choices of Poincar\'e insertions. We set $\bfm_P(\Gamma)$ to be the summand of the above sum. To conclude, one has to prove the MCF for the multiplicity function, which can be done restricting to the orbit of a single primitive tropical curve and a fixed choice of Poincar\'e insertion.

\medskip

\textbf{$\bullet$ Step 2: Proving the MCF for the multiplicity function.} Let $\phi\colon\Gamma\to\Sigma_0$ be a primitive tropical curve of degree $B$. Its norm $|\Gamma|=|B|$ may not be $1$. Let $\delta_\Gamma$ be the g.c.d. of edge weights, which may not be $1$: $\Gamma$ being primitive only means that $\gcd(\delta_\Gamma,(a_v)_v)=1$. Let $P$ be a choice of Poincar\'e insertion for the diagonal. We wish to prove the MCF for $\delta\mapsto\bfm_P(\delta\Gamma)$.

To do so, by Lemma \ref{lem:MCF-for-N-k}, one has to prove that $\delta\mapsto\bfm_P(\delta\Gamma)$ is invariant by $|B|$-torsion and check the MCF for $\delta\mapsto\m{|B|}\bfm(\delta\Gamma)$.

For the torsion-invariance, we use the invariance for each vertex. The vertex multiplicity is a correlated invariant for the ramification profile $(a_v,\bfw_v)$, of g.c.d. $\gcd(a_v,\delta_v)$ and norm $\delta_v=\gcd(\bfw_v)$
$$\operatorname{Mult}_{\gamma_v,P}(a_v,\bfw_v)=\ft_*\left( \fvir{\M_v}{\delta_v}\cap\ev^*(\gamma_v\cup\Delta_v(P)) \right).$$
It is therefore invariant invariant by multiplication by $T_{\delta_v/\gcd(a_v,\delta_v)}$. We have the following equality:
$$\frac{\delta_v}{\gcd(a_v,\delta_v)\gcd(\frac{\delta_v}{\gcd(a_v,\delta_v)},\frac{\delta_v}{\delta_\Gamma})} = \frac{\delta_\Gamma}{\gcd(a_v,\delta_\Gamma)}.$$
Hence, we deduce that $\m{\delta_v/\delta_\Gamma}(\operatorname{Mult}_{\gamma_v,P}(a_v,\bfw_v))$ is invariant by multiplication by $T_{\delta/\gcd(a_v,\delta_\Gamma)}$. Therefore, their product is invariant by multiplication by
$$\prod_v T_{\delta_\Gamma/\gcd(\delta_\Gamma,a_v)} = T_{\mathrm{lcm}(\frac{\delta_\Gamma}{\gcd(a_v,\delta_\Gamma)})_v} = T_{\delta_\Gamma},$$
since $\gcd(\delta_\Gamma,(a_v)_v)=1$. Composing with $\d{\frac{1}{|\Gamma|/\delta_\Gamma}}$, we finally deduce the multiplicity is invariant by $|\Gamma|$-torsion.

We now have to prove the MCF after composition via $\m{|\Gamma|}$. Actually, we have the following composition:
$$\m{|\Gamma|}\d{\frac{\delta_\Gamma}{|\Gamma|}}\m{\frac{\delta_v}{\delta_\Gamma}} = \m{\delta_v}.$$
Therefore,
$$\m{|\Gamma|}(\bfm_P(\delta\Gamma)) = \frac{\prod_e (\delta w_e)}{|\operatorname{Aut}(\Gamma)|}\frac{|\delta B|^2}{(\delta w_{e_0})^2}\prod_v \m{\delta_v}\operatorname{Mult}_{\gamma_v,P}(\delta a_v,\delta\bfw_v).$$
The function $(a,\bfw)\mapsto\operatorname{Mult}_{\gamma_v,P}(a,\bfw)$ satisfies the MCF by Theorem \ref{theo-MCF-EP1-general}. The element $(a_v,\bfw_v)$ has norm $\delta_v$. By Lemma \ref{lem:restriction-suborbit}, $\delta\mapsto T_{\delta_v}\operatorname{Mult}_{\gamma_v,P}(\delta a_v,\delta\bfw_v)$ thus satisfies the MCF with norm $\delta_v$. Applying $\m{\delta_v}$, which cancels the $T_{\delta_v}$, we deduce that the factor of the product satisfies the MCF. Therefore, $\m{|\Gamma|}(\bfm_P(\delta\Gamma))$ is the product of a monomial by functions satisfying the MCF, and thus alaso satisfies the MCF.

\medskip

\textbf{$\bullet$ Step 3: exponent computation.} To conclude, we merely have to compute the exponent. Let $\phi\colon\Gamma\to\Sigma_0$ be a tropical curve with set of bounded edges $E_b$.
The exponent of the MCF is $|E_b|$ plus the sum of the exponents of the MCF for each vertex, provided by Theorem \ref{theo-MCF-EP1-general}:
$$|E_b|+\sum_v (2g_v-2+\deg\gamma_v) + |E_b|-1,$$
where the last term is due to the diagonal. By computing the Euler characteristic of $\Gamma$, we have $|V|-|E_b|=1-b_1(\Dfk)$. Therefore, the exponent is equal to
$$2\sum g_v+2(|E_b|-|V|)+\deg\gamma-1 = 2g-3+\deg\gamma.$$
\end{proof}

\begin{rem}
      Similarly to the MCF for abelian surfaces being compatible with the reduced decomposition formula, the correlated decomposition formula from \cite{blomme2024correlated} applied for $E\times\PP^1$ is also compatible with the MCF in the same way: it is true diagram by diagram and the MCF comes from a vertex level. However, there is little gain to proceed like that because the vertex contributions still are correlated GW-invariants of $E\times\PP^1$. This is the reason why the proof of Theorem \ref{theo-MCF-EP1-general} goes through the correlated DR-cycle
\end{rem}

\bibliographystyle{alpha}
\bibliography{biblio}

\end{document}